\documentclass[reqno]{amsart}

\usepackage{eucal}
\usepackage{latexsym}
\usepackage{amsmath}
\usepackage{amssymb}
\usepackage{amscd}
\usepackage{multicol}
\usepackage{amsfonts}
\usepackage{amsthm}
\usepackage{color}
\usepackage{graphicx}
\usepackage[all]{xy}
\usepackage{epsfig}
\usepackage{psfrag}
\usepackage{upgreek}
\usepackage{courier}

\input xy
\xyoption{all}

\theoremstyle{plain}
\newtheorem{thm}{Theorem}
\newtheorem{theorem}[thm]{Theorem}
\newtheorem{corollary}[thm]{Corollary}
\newtheorem{lemma}[thm]{Lemma}
\newtheorem{prop}[thm]{Proposition}
\newtheorem{defin}[thm]{Definition}

\newtheoremstyle{exm}
{9pt}{9pt}{}{}{\bfseries}{}{.5em}{}
\theoremstyle{exm}
\newtheorem{exm}[thm]{Example}
\newtheoremstyle{rmk}
{9pt}{9pt}{}{}{\bfseries}{}{.5em}{}
\theoremstyle{rmk}
\newtheorem{rmk}[thm]{Remark}
\theoremstyle{alg}
\newtheorem{alg}{Algorithm}
\newtheoremstyle{question}
{9pt}{9pt}{}{}{\bfseries}{}{.5em}{}
\theoremstyle{question}
\newtheorem{question}[thm]{Question}
\numberwithin{equation}{section}
\numberwithin{thm}{section}
\numberwithin{figure}{section}

\newcommand{\w}{\mathsf{w}}
\newcommand{\m}{\mathsf{m}}
\newcommand{\ii}{\mathsf{i}}
\newcommand{\te}{\mathsf{t}}
\newcommand{\M}{\mathsf{M}}
\newcommand{\N}{\mathsf{N}}
\newcommand{\lu}{\mathsf{l}}

\newcommand{\W}{\mathsf{W}}

\newcommand{\n}{\mathsf{n}}
\newcommand{\cc}{\mathsf{c}}
\newcommand{\ce}{\mathsf{c}_{1}^{S^1}}
\newcommand{\taut}{\widetilde{\tau}}
\newcommand{\LL}{\mathbb{L}}
\newcommand{\J}{\mathsf{J}}

\newcommand{\R}{\mathbb{R}}
\newcommand{\Z}{{\mathbb{Z}}}
\newcommand{\C}{{\mathbb{C}}}

\newcommand{\Q}{\mathbb{Q}}

\newcommand{\LS}{\mathbb{L}^{S^1}}
\DeclareMathOperator{\td}{Todd}
\DeclareMathOperator{\ttot}{\mathcal{T}}
\DeclareMathOperator{\sign}{sign}
\DeclareMathOperator{\ind}{ind}
\DeclareMathOperator{\ch}{Ch}
\DeclareMathOperator{\ns}{Null} 
\DeclareMathOperator{\diag}{diag}

\title[New tools]{New tools for classifying Hamiltonian circle actions with isolated fixed points}

\author[L. Godinho]{Leonor Godinho}
\address{Departamento de Matem\'{a}tica, Centro de An\'{a}lise Matem\'{a}tica, Geometria e Sistemas Din\^{a}micos-LARSYS, Instituto Superior T\'ecnico, Av. Rovisco Pais 1049-001 Lisbon, Portugal}
\email{lgodin@math.ist.utl.pt}

\author[S. Sabatini]{Silvia Sabatini}
\address{Department of Mathematics, EPFL, Lausanne, Switzerland}
\email{silvia.sabatini@epfl.ch}

\thanks{Partially supported by Funda\c{c}\~{a}o para a Ci\^{e}ncia e Tecnologia (FCT/Portugal)}

\date{June 14th, 2012}

\begin{document}

\begin{abstract}
For every compact almost complex manifold $(\M,\J)$
equipped with a $\J$-preserving circle action with isolated
fixed points, a simple algebraic identity involving the first Chern class is derived.
This enables us to construct an algorithm to obtain linear relations among the isotropy
weights at the fixed points.
Suppose that $\M$ is symplectic and the action is Hamiltonian.
If the manifold satisfies an extra ``positivity condition" this algorithm
determines 
a family of vector spaces which contain the admissible lattices of weights. 
 
When the number of fixed points is minimal, 
this positivity condition is necessarily satisfied whenever $\dim(\M)\leq 6$,
and, when $\dim(\M)=8$, whenever the $S^1$-action extends to an effective
Hamiltonian $T^2$-action, or none of the isotropy weights is $1$.
Moreover there are no known examples with a minimal number of fixed points
contradicting this condition, and their existence is related to interesting
questions regarding fake projective spaces \cite{Y}. 
We run the algorithm for $\dim(\M)\leq 8$,
quickly obtaining all the possible families of isotropy weights.
In particular, we simplify the proofs of Ahara and Tolman  for $\dim(\M)=6$ \cite{Ah,T1}
and, when $\dim(\M)=8$, we prove that the equivariant cohomology ring, Chern classes and
isotropy weights agree with the ones of $\C P^4$ with the standard $S^1$-action
(thus proving the symplectic Petrie conjecture \cite{T1} in this setting). 
\end{abstract}
\maketitle

\tableofcontents

\section{Introduction}
Given a compact smooth manifold $\M$ and a compact Lie group $G$
~\begin{footnote}
{In this work we will always assume $\M$ to be \emph{connected}, and the $G$-action to be \emph{smooth} and \emph{effective},
i.e. $\cap_{x\in \M}G_x=\{e\}$, where $G_x=\{g\in G\mid g\cdot x=x\}$ is the stabilizer subgroup of $x\in \M$.}
\end{footnote}, it is a hard problem to determine 
whether $\M$ admits a $G$-action, and if it does, how many actions it can have.
Petrie \cite{P1,P2} first addressed these questions when $G$ is a torus and $\M$ is a homotopy complex
projective space, meaning that it is homotopically equivalent to $\C P^n$.
When the fixed points  are isolated points, he showed that
it is crucial to understand the torus 
representations on the normal bundle to these points. Indeed, he proved that the Pontrjagin classes are determined 
by these representations and, in particular, when $G=T^n$, these 
classes agree with the ones of $\C P^n$. 
Motivated by this, he 
conjectured that the same conclusion would hold
if $G$ were simply $S^1$.
This conjecture has been proved in many particular situations, but the complete 
proof is still missing.

Suppose now that $(\M,\J)$ is a compact almost complex manifold of dimension $2n$,
and that $S^1$
acts on $\M$ preserving $\J$ with a
discrete
fixed-point set $\M^{S^1}$. Then, for each $P_i\in \M^{S^1}$, there is a well defined multiset of integers
$\{w_{i1},\ldots, w_{in}\}$, the weights of the $S^1$-representation on $T_{P_i}\M$, which determine this representation.
Indeed, there is an identification of $T_{P_i}\M$ with $\C^n$ such
that the $S^1$-action on $T_{P_i}\M$ is given by 
$$\alpha\cdot(z_1,\dots,z_n)=(\alpha^{w_{i1}}z_1,\ldots, \alpha^{w_{in}}z_n).$$
In this setting, motivated by the previous discussion,
we raise the following question.
\begin{question}\label{q3}
Let $(\M,\J)$ be a compact almost complex manifold. 
What are all the possible multisets of integers that can arise as weights of a $\J$-preserving $S^1$-action with isolated
fixed points? 
\end{question} 
After Petrie, there has been extensive work in this direction (see also \cite{Ah, GuZ1, GuZ2, Ha} and more recently \cite{L2, LL1, LL2, PT, T1}).
In particular,
let
$(\M,\omega)$
be a compact symplectic manifold of dimension $2n$, and suppose that the $S^1$-action is Hamiltonian with a minimal number of fixed points (i.e. $\lvert \M^{S^1}\rvert=n+1$).
In this case, since the set of compatible almost complex structures $\J\colon T\M\to T\M$ 
 is contractible, for each fixed point $P$ the set
of weights of the $S^1$-representation on $T_P\M$ and the total Chern class $\cc\in H^*(\M;\Z)$ do not depend on $\J$.
Then in \cite{T1}, Tolman  proves that  \emph{the $S^1$-representation on $T\M|_{\M^{S^1}}$ completely
determines the (equivariant) cohomology ring of $\M$ and the (equivariant) total Chern class.} 
In the case in which $\lvert \M^{S^1}\rvert$ is not minimal, determining the (equivariant) cohomology ring
from the fixed point set data is much harder, but it
has been investigated in several cases (cf. \cite{GT, GZ, ST}).

The multiset of weights of the $S^1$-representation on $T\M|_{\M^{S^1}}$ has to satisfy many rigid conditions,
coming, for example, from localization
theorems in equivariant cohomology and equivariant $K$-theory.
\textit{However, the type of equations that one gets from these techniques are
\textbf{``high-degree" polynomial equations}} that cannot be solved directly.
Nevertheless, they can be used to check whether a specific multiset of weights is admissible.
For example, as a consequence of the ABBV Localization formula \cite{AB,BV}, if the fixed point set $\M^{S^1}$ is discrete, we obtain
\begin{equation}\label{localized chern}
\sum_{i=0}^N \frac{\sigma_{j_1}(w_{i1},\ldots,w_{in})\cdots \sigma_{j_r}(w_{i1},\ldots,w_{in})}{\sigma_n(w_{i1},\ldots,w_{i_n})}=0\quad \mbox{for all}\quad j_1+\cdots +j_r<n\;,
\end{equation}
where $\sigma_j(x_1,\ldots,x_n)$ denotes the $j$-th elementary symmetric polynomial in $x_1,\ldots,x_n$, and $N+1$ is the number of fixed points (cf. Section~\ref{ec}).
If the dimension of $\M$ is ``small", equations \eqref{localized chern} can be still used to find all
the possible $S^1$-representations on $T\M|_{\M^{S^1}}$ (cf. \cite{Ah, T1}). However, as the dimension of $\M$ increases and/or
the number of fixed points gets larger, these equations become really unhandy. 
\\$\;$

In this paper we
introduce a new approach to the problem and give an explicit \emph{algorithm} that yields powerful
\emph{linear relations} among the weights.  The first important result in this direction concerns the Chern number $\cc_1\cc_{n-1}[\M]$.
\begin{theorem}\label{thm1}
Let $(\M,\J)$ be an almost complex manifold equipped with an $S^1$-action which preserves the almost complex structure $\J$ and has isolated fixed points.
For every $p=0,\ldots,n$, let $N_p$ be the number of fixed points with exactly $p$ negative weights.
Then
\begin{equation*}
\int_\M \cc_1\cc_{n-1}=\sum_{p=0}^n N_p[6p(p-1)+\frac{5n-3n^2}{2}].
\end{equation*}
In particular, if 
$N_p=1$ for every $p=0,\ldots,n$, then
\begin{equation*}
\int_\M \cc_1\cc_{n-1}=\frac{1}{2}n(n+1)^2.
\end{equation*}
\end{theorem}

Based on this fact, we are able to determine a family of vector spaces which contain the admissible lattices of
weights for $S^1$-actions satisfying  a given upper bound on the absolute value of the sum of the weights at each fixed point.
This upper-bound condition can be removed whenever $(\M,\omega)$ is a compact symplectic manifold
with a Hamiltonian $S^1$-action which satisfies a certain ``positivity condition",
which we will refer to as
$(\mathcal{P}_0^+)$ (see Definition \ref{P0+}).
This is satisfied  for example when, for each fixed point $P$, there exist exactly $\dim(\M)/2$
spheres containing $P$ which are fixed by a nontrivial subgroup of $S^1$, and $\cc_1$ is positive on each of them
(see  Remark~\ref{monotone..} and Section~\ref{refinement}). 

When the number of fixed points is minimal, $(\mathcal{P}_0^+)$
does not seem to be restrictive; indeed
there are no known examples of manifolds
which do not satisfy it (see Section~\ref{known examples}).
In particular, when $\dim(\M)$ is $4$ or $6$, it is known that $(\mathcal{P}_0^+)$ is always satisfied, and our algorithm
quickly determines all the possible families of isotropy weights, recovering the results of Ahara and Tolman
\cite{Ah, T1} for dimension $6$.
When $\dim(\M)=8$ and the action has exactly 5 fixed points, we prove that
\begin{itemize}
\item[(*)]  \emph{$(\mathcal{P}_0^+)$ is satisfied whenever the action extends to
an effective Hamiltonian $T^2$-action or none of the weights is 1.}
\end{itemize} 
(In fact, in this case, the action satisfies even a stronger condition  as it shown in Propositions~\ref{t2} and \ref{not pm1}.)

When $(\mathcal{P}_0^+)$ is not necessarily satisfied, we explore Hattori's results \cite{Ha} in the symplectic category.
In particular, we are able to get equations involving the Chern numbers of the manifold and the integral
of $y^n$, where $y$ is the generator of $H^2(\M;\Z)\simeq \Z$ for which $\cc_1=C_1y$, with 
$C_1\in \Z_{>0}$. When $\dim(\M)=8$,
using these equations, together with the algorithm mentioned above, we prove the following result.
\begin{thm}\label{RS1}
Let $(\M,\omega)$ be a compact symplectic manifold of dimension $8$, with a Hamiltonian $S^1$-action and $5$ fixed points.
Let $y$ be the generator of $H^2(\M;\Z)$ such that $\cc_1=C_1y$, for some positive integer $C_1$. Then $C_1$ can only be $1$ or $5$. Moreover the following are equivalent: 
\begin{itemize}
\item[(i)] $(\mathcal{P}_0^+)$ is satisfied.
\item[(ii)] $C_1=5$.
\item[(iii)] The cohomology ring agrees with the one of $\C P^4$, i.e. 
$$H^*(\M;\Z)=\Z[y]/(y^5),$$ 
where $y$ is of degree two.
\item[(iv)] The total Chern class agrees with the one of $\C P^4$, i.e. $\cc(T\M)=(1+y)^5$.
\item[(v)] The isotropy weights at the fixed points agree with the ones
of $\C P^4$ with the standard $S^1$-action.
\end{itemize}
\end{thm}

Combining Theorem~\ref{RS1} with (*) we obtain the following result.
\begin{corollary}\label{main dim 8}
Let $(\M,\omega)$ be a compact symplectic manifold of dimension $8$, with a Hamiltonian $S^1$-action
and $5$ fixed points, satisfying one of the following two conditions:
\begin{itemize}
\item[i)] the $S^1$-action extends to an effective Hamiltonian $T^2$-action;
\item[ii)] none of the weights of the action is $1$. 
\end{itemize}
Then the isotropy weights agree with the ones of the standard $S^1$-action on $\C P^4$. Moreover, the cohomology ring and Chern classes
agree with the ones of $\C P^4$, i.e.
$$
H^*(\M;\Z)=\Z[y]/(y^5)\quad\mbox{and}\quad \cc(T\M)=(1+y)^5\,,
$$
where $y$ has degree 2.
\end{corollary}
If $\M^{2n}$ is a connected compact K\"{a}hler manifold and has the same Betti numbers of $\C P^n$ but is different
from $\C P^n$, then it is called a \emph{fake projective space}. There has been extensive work towards the classification
of these spaces. In particular, it is known
\cite{W,Y} that, if $\M$ is a fake projective space of complex dimension $4$, then its Chern numbers
$(\cc_1^4, \cc_1\cc_3,\cc_1^2\cc_2,\cc_2^2,\cc_4)[\M]$ can only take one of the following two sets of values: 
$$(625,50,250,100,5)\quad \text{and} \quad (225,50,150,100,5).$$ 
However, it is not known whether there exists a connected K\"ahler manifold with the second
set of Chern numbers. 
This set falls into the case $C_1=1$ described above.
In current work in progress, we are trying to find a multiset of weights which would give this list of Chern numbers.
This would, in principle, allow us to construct a K\"ahler manifold with such values, as in \cite{M}. 

In what follows we give a brief description of the structure of the article.
In Section~\ref{BM}, we review some background material, establish some notation and recall fundamental facts about equivariant
cohomology and $K$-theory, and equivariant line bundles. 
In Section~\ref{c1cn-1}, we recall some facts about the Hirzebruch genus
of a compact almost complex manifold $(\M,\J)$ with a $\J$-preserving $S^1$-action with isolated fixed points, and prove Theorem~\ref{thm1}.
In Section~\ref{dmg}, we introduce a combinatorial object, called multigraph, which plays a crucial role 
in the algorithm, and encodes important information about the $S^1$-action (see Lemma~\ref{set of weights} and Proposition~\ref{magnitude}).
The main result of this section is Theorem~\ref{sum of m}, which 
gives a link between the combinatorics of the multigraph and the Chern number $\cc_1\cc_{n-1}[\M]$.
In Section~\ref{algo}, we introduce our algorithm, and discuss the positivity condition $(\mathcal{P}_0^+)$ mentioned above.
Section~\ref{mnfp} specializes to the case in which $(\M,\omega)$ is a compact symplectic manifold
with a Hamiltonian $S^1$-action with a minimal number of fixed points.
In particular, in Section~\ref{hamiltonian minimal}, we review some important results on the
equivariant cohomology ring and Chern classes, and derive Proposition~\ref{C1p} which will play an important role in the efficiency of our algorithm. In Section~\ref{refinement} we apply the preceding results to refine the algorithm and give conditions under which $(\mathcal{P}_0^+)$ is satisfied.
In Section~\ref{known examples} we give a list of known examples of $S^1$-Hamiltonian manifolds with a minimal number of fixed points,
and analyze some of their properties. In Section~\ref{hattori results} we analyze in detail the consequences of \cite{Ha} in the symplectic category,
and, in particular, what happens when $(\M,\omega)$ is  $8$-dimensional with $5$ fixed points but does not necessarily satisfy $(\mathcal{P}_0^+)$ (see Theorem~\ref{m not 2}).
Section~\ref{classification results} contains the classification results obtained using our algorithm on manifolds
satisfying $(\mathcal{P}_0^+)$. Its main result is Theorem~\ref{thm dim8}. Finally, at the end of this section we prove Theorem~\ref{RS1}, combining the results of the classification in Section~\ref{classification results} with the ones of Section~\ref{hattori results}. 

The accompanying software, based on the algorithms presented in this paper, can be found at 
{\bf \texttt{http://www.math.ist.utl.pt/$\sim$lgodin/MinimalActions.html}}.
 
\vspace{.5cm} 
\textbf{Acknowledgements.}\, We thank Tudor Ratiu for his support, Susan Tolman for introducing us to this problem, Victor Guillemin for useful discussions and Manuel Racle, Jos\'{e} Braga, Carlos Henriques and the students Filipe Casal, Francisco Pav\~{a}o Martins and Diogo Po\c{c}as  for helping us giving our first steps in C++ and Mathematica. 
\section{Background material}\label{BM}
In this section we will review some basic material and important results needed in this work and establish some notation.
\subsection{Equivariant cohomology and equivariant Chern classes}\label{ec}
(For a detailed discussion see, for instance, \cite{AB,GS}.)
Let $\M$ be a manifold endowed with a differentiable $S^1$-action. Let $ES^1$ be a contractible space on which
$S^1$ acts freely, and let $BS^1=ES^1/S^1$ be the classifying space. Then $ES^1$ can be identified with the unit sphere $S^{\infty}$ inside
$\C^{\infty}$ and $BS^1$ with $\C P^{\infty}$. Since the $S^1$-action on $ES^1$ is free, the diagonal action on
$\M \times ES^1$ is also free.
By the \textit{Borel construction}, the $S^1$-equivariant cohomology ring $H_{S^1}^*(\M)$ is defined to be the ordinary
cohomology of the orbit space $\M\times_{S^1} S^{\infty}$, 
$$
H_{S^1}^*(\M)=H^*(\M\times_{S^1}S^{\infty})\;.
$$
In particular, the $S^1$-equivariant cohomology of a point is given by 
$$H_{S^1}^*(pt;A)=H^*(BS^1;A)=H^*(\C P^{\infty};A)=A[x],$$
where $A$ is the coefficient ring and $x$ is of degree $2$.
The unique map $p\colon \M \to pt$ induces a map in equivariant cohomology 
\begin{equation}\label{module}
p^*\colon H_{S^1}^*(pt)\to H_{S^1}^*(\M),
\end{equation}
which gives $H_{S^1}^*(\M)$ the structure of an $H_{S^1}^*(pt)$-module.
Observe that the homomorphism $\{e\}\to S^1$ induces a restriction map in cohomology
\begin{equation}\label{restriction}
r\colon H_{S^1}^*(\M)\to H^*(\M)\;.
\end{equation}
Hence, as long as one knows the kernel and cokernel of $r$, the equivariant cohomology ring recovers information on the ordinary cohomology ring.

The projection onto the second factor $\pi\colon \M\times_{S^1}S^{\infty}\to \C P^{\infty}$ gives rise to a
push-forward map 
\begin{equation}\label{push forward}
\pi_*\colon H^*_{S^1}(\M) \to H^{*-\dim(\M)}(\C P^{\infty})
\end{equation}
which can be thought as an integration along the fibers of $\pi$. We will denote it by $\int_{\M}$.

Let $E\to \M$ be an $S^1$-equivariant  vector bundle. Then the \emph{equivariant Euler class} $e^{S^1}(E)$ of $E$ is defined as the Euler class
of the bundle $E\times_{S^1}S^{\infty}\to \M \times_{S^1}S^{\infty}$.
If $E$ is a complex vector bundle, the \emph{equivariant Chern classes} $\cc^{S^1}_i(E)$ are the Chern classes of $E\times_{S^1}S^{\infty}\to \M\times_{S^1}S^{\infty}$,
and $r(\cc_i^{S^1})=\cc_i$ for every $i$, where $r$ is the restriction map \eqref{restriction}.

Let $\M^{S^1}$ be the set of the  $S^1$-fixed points, and $F\subset \M^{S^1}$ one of its connected components. Since $F$ is $S^1$-invariant, the inclusion $\iota_F\colon F \hookrightarrow \M$
induces a restriction map in equivariant cohomology $\iota_F^*\colon H_{S^1}^*(\M)\to H_{S^1}^*(F)$.
The Atiyah-Bott-Berline-Vergne Localization formula for $S^1$-actions allows to compute 
the push-forward map \eqref{push forward} in terms of the fixed point set data (cf. \cite{AB,BV}).
\begin{theorem}[ABBV Localization formula]
Let $\M$ be a compact oriented manifold endowed with a smooth $S^1$-action.
Given $\mu\in H_{S^1}^*(\M;\Q)$
\begin{equation*}
\int_\M \mu= \sum_{F}\int_F\frac{\iota_F^*(\mu)}{e^{S^1}(\N_F)}\;,
\end{equation*}
where the sum is over all the fixed-point set components of the action, and
$e^{S^1}(\N_F)$ is the equivariant Euler class of the normal bundle to $F$.
\end{theorem} 
This localization formula becomes particularly easy when the fixed point set is discrete, i.e.  
$\M^{S^1}=\{P_0,\ldots,P_N\}$. In this case, the normal bundle to a fixed point $P$ is just $T_P\M$ and, since $\mathbf{0}$ is the only fixed point of the isotropy representation of $S^1$ on $T_P\M$, it follows that
$T_P\M$ has a canonical orientation and is an even dimensional vector space of dimension $2n$.
 Let 
$$e^{S^1}(T_P\M)\in H_{S^1}^{*}(\{P\};\Z)\simeq \Z[x]$$ 
be the equivariant Euler class of $T_P\M$. 
Let $w_{1P},\ldots,w_{nP}$ be the weights of the isotropy representation of $S^1$ on $T_P\M$. Even if the sign of the individual weights is not well defined, the sign of the product is, and
a standard computation shows that 
$$e^{S^1}(T_P\M)=(\displaystyle\prod_{i=1}^nw_{iP})x^n,$$ 
where $x$ is of degree $2$.
Then the ABBV Localization formula reduces to the following.
\begin{corollary}\label{abbv discrete}
Let $\M$ be a compact oriented manifold endowed with a smooth $S^1$-action
such that $\M^{S^1}=\{P_0,\ldots,P_N\}$.
Given $\mu\in H_{S^1}^*(\M;\Q)$
\begin{equation}\label{abbv dis}
\int_{\M} \mu= \sum_{i=0}^N\frac{\mu(P_i)}{(\prod_{j=1}^nw_{ij})x^n}\;,
\end{equation}
where $w_{i1},\ldots,w_{in}$ are the weights of the $S^1$-isotropy representation on $T_{P_i}\M$
and $\mu(P)=\iota_{\{P\}}^*(\mu)$ for all $P\in \M^{S^1}$.
\end{corollary}
\noindent Notice that on the right hand side of \eqref{abbv dis}, each term is not an element of $\Q[x]$, but their sum is.

Suppose that $(\M,\J)$ is an almost complex manifold equipped with a $\J$-preserving  $S^1$-action with isolated fixed points $P_0,\ldots,P_N$. Then the signs of the individual weights of the $S^1$-representation on $T_{P_i}\M$
are well defined, and a standard computation shows that
the restriction of the $j$-th equivariant Chern class to $P_i$ is given by
$$\cc_j^{S^1}(P_i)=\sigma_j(w_{i1},\ldots,w_{in})x^j\in H_{S^1}^{2j}(\{P_i\};\Z)\;,$$ where $\sigma_j$ denotes the $j$-th elementary symmetric polynomial.
\\$\;$

\subsubsection{The symplectic case}\label{ec symplectic}
Let us now assume that $(\M,\omega)$ is a compact symplectic manifold of dimension $2n$, endowed with a symplectic $S^1$-action with isolated fixed points $P_0,\ldots,P_N$. 
Let $\J\colon T\M\to T\M$ be an almost complex structure compatible with $\omega$, i.e. $\omega(\cdot,\J \cdot)$ 
is a Riemannian metric. Since the set of such structures is contractible, the set of
weights of the isotropy representation of $S^1$ on $T_{P_i}\M$ is well defined for every $P_i\in \M^{S^1}$. Let $w_{i1},\ldots,w_{in}$ be the multiset of these weights. Then we can identify $T_{P_i}\M$ with $\C^n$, and the
$S^1$-action on $T_{P_i}\M$ with the $S^1$-action on $\C^n$ given by 
$$\alpha\cdot(z_1,\dots,z_n)=(\alpha^{w_{i1}}z_1,\ldots, \alpha^{w_{in}}z_n).$$ 
Hence $T_{P_i}\M\simeq\bigoplus_{j=1}^n \C_{w_{ij}}$, where $\C_{w_{ij}}$ is the
one dimensional complex subspace on which $S^1$ acts with weight $w_{ij}$.
For each $i=0,\ldots,N$, we denote $\bigoplus_{w_{ij}<0}\C_{w_{ij}}$ by
 $\N_{P_i}^-$.  

Consider the case in which the action is also
\textbf{Hamiltonian}, i.e.  there exists an $S^1$-invariant function $\psi\colon \M \to \R$, called \textbf{moment map},
satisfying
$$
d \psi = -\iota_{\xi^{\#}}\omega,
$$
where $\xi^{\#}$ denotes the vector field generated by the $S^1$-action, and $\iota_{\xi^{\#}}$ is the interior derivative.
Then $\psi\colon \M \to \R$ is a {\bf perfect Morse function} whose critical set coincides with the fixed point set $\M^{S^1}$. Hence, for every $P_i\in \M^{S^1}$, the negative normal bundle at $P_i$ is precisely $\N_{P_i}^-$,
 and the (Morse) index at $P_i$ is $2\lambda_i$, where $\lambda_i$ is the \emph{number of negative weights at $P_i$}.
 
Notice that the existence of a moment map $\psi$ gives rise to a
natural equivariant extension of $\omega$, i.e. 
$\omega - \psi \otimes x$, which is an $S^1$-invariant form, closed under the differential $d_{S^1}=d\otimes 1-\iota_{\xi^{\#}}\otimes x$ of the Cartan complex. Consequently,
it represents a class 
$$[\omega -\psi\otimes x]\in H_{S^1}^2(\M;\R).$$
The invariant form $\omega-\psi \otimes x$ is called \emph{equivariant symplectic form}.
  
Kirwan \cite{Ki} uses the existence of such a map to prove very nice properties for the equivariant cohomology ring $H_{S^1}^*(\M;\Z)$.
First of all, if $\iota\colon \M^{S^1}\to \M$ denotes the inclusion of the fixed point set into $\M$,
the map
\begin{equation}\label{inj}
\iota^*\colon H_{S^1}^*(\M;\Z)\to H_{S^1}^*(\M^{S^1};\Z)
\end{equation}
is \emph{injective}.
Hence, any equivariant cohomology class $\gamma\in H_{S^1}^*(\M;\Z)$ is completely determined
by its restriction to the fixed points.

Moreover, the restriction map \eqref{restriction} to the ordinary cohomology ring
is \emph{surjective}, and the kernel is given by the ideal generated by $p^*(x)$, where $p^*$
is the map \eqref{module} and $H^*(\C P^{\infty};\Z)\simeq \Z[x]$. In the following we will simply denote $p^*(x)$ by $x$.

In addition, the number of fixed points of index $i$ equals the rank of $H^i(\M;\Z)$, which is the $i$-th Betti number  $b^i(\M)$ of $\M$. More precisely, if $N_p$ is the number of fixed points of index $2p$ for every $p=0,\ldots,n$, then $H^{2p}(\M;\Z)=\Z^{N_p}$
 and  zero otherwise. It is then easy to see, by reversing the circle action, that we must have
 \begin{equation}\label{eq:reverse}
 N_p=N_{n-p}.
 \end{equation}
 Suppose that $[\omega]$ belongs to the image of the map $H^2(\M;\Z)\to H^2(\M;\R)$.
 Since $[\omega]^k\neq 0$ for every $k=0,\ldots,n$, it follows, by the above result,  that $N_p\neq 0$ for every $p=0,\ldots,n$.
Thus 
\begin{quote}
\emph{the minimal number of fixed points on a compact Hamiltonian \\manifold $\M$ of dimension $2n$ is $n+1$}, and, in this
case, 
$$H^i(\M,\Z)=H^i(\C P^n;\Z) \quad i=0,\ldots,n. $$ 
\end{quote}
A basis for the equivariant cohomology of $\M$ is given as follows.
Let us define
\begin{equation}\label{eq:lambda-}
\Lambda_i^-=(\prod_{w_{ij}<0}w_{ij})x^{\lambda_i},
\end{equation}
so that $\Lambda_i^-$ coincides with $e^{S^1}(\N_{P_i})$, the equivariant Euler class of the negative normal bundle at $P_i$. Accordingly, we also define 
\begin{equation}\label{eq:lambda+}
\Lambda_i^+=(\prod_{w_{ij}>0}w_{ij})x^{n-\lambda_i}
\end{equation}
and $\Lambda_i=e^{S^1}(T\M_{P_i})=\Lambda_i^+\Lambda_i^-$, for every $i=0,\ldots,N$.
\begin{lemma}[Kirwan, \cite{Ki}]\label{kirwan}
Let $(\M,\omega)$ be a compact symplectic manifold endowed with a Hamiltonian $S^1$-action with isolated
fixed points $P_0,\ldots,P_N$. Let $\psi\colon \M \to \R$ be the corresponding moment map.
For every fixed point $P_i\in \M^{S^1}$ there exists a class $\gamma_i\in H_{S^1}^{2\lambda_i}(\M;\Z)$ such that
\begin{itemize}
\item[(1)] $\gamma_i(P_i)= \Lambda_i^-$;
\item[(2)] $\gamma_i(P_j)=0$ for every $P_j\in \M^{S^1}\setminus\{P_i\}$ such that $\psi(P_j)\leq\psi(P_i)$.
\end{itemize}
Moreover, for any such classes,  $\{\gamma_i\}_{i=0}^N$ is a basis for $H_{S^1}^*(\M;\Z)$ as a module over $H^*(\C P^{\infty};\Z)=\Z[x]$.
\end{lemma} 
Notice that the set of classes satisfying $(1)$ and $(2)$ is not unique. In fact, if there exist $P_l$ and $P_m$ such that
$\psi(P_l)< \psi(P_m)$ and $\lambda_l=\lambda_m$, then the class $\gamma^\prime_l=\gamma_l+k\gamma_m$ satisfies
the same properties $(1)$ and $(2)$ satisfied by $\gamma_l$, for any $k\in \Z$.  

Fix a set of classes $\{\gamma_i\}_{i=0}^N$ satisfying $(1)$ and $(2)$ of Lemma~\ref{kirwan}. Since they form a basis for $H_{S^1}^*(\M;\Z)$
as a module over $H^*(\C P^{\infty};\Z)$, given any class $\alpha\in H_{S^1}^*(\M;\Z)$ there exist
$\alpha^0,\ldots, \alpha^N\in H^*(\C P^{\infty};\Z)$ such that
$
\alpha=\sum_{j=0}^N\alpha^i\gamma_i\;.
$
The next Lemma gives a recursive formula that computes the coefficients $\alpha^i$s in terms
of $\iota^*(\alpha)$ and $\iota^*(\gamma_i)$, for $i=0,\ldots,N$, which  is an immediate consequence of properties $(1)$ and $(2)$.
\begin{lemma}\label{coefficients} Let us order the fixed points $P_0,\ldots,P_N$ in such a way that
$$\psi(P_0)<\psi(P_1)\leq \psi(P_2)\leq \cdots \leq\psi(P_{N-1})<\psi(P_N).$$
Then the coefficients $\alpha^i$s can be computed recursively as
$$
\alpha^i = \frac{\alpha(P_i)-\sum_{h:\;\psi(P_h)<\psi(P_i)}\alpha^h\gamma_{h}(P_i)}{\Lambda_i^-}.
$$
\end{lemma}
For every $0\leq j\leq l\leq N$ and $0\leq i \leq N$, the elements $\alpha_{jl}^i$ in $H^*(\C P^{\infty};\Z)$ such that
$$
\gamma_j\gamma_l=\sum_{i=0}^N\alpha_{jl}^i\gamma_i\;
$$
are 
called the \emph{equivariant structure constants} of $H_{S^1}^*(\M;\Z)$ with respect to the basis $\{\gamma_0,\ldots,\gamma_N\}$. 
In order to compute them, by Lemma~\ref{coefficients}, it is sufficient to compute $\iota^*(\gamma_i)$
for every $i$. This problem has been extensively studied in the literature, and  it is possible to get an explicit
formula for these restrictions only in very special cases (see for example \cite{GZ,GT,ST}).
However, when the number of fixed points is minimal,
Tolman \cite{T1} shows that one can give an explicit basis for $H_{S^1}^*(\M;\Z)$ whose restriction to the fixed point set
can be completely recovered from the fixed point set data (see Section~\ref{hamiltonian minimal} for details).

\subsection{$K$-theory and equivariant $K$-theory}
(For a detailed discussion see, for example, \cite{At, AS, AS3}.)
Let $(\M,\J)$ be a compact almost complex manifold endowed with a $\J$-preserving $S^1$-action with isolated fixed points $P_0,\ldots,P_N$.
We recall that $K(\M)$ (resp. $K_{S^1}(\M)$) is the abelian group associated to the
semigroup of isomorphism classes of complex vector bundles (resp. complex $S^1$-vector bundles) over $\M$,
endowed with the direct sum operation $\oplus$. This also has a
ring structure, given by the tensor product $\otimes$.
If $\M$ is a point, we have 
$$K(pt)\simeq \Z\quad \text{and} \quad K_{S^1}(pt) \simeq R(S^1),$$ 
the character ring of $S^1$. This last ring can be simply identified with the Laurent
polynomial ring $\Z[t,t^{-1}]$, where $t$ denotes the standard $S^1$-representation. 

In analogy with what happens in cohomology, the unique map
$ \M\to pt$ induces  maps in equivariant and in ordinary $K$-theory,
\begin{equation*}
K_{S^1}(pt)\to K_{S^1}(\M)\quad\text{and} \quad K(pt)\to K(\M),
\end{equation*}
which give $K_{S^1}(\M)$ the structure of a $\Z[t,t^{-1}]$-module, and $K(\M)$ the structure of a $\Z$-module.

Moreover, the inclusion homomorphism $\{e\}\hookrightarrow S^1$ induces a restriction map in
$K$-theory
\begin{equation*}
r\colon K_{S^1}(\M)\to K(\M)\;,
\end{equation*}
which, when $\M$ is a point,  is just the evaluation map, $r\colon \Z[t,t^{-1}]\to \Z$, at $t=1$.

Consider the $K$-theoretic push-forward map in equivariant and ordinary $K$-theory,
namely the index homomorphisms
\begin{equation}\label{equiv ind}
\ind_{S^1}\colon K_{S^1}(\M)\to K_{S^1}(pt)\simeq \Z[t,t^{-1}]
\end{equation}
and
\begin{equation}\label{index}
\ind\colon K(\M)\to K(pt)\simeq \Z\;.
\end{equation} 
By the Atiyah-Singer formula, \eqref{index} can be computed as 
\begin{equation}\label{AT formula}
\ind(\eta)=\int_\M \ch(\eta)\ttot(\M)\;,\quad\mbox{for every}\quad\eta\in K(\M),
\end{equation}
where $\ch(\cdot)$ denotes the Chern character $\ch\colon K(\M)\to H^*(\M;\Q)$, and $\ttot(\M)$ is the total Todd class of $\M$.

On the other hand, by the Atiyah-Segal formula \cite{AS}, the  map \eqref{equiv ind} can be computed in terms of the fixed-point set data. Namely, if $w_{i1},\ldots,w_{in}$ are the weights of the $S^1$-representation on $T_{P_i}\M$,
\begin{equation}\label{AS formula}
\ind_{S^1}(\eta^{S^1})=\sum_{i=0}^N \frac{\eta^{S^1}(P_i)}{\prod_{j=1}^n (1-t^{-w_{ij}})}\,,\quad\mbox{for every}\quad \eta^{S^1}\in K_{S^1}(\M)\;,
\end{equation} 
where $\eta^{S^1}(P_i)\in \Z[t,t^{-1}]$ denotes the restriction of $\eta^{S^1}$ to the fixed point $P_i$.
Notice that the sum on  right hand side of \eqref{AS formula} is in $\Z[t,t^{-1}]$, despite the fact that each of its terms is not.
Observe also  that  we have the following commutative diagram
\begin{equation}\label{K commutes}
\xymatrix{
K_{S^1}(\M) \ar[r]^{r} \ar[d]_{\ind_{S^1}} & K(\M) \ar[d]_{\ind} \\
 \Z[t,t^{-1}] \ar[r]^{r} &   \Z.
 } \
\end{equation}
Thus, the value of $\ind_{S^1}(\eta^{S^1})$ at $t=1$ can be computed using \eqref{AT formula} and \eqref{AS formula}:
\begin{equation*}
\ind_{S^1}(\eta^{S^1})_{\rvert_{t=1}}=
\int_\M \ch(r(\eta^{S^1}))\ttot(\M)\;.
\end{equation*}
\subsection{Equivariant complex line bundles}\label{eclb}
(For a detailed discussion see, for example, \cite{Ha, HL, HY, Mu} and \cite[Appendix C]{GKS}.)
Let $\M$ be a compact manifold with a smooth $S^1$-action and
let $\LL$ be a complex line bundle over $\M$. Then $\LL$ is called \emph{admissible} if the $S^1$-action on $\M$
lifts to an $S^1$-action on $\LL$ which makes the projection map $\LL\to \M$ equivariant.
We denote by $\LL^{S^1}$ the line bundle $\LL$ endowed with the lifted $S^1$-action.
\begin{theorem}\label{lift clb}
Let $\;\LL^{S^1}\to \M$ be an $S^1$-equivariant complex line bundle over $\M$, and let $\cc_1^{S^1}(\LL^{S^1})\in H^2_{S^1}(\M;\Z)$
be its equivariant first Chern class. Then $\cc_1^{S^1}$ determines a one-to-one correspondence between equivalence
classes of  $S^1$-equivariant complex line bundles over $\M$ and elements of $H_{S^1}^2(\M;\Z)$.
\end{theorem}
Consequently, a complex line bundle $\LL\to \M$ is admissible if and only if $\cc_1(\LL)$ is in the image of the restriction
map $r\colon H_{S^1}^2(\M;\Z)\to H^2(\M;\Z)$.
Moreover, all different liftings of the $S^1$-action on $\LL$ are parametrized by $H^2(\C P^{\infty};\Z)\simeq \Z$. 
In particular, suppose that the $S^1$-action has isolated fixed points $P_0,\ldots,P_N$ and
let $\LL$ be an admissible complex line bundle. Then, as mentioned above, the lift $\LL^{S^1}$ is not
uniquely determined, but, for any such lift,
there exists an integer $a$ such that
the restriction of $\LL^{S^1}$ to a fixed point $P_i$ is of the form
$$
\LL^{S^1}(P_i)=t^{a_i+a}\quad\mbox{for every}\quad i=0,\ldots,N,
$$
for fixed integers $a_0,\ldots,a_N$, where $t$ denotes the standard 1-dimensional $S^1$-representation.
More precisely, the integer $a_i$ is given by 
$$ a_i = \frac{(\cc_1^{S^1}(\LL^{S^1}))(P_i)}{x}, \quad i=0,\ldots,N. $$
Hence, the values of $\LL^{S^1}$ at the fixed points  are determined up to a constant representation. 
\section{The Chern number $\cc_1\cc_{n-1}[\M]$}\label{c1cn-1}

Let $(\M,\J)$ be a compact almost complex manifold, and let $S^1$ be a circle acting on $\M$ preserving the
almost complex structure $\J$. In this section we show that the Chern number $\int_\M \cc_1\cc_{n-1}$
can be completely determined by the fixed-point set data. When $\M$  is a compact complex manifold,
 Libgober and Wood  \cite{LW} proved that this integral is determined by the Hodge numbers
of $\M$. The same fact was also shown later  by Borisov in \cite{Bo}, inspired by previous work of Eguchi, Hori and Xiong
on Fano varieties \cite{EHX}. 
\\$\;$

Let $(\M,\J)$ be a compact almost complex manifold of dimension $2n$, and let 
$\cc_j\in H^{2j}(\M;\Z)$ be the Chern classes of the tangent bundle, for every $j=0,\ldots,n$.
Let $\chi_y(\M)$ be the Hirzebruch genus of $\M$, i.e. the genus corresponding to the power series
$$
Q_y(x)=\frac{x(1+ye^{-x(1+y)})}{1-e^{-x(1+y)}}
$$
(cf. \cite{Hi} for details). Then 
$$\chi_y(\M)=\displaystyle\sum_{i=0}^n\left(\int_\M T_i^n\right) y^i,$$ 
where $T_i^n$ is a rational combination of products of Chern classes
 $$\cc_{j_1}\cdots \cc_{j_r}\;,\;\;\mbox{ with }\;\;j_1+\cdots +j_r=n.$$ 
For example, up to order $4$, the $T_i^j$s are
\begin{align*}
T_0^0 &=1, \quad \quad \quad \quad \,T_0^1 =-T_1^1=\displaystyle\frac{1}{2} \cc_1, \quad \,\, T_0^2 =T_2^2=\displaystyle \frac{\cc_1^2+\cc_2}{12}, \\ T_2^1  & = \displaystyle\frac{\cc_1^2-5\cc_2}{6}, \quad 
T_0^3 =-T_3^3=\displaystyle\frac{\cc_1\cc_2}{24},    \quad  T_1^3   =-T_2^3=\displaystyle\frac{\cc_1\cc_2-12\cc_3}{24}, \\
T_0^4 &= T_4^4=\displaystyle\frac{-\cc_1^4+4\cc_1^2\cc_2+3\cc_2^2+\cc_1\cc_3-\cc_4}{720},  \\ T_1^4  & =T_3^4=\displaystyle\frac{-\cc_1^4 + 4 \cc_1^2 \cc_2 + 3 \cc_2^2 - 14 \cc_1 \cc_3 - 31 \cc_4}{180}, \\
T_2^4 &=\displaystyle\frac{-\cc_1^4 + 4 \cc_1^2 \cc_2 + 3 \cc_2^2 - 19 \cc_1 \cc_3 + 79 \cc_4}{120}.
\end{align*}
We recall that the Hirzebruch genus recovers three important topological invariants of the manifold $\M$:
\begin{itemize}
\item {\bf If ${\bf y=0}$} then  $\chi_0(\M)=\td(\M)$ is the Todd genus of $\M$, i.e. the genus associated to the power series $Q_0(x)=\displaystyle\frac{x}{1-e^{-x}}$. Hence   $T_0^n$ is just the Todd polynomial of degree $n$. 
\item {\bf If ${\bf y=1}$} then $\chi_1(\M)=\sign(\M)$ is the signature of $\M$.
\item {\bf If ${\bf y=-1}$} then  $\chi_{-1}(\M)=\int_\M\cc_n$ is the Euler characteristic of $\M$. Moreover, if $\M$ is a complex manifold, $\chi_{-1}(\M)=\ind(\overline{\partial})$ is the index of the $\overline{\partial}$ operator. 
\end{itemize}

Let $T^*\M\otimes \C=T^{1,0}\M\oplus T^{0,1}\M$ be the splitting of the complexified cotangent bundle  induced by $\J$, into its holomorphic
and antiholomorphic parts.  For every $p=0,\ldots,n$, let $\chi_p(\M)$ be the topological index of the bundle $\Lambda^pT^{1,0}(\M)$, regarded as an element of $K(\M)$, i.e. $\chi_p(\M)=\ind(\Lambda^pT^{1,0}\M)\in K(pt)$. Then by the Atiyah-Singer formula  \cite{AS3}
$$
\chi_p(\M)=\int_\M\ch(\Lambda^pT^{1,0}\M)\ttot(\M)\;,
$$
where the orientation on $\M$ is the one induced by $\J$, $\ch(\cdot)$ denotes the Chern character, and $\ttot(\M)$ is the total Todd class of $\M$, i.e.
\begin{equation}\label{total todd class}
\ttot(\M)=T_0^0+T_0^1+\cdots +T_0^n\quad\in H^*(\M;\Q)\;.
\end{equation}
If $\cc(T\M)=\prod_{i=1}^n(1+x_i)$ is a formal factorization of the total Chern class, then a standard computation shows that
$$\displaystyle\sum_{p=0}^n \ch(\Lambda^pT^{1,0}\M)y^p=\displaystyle \prod_{i=1}^n(1+ye^{-x_i})$$ 
and
\begin{equation}\label{chip}
\chi_y(\M)=\sum_{p=0}^n\chi_p(\M)y^p.
\end{equation}

Now suppose that $\M$ is equipped with a circle action preserving the almost complex structure $\J$ such that the fixed point
set $\M^{S^1}$ is discrete. Then, for all $p=0,\ldots,n$, the bundle $\Lambda^pT^{1,0}(\M)$ can be regarded as an element of $K_{S^1}(\M)$. Let us define
$$\chi_y(\M,t):=\sum_{p=0}^n \ind_{S^1}\left(\Lambda^pT^{1,0}(\M)\right)y^p.$$
Then, following an idea of Atiyah-Hirzebruch \cite{AH}, Kosniowski \cite{K} and  Lusztig \cite{L},  Li proves  in \cite{L2} that
when $(\M,\J)$ is a compact almost complex manifold with a $\J$-preserving $S^1$-action with isolated fixed points $P_0,\ldots,P_N$,
then $\chi_y(\M,t)$ is independent of $t$, and has the following explicit expression.
Let $\lambda_i$ denote the number of negative weights at
$P_i$, then
$$
\chi_y(\M,t)=\sum_{i=0}^N (-y)^{\lambda_i}
$$
(see also Section 5.7 of \cite{Hi}).
Since $\chi(\M)=\chi(\M,1)$, we have that   
\begin{equation}\label{chi}
\chi_y(\M)=\sum_{i=0}^N (-y)^{\lambda_i}.
\end{equation}
Let $N_p$ be the number of fixed points with exactly $p$ negative weights. Then we can rewrite \eqref{chi} as
\begin{equation}\label{chi2}
\chi_y(\M)=\sum_{p=0}^n N_p(-y)^p\;,
\end{equation}
and so by \eqref{chip} it follows that $\chi_p(\M)=(-1)^pN_p$.
Moreover, since 
$$\sum_{i=0}^N(-y)^{\lambda_i}=\sum_{i=0}^N(-y)^{n-\lambda_i},$$ 
we have
$$\chi_p(\M)=(-1)^n\chi_{n-p}(\M).$$
From $\chi_{-1}(\M)=\int_\M \cc_n$ and \eqref{chi2} it immediately follows  that
\begin{equation}\label{cn}
\int_\M \cc_n=\sum_{p=0}^n N_p\;.
\end{equation}

We are now able to prove Theorem~\ref{thm1}, which
 shows that there exists another Chern number which only depends on the fixed-point set
data. For that, we adapt the results in \cite{LW} and \cite{Bo} to the almost complex case, and combine them with the localization formula \eqref{chi2} for the Hirzebruch genus.
\begin{proof}[Proof of Theorem~\ref{thm1}]
The proof follows closely the argument found in \cite[Proposition 2.2.]{Bo}.
By \eqref{chip},  we have
\begin{equation}\label{dchi}
\displaystyle\frac{d^2\chi_y(\M)}{dy^2}|_{y=-1}=\displaystyle\sum_{p=0}^n\chi_p(\M)(-1)^pp(p-1).
\end{equation}
On the other hand, we can use the Atiyah-Singer formula to express 
$$\displaystyle\frac{d^2\chi_y(\M)}{dy^2}|_{y=-1}$$ 
as a combination of Chern numbers. Indeed,
\begin{gather*}
\frac{d^2\chi_y(\M)}{dy^2}=\int_\M \frac{d^2}{dy^2}\left( \sum_{p=0}^n \ch(\Lambda^pT^{1,0}\M)y^p\right)\ttot(\M)=\int
_\M \frac{d^2}{dy^2}\left(\prod_{i=1}^n \frac{x_i(1+ye^{-x_i})}{1-e^{-x_i}}\right).\\
\end{gather*}
Hence,
\begin{align*}
\frac{d^2\chi_y(\M)}{dy^2}|_{y=-1} & =2\sum_{j<k}\int_\M \left(\frac{x_je^{-x_j}}{1-e^{-x_j}}\right) \left(\frac{x_ke^{-x_k}}{1-e^{-x_k}}\right)\prod_{h\neq j,k}x_h \\ 
& = 2\sum_{j<k}\int_\M \left(1-\frac{x_j}{2}+\frac{x_j^2}{12}\right)\left(1-\frac{x_k}{2}+\frac{x_k^2}{12}\right)\prod_{h\neq j,k}x_h \\
& =   2\sum_{j<k}\int_\M \left(\frac{x_jx_k}{4}+\frac{x_j^2+x_k^2}{12}\right) \prod_{h\neq j,k}x_h\;\;.
\end{align*}
Moreover, it is  easy to see that
$$
\int_\M \cc_1\cc_{n-1}=n\int_\M\cc_n + \sum_{j<k}\int_\M(x_j^2+x_k^2)\prod_{h\neq j,k}x_h\;,
$$
and so, the previous computation along with \eqref{dchi}, yields
$$
\displaystyle\sum_{p=0}^n \chi_p(\M)(-1)^pp(p-1)=\frac{1}{6}\int_\M \cc_1\cc_{n-1}+\frac{3n^2-5n}{12}\int_\M \cc_n\;.
$$
The result then follows from the fact that $\chi_p(\M)=(-1)^pN_p$ and from \eqref{cn}.
\end{proof}
In virtue of the observations made in Section~\ref{ec}, Theorem~\ref{thm1} has the following immediate
corollary.
\begin{corollary}\label{corollary1}
Let $(\M,\omega)$ be a compact symplectic manifold with a Hamiltonian $S^1$-action and isolated fixed points.
Let $b^i(\M)$ be the $i$-th Betti number of $\M$.
Then
\begin{equation*}
\int_\M \cc_1\cc_{n-1}=\sum_{p=0}^n b^{2p}(\M)[6p(p-1)+\frac{5n-3n^2}{2}].
\end{equation*}
In particular, if the number of fixed points is minimal,
\begin{equation}\label{eq:c1cn-1minimal}
\int_\M \cc_1\cc_{n-1}=\frac{1}{2}n(n+1)^2.
\end{equation}
 \end{corollary}

\begin{tiny}
\end{tiny}\section{Defining the multigraphs}\label{dmg}
\subsection{Abstract multigraphs}
An \textbf{oriented multigraph} $\Gamma$ is an oriented graph with multiple oriented edges. More precisely, 
it is an ordered pair $\Gamma=(V,E)$, where
\begin{itemize}
\item $V$ is a set of \textit{vertices},
\item $E$ is a multiset of ordered pairs of vertices which we will call  the \textit{edges} of $\Gamma$.
\end{itemize}
Let $\ii\colon E \to V$ (resp. $\te\colon E \to V$) be the map which associates to each edge $e$ its \emph{initial} (resp. \emph{terminal}) point. Note that we allow a multigraph to contain cycles, i.e. edges $e$ with $\ii(e)=\te(e)$.
We denote by $E^{\circlearrowleft}$ the subset of $E$ formed by cycles, and by $E^{\circlearrowleft}_P$ the set of cycles that start and end at $P\in V$.

For every $P\in V$, let $E_{P,\ii}$ (resp. $E_{P,\te}$ ) be the set of edges whose initial (resp. terminal) point is $P$, i.e. $E_{P,\ii}=\{e\in E \mid \ii(e)=P\}$ (resp. $E_{P,\te}=\{e\in E \mid \te(e)=P\}$), and define $E_P$ to be the multiset
$E_{P,\ii}\cup E_{P,\te}$. This is a multiset in the sense that, if $e$ is a cycle, it appears twice in
$E_P$: once as an element of $E_{P,\ii}$ and once as an element of $E_{P,\te}$. 
We say that the multigraph $\Gamma=(V,E)$ has \emph{degree} $n$ if $\lvert E_P\rvert =n$ for all $P\in V$. Moreover,
for every $P\in V$, we define $\lambda(P)$ to be the number of edges ending at $P$, i.e. $\lambda(P)=\lvert E_{P,\te}\rvert$.

Define $\delta\colon V\times E\to \{-1,0,1\}$ to be the map 
$$
\delta(P,e)=
\begin{cases} 
\;\;1 & \mbox{if }e\in E_{P,\ii}\setminus E_{P,\te}\\
-1 & \mbox{if }e\in E_{P,\te}\setminus E_{P,\ii}\\
0 & \mbox{otherwise}.\\
\end{cases}
$$
Assume that the edge set $E$ contains finitely many elements, and let us fix an ordering on 
$E=(e_1,\ldots,e_{\lvert E \rvert })$. To every oriented multigraph $\Gamma$ we associate an $(\lvert E\rvert \times \lvert E \rvert)$-\emph{matrix} $A(\Gamma)$
whose element at position $(h,m)$ is given by 
\begin{equation}\label{matrix}
a_{h,m}:=\delta \left(\ii(e_h),e_m\right)-\delta\left(\te(e_h),e_m\right).  
\end{equation}
We say that a multigraph $\Gamma=(V,E)$ is \textbf{labeled} by $\w$ (and denote it by $(\Gamma,\w)$) if we associate to $\Gamma$ a map $\w\colon E\to \Z$.
Finally, suppose that $(\Gamma,\w)$ satisfies $\w(e)\neq 0$ for all $e\in E$. Then we define the \textbf{magnitude}  $\m\colon E\to \Z$  of $(\Gamma,\w)$ as
$$
\m(e):=\frac{\displaystyle\sum_{f\in E}\left(\delta\left(\ii(e),f\right)-\delta\left(\te(e),f\right)\right)\w(f)}{\w(e)}\;.
$$ 
The subset of edges $e\in E$ with $\m(e)>0$ (resp. $\m(e)<0$) is denoted by $E^+$ (resp. $E^-$), and its elements are called 
\textbf{positive edges} (resp. 
\textbf{negative edges}). Finally we denote by $E^0$ the subset of edges with $\m(e)=0$.

For every multigraph $\Gamma$ we can introduce an equivalence relation on the set of vertices: given $P,Q\in V$ we say that
$P\sim Q$ if and only if there exists a sequence of unoriented edges connecting $P$ to $Q$. The corresponding equivalence classes are
called the \emph{connected components} of $\Gamma$. 

\subsection{The multigraphs associated to an $S^1$-action}
Let $(\M,\J)$ be a compact almost complex manifold of dimension $2n$, equipped with a $\J$-preserving circle action
with a discrete fixed point set $\M^{S^1}$.
In this section we will associate a  family of oriented multigraphs to the  $S^1$-space  $(\M,\J,S^1)$,  which will encode information about the fixed-point set data.

Let $\M^{S^1}=\{P_0,\ldots,P_N\}$, and $w_{i1},\ldots,w_{in}$ the weights of the isotropy representation
of $S^1$  
on $T_{P_i}\M$ (repeated with multiplicity), for $i=0,\ldots, N$.
\begin{defin}
The \textbf{multiset of weights} $\W$ associated to the $S^1$-action on $(\M,\J)$ is the multiset 
$$\displaystyle\biguplus_{P_i\in \M^{S^1}}\{w_{i1},\ldots,w_{in}\}$$ and
the \textbf{multiset of positive weights} $\W_+\subset \W$ (resp. \textbf{negative weights} $\W_-\subset \W$) is the multiset 
$$\displaystyle\biguplus_{P_i\in \M^{S^1}}\{w_{ik}\mid w_{ik}>0\}\quad \mbox{(resp.} \displaystyle\biguplus_{P_i\in \M^{S^1}}\{w_{ik}\mid w_{ik}<0\}\mbox{)},$$ 
where $\biguplus$ denotes the union of multisets.
\end{defin}
Note that none of the elements of $\W$ is zero, since we are assuming the fixed points to be isolated. Hence, $\W=\W_+\cup \W_-$.

Our definition of multigraph associated to the $S^1$-action on $(\M,\J)$ relies on a crucial property of $\W$
which was first
proved by Hattori in the almost complex case \cite[Proposition 2.11]{Ha} (see also
\cite[Theorem 3.5]{L2} and \cite[Lemma 13]{PT}).
\begin{lemma}\label{pairs}
Let $(\M,\J)$ be an almost complex manifold equipped with a $\J$-preserving circle action with isolated fixed points.
Let $\W_+$ and $\W_-$ be the multisets of positive and negative weights.
Then $\W_+=-\W_-$. 
\end{lemma} 
In other words, $\W$ always contains pairs of integers of opposite signs.
Let us consider the index sets
$$
I:=\{ (i,k)\in\{0,\ldots,N\}\times \{1,\ldots,n\}\mid w_{ik}\in \W_+ \}
$$
and 
$$
J:=\{ (i,k)\in\{0,\ldots,N\}\times \{1,\ldots,n\}\mid w_{ik}\in \W_- \}.
$$
Since $\W_+=-\W_-$,
we can always choose a bijection between $\W_+$ and $\W_-$ such that the corresponding bijection 
$f:I\to J$ between the two index sets satisfies
\begin{equation}\label{cond action}
f(i,k)=(j,l) \quad \text{and} \quad w_{ik}=-w_{jl} .
\end{equation}
Hence, to every fixed point $P_i$, we are associating a fixed point $P_j$ such that one of the
weights at $P_i$ is $w_{ik}\in \W_+$ and one of the weights at $P_j$ is $w_{jl}=-w_{ik}\in \W_-$.
Let $\rho:I\to \{0,\ldots, N\}$ be defined by $\rho(i,k):=j$, where $j$ is such that $f(i,k)=(j,l)$.
For every choice of a bijective pairing $f:I \to J$ satisfying \eqref{cond action}, we can associate an
oriented 
multigraph $\Gamma=(V,E)$ to the $S^1$-action on $(\M,\J)$ as follows:
\begin{itemize}
\item The vertex set $V$ is the fixed point set.
\item The edge set is determined by the bijection between $\W_+$ and $\W_-$ chosen above. More precisely,
 for every $(i,k)\in I$ such that $f(i,k)=(j,l)$,
there corresponds an oriented edge $e_{ik}$ from $P_i$ to $P_j$.
The edge set $E$ is then $E=\{e_{ik}\mid (i,k)\in I\}$. 
\end{itemize}
Consequently, the elements of $E$ are in bijection with those of  $I$. Note that, in the notation above, if $e_{ik}$ is the edge associated to the pair $(i,k)\in I$, then $\ii(e_{ik})=P_i$ and $\te(e_{ik})=P_{\rho(i,k)}$.

We label the multigraph $\Gamma=(V,E)$ by the map
 $
\w\colon E \to \Z_{>0}
$ 
which, to each edge $e=e_{ik}$, associates  the weight $\w(e)=w_{ik}$ (which, by definition, is  always positive).
We call  this map the \emph{weight map}. 

Let us  denote by $\mathcal{W}$ the family of multigraphs
associated to the $S^1$-action on $\M$ labeled by the weight map, and denote its elements 
by  pairs $(\Gamma,\w)$.
\begin{lemma}\label{set of weights}
Let $(\Gamma,\w)$ be an element of $\mathcal{W}$. Then
\begin{itemize}
\item[1)] $(\Gamma,\w)$ determines $\W$.
More precisely, for every $P_i\in \M^{S^1}$, let
$\{w_{i1},\ldots,w_{in}\}$ be the multiset of weights of the $S^1$-representation on $T\M|_{P_i}$. Then
\begin{equation}\label{weight}
\{w_{i1},\ldots,w_{in}\}=\{\delta(P_i,e)\, \w(e),e\in E_{P_i}\setminus E^{\circlearrowleft}_{P_i}\}\cup \{\pm \w(e), e\in E^{\circlearrowleft}_{P_i}\}\;.
\end{equation}
\item[2)]
The magnitude is given by
\begin{equation}\label{mc}
\m(e)=\frac{\cc^{S^1}_1\left(\ii(e)\right)-\cc^{S^1}_1\left(\te(e)\right)}{\w(e)x}\quad\mbox{for all}\quad e\in E.
\end{equation}
\end{itemize}
\end{lemma}
\begin{proof}
The equality of multisets \eqref{weight} is an easy consequence of the definitions.
Moreover, \eqref{mc} comes from the definition of $\m(e)$, \eqref{weight} and the fact that $\cc_1^{S^1}(P_i)=(\sum_{j=1}^nw_{ij})x$
 for every $P_i\in \M^{S^1}$.
\end{proof}
Observe that for all $P\in V$, $\lambda(P):=\lvert E_{P,\te}\rvert$ is just the number of negative weights of the $S^1$-representation on $T_P\M$. Then 
$$\lvert \W_-\rvert=\sum_{P\in \M^{S^1}}\lambda(P)$$ 
which, by Lemma~\ref{pairs},
must be equal to
$$\lvert \W_+\rvert =\sum_{P\in \M^{S^1}}(n-\lambda(P)),$$ 
and so 
$$\lvert \W_+\rvert =\lvert \W_-\rvert =\frac{\lvert \W\rvert }{2}=\frac{(N+1)n}{2}.$$ 
From the definition of $E$ we then have
\begin{equation*}
\lvert E\rvert =\frac{(N+1)n}{2}.
\end{equation*}

We will now explore the properties satisfied by the magnitude $\m$ associated to $(\Gamma,\w)\in \mathcal{W}$.

\begin{rmk}
Suppose that there exists an edge $e\in E$ such that the points $\ii(e)$ and $\te(e)$ are the $S^1$-fixed points of an $S^1$-invariant sphere $S^2_e$ embedded
in $\M$, where $\w(e)$ (resp. $-\w(e)$) is the weight of the $S^1$ representation on $T_{\ii(e)}S_e^2$ (resp. $T_{\te(e)}S_e^2$). Then $\m(e)$ is just  the integral
of the first Chern class $\cc_1$ on the sphere $S_e^2$. In fact, by \eqref{mc} and the ABBV localization formula (Corollary~\ref{abbv discrete}), 
we have 
$$
\m(e)=\frac{\cc^{S^1}_1(\ii(e))-\cc^{S^1}_1(\te(e))}{\w(e)x}=\int_{S_e^2} \iota^*(\cc_1^{S^1})=\int_{S_e^2}\iota^*(\cc_1),
$$
where $\iota: S^2 \to \M$ is the inclusion map.
\end{rmk}
Let $(\Gamma,\w)$ be an element of $\mathcal{W}$, and let us fix an ordering
$e_1,\ldots,e_{\lvert E\rvert}$  of the edge set $E$. Then, by definition of $\m$,
we get
\begin{equation}\label{system}
\sum_{m=1}^{\lvert E\rvert}\big(\delta\left(\ii(e_h),e_m\right)-\delta\left(\te(e_h),e_m\right) \big)\w(e_m)-\w(e_h)\m(e_h)=0\;
\end{equation}
for all $h=1,\ldots,\lvert E \rvert$.

Denoting by $\w(E)$ the vector $\left(\w(e_1),\ldots,\w(e_{\lvert E \rvert})\right)$,
and by $\m(E)$ the vector $\left(\m(e_1),\ldots,\m(e_{\lvert E \rvert})\right)$,
from  \eqref{matrix} and \eqref{system} we obtain 
the  homogeneous system
\begin{equation}\label{system2}
\Big(A(\Gamma)-\diag(\m(E))\Big)\cdot \w(E)^t=0\;,
\end{equation}
where $\diag(\m(E))$ is the $(\lvert E\rvert \times \lvert E\rvert)$-diagonal matrix $\diag\left(\m(e_1),\ldots,\m(e_{\lvert E \rvert})\right)$.

Let us study  \eqref{system2} in more detail.
Suppose that $e=e_j$ is a cycle, i.e. that  $\ii(e)=\te(e)$. Then the $j-$th row and column of
$A(\Gamma)-\diag\left(\m(E)\right)$  are zero vectors. Let $\Gamma^\prime$ be the graph obtained from $\Gamma$ by deleting the cycles $e\in E^{\circlearrowleft}$,
i.e. $\Gamma^\prime=(V,E^\prime)$, where $E^\prime=E\setminus E^{\circlearrowleft}$, and
 $\Gamma^\prime=\Gamma_1\cup\cdots\cup\Gamma_l$ is the decomposition of $\Gamma^\prime$ into its connected components $\Gamma_i$.
Pick the following ordering of the edges of $\Gamma$:
\begin{equation}\label{ordering of E}
E=(e^1_{1},\ldots,e^1_{\lvert E_1\rvert},e^2_1,\ldots,e^2_{\lvert E_2\rvert },\ldots,e^l_{1},\ldots,e^l_{\lvert E_l\rvert},e_1^{\circlearrowleft},\ldots,e^{\circlearrowleft}_{\lvert E^{\circlearrowleft}\rvert })\;,
\end{equation}
where, writing  $\Gamma_i=(V_i,E_i)$ for every $i=1,\ldots,l$, we have $e^i_j\in E_i$ for all $j=1,\ldots,\lvert E_i\rvert$, and $e_j^{\circlearrowleft}\in E^{\circlearrowleft}$ for every $j=1,\ldots,\lvert E^{\circlearrowleft}\rvert$.
It is easy to see that $A(\Gamma)$ is a block diagonal matrix of the form
$$
A(\Gamma)=
\left(
\begin{array}{cccc}
A(\Gamma_1) &  & 0 & 0\\
\vdots & \ddots & \vdots &\vdots\\
 &  & A(\Gamma_l) & 0 \\
0 & \cdots & 0 & \mathbf{0}\\ 
\end{array}
\right),
$$
where the bottom right $(\lvert E^{\circlearrowleft}\rvert \times \lvert E^{\circlearrowleft}\rvert $)-zero matrix corresponds to the cycles.
For every $i=1\ldots,l$, let 
\begin{equation}\label{eq:Matricesi}
\ns \Big(A(\Gamma_i)-\diag(\m(E_i))\Big)
\end{equation}
be the null space of $A(\Gamma_i)-\diag\left(\m(E_i)\right)$, and let  $\Z_{>0}^{\lvert E_i\rvert }$ be the vectors in $\R^{\lvert E_i\rvert}$ with positive integer entries.
Since $\w(e)$ is a positive integer for every $e\in E$, we have the following result.

\begin{prop}\label{magnitude}
Let $(\Gamma,\w)$ be any multigraph associated to a fixed $S^1$-action on $(\M,\J)$, and let
$\Gamma^\prime=(V,E^\prime)$ be the graph obtained from $\Gamma=(V,E)$ by deleting its cycles,
i.e. $E^\prime=E\setminus E^{\circlearrowleft}$. Moreover,
let $\Gamma_1\cup\cdots\cup\Gamma_l$ be the decomposition of $\Gamma^\prime$ into
its connected components, and let  $A(\Gamma_i)$ be the matrices associated to $\Gamma_i$.
Then
\begin{equation*}
\Big( \ns \big(A(\Gamma_i)-\diag\left(\m(E_i)\right)\big) \Big)\cap \Z_{>0}^{\lvert E_i\rvert }\neq \emptyset\;,\quad \mbox{for every} \;\;\;i=1,\ldots,l,
\end{equation*} 
where $\m$ is the magnitude associated to $\Gamma$. 
In particular, 
\begin{equation*}
\det \Big(A(\Gamma_i)-\diag(\m(E_i))\Big)=0\;, \quad \mbox{for every} \;\;\;i=1,\ldots,l.
\end{equation*}
\end{prop}

Notice that the magnitude $\m\colon E \to \Q$ associated to $(\Gamma,\w)\in\mathcal{W}$ clearly depends on the labeled multigraph chosen. However, the sum of the elements in the image of $\m$ does not. Indeed, we have the following result.
\begin{prop}\label{thm2}
For any $(\Gamma,\w)\in\mathcal{W}$ the associated magnitude $\m\colon E\to \Q$ satisfies
$$
\sum_{e\in E} \m(e) = \int_\M \cc_1 \cc_{n-1}.
$$
\end{prop}
\begin{proof}
By the ABBV localization formula (Corollary~\ref{abbv discrete}) and \eqref{mc} we have
\begin{align*}
 \int_\M \cc_1 \cc_{n-1} & = \int_\M \cc^{S^1}_1 \cc^{S^1}_{n-1} = \sum_{i=0}^{N} \frac{\cc^{S^1}_1(P_i) \cc^{S^1}_{n-1}(P_i)}{e^{S^1}(\N_{P_i})} =\frac{1}{x} \sum_{i=0}^{N} \frac{\cc^{S^1}_1(P_i) \left(\sum_{l=1}^n \prod_{k\neq l} w_{ik}\right) }{\prod_{k=1}^n w_{ik}} \\ & =\frac{1}{x} \sum_{i=0}^{N}\sum_{k=1}^{n} \frac{\cc^{S^1}_1(P_i)}{w_{ik}} =  \frac{1}{x}\sum_{(i,k)\in I} \frac{\cc^{S^1}_1(P_i)}{w_{ik}}+\frac{1}{x}\sum_{(j,l)\in J} \frac{\cc^{S^1}_1(P_j)}{w_{jl}}=\\
& =  \sum_{(i,k)\in I} \frac{\cc^{S^1}_1(P_i)-\cc^{S^1}_1(P_{\rho(i,k)})}{w_{ik}x} =\sum_{e\in E}\frac{\cc^{S^1}_1(\ii(e))-\cc^{S^1}_1(\te(e))}{\w(e)x}=  \sum_{e\in E} \m(e).
\end{align*}
\end{proof}
Combining Theorem~\ref{thm1} and Proposition~\ref{thm2} we obtain the following result.
\begin{theorem}\label{sum of m}
Let $(\Gamma,\w)$ be an element of $\mathcal{W}$, and let $\m$ be the associated magnitude.
Then the sum of the magnitudes of the edges is an invariant of $\mathcal{W}$.
More precisely, let $n$ be the degree of $\Gamma$ and let $N_p$ be the number of vertices $P$ with $\lambda(P)=p$, for all $p=0,\ldots,n$.
Then
\begin{equation}\label{sum m}
\sum_{e\in E} \m(e) =\sum_{p=0}^n N_p[6p(p-1)+\frac{5n-3n^2}{2}].
\end{equation}
In particular, if $N_p=1$ for every $p=0,\ldots,n$, then
\begin{equation}\label{sum m min}
\sum_{e\in E} \m(e) =\frac{1}{2}n(n+1)^2\;.
\end{equation}
\end{theorem}
Let us now restrict our attention to a class of multigraphs for which the magnitude $\m\colon E\to \Q$ has integer values.

For every integer $l>1$, let $\M^{\Z_l}$ be the submanifold of $\M$ fixed by the subgroup $\Z_l\subset S^1$. Then we have the following lemma.
\begin{lemma}
There exists a labeled multigraph $(\Gamma,\w)\in\mathcal{W}$ such that, for every $e\in E$ with $\w(e)>1$, 
\begin{equation}\label{cc}
\ii(e)\quad\mbox{and}\quad\te(e)\quad\mbox{ lie in the same connected component of}\quad \M^{\Z_{\w(e)}}.
\end{equation}
\end{lemma}
\begin{proof} For every fixed point $P$, let $w_1,\ldots,w_n$ be the weights of the $S^1$-representation on $T_P\M$. For every
$w_i>1$, the set 
$\M^{\Z_{w_{i}}}$ is an almost complex submanifold of $\M$ of dimension greater than zero. Let $N$ be the connected component
of $\M^{\Z_{w_{i}}}$ containing $P$. Then, by
applying Lemma~\ref{pairs} to $N$, we know that  there exists a fixed point $Q\in N^{S^1}$ with one weight equal to $-w_i$.
Since $N^{S^1}\subset \M^{S^1}$, the conclusion follows immediately.
\end{proof}
\begin{defin}\label{img}
A labeled multigraph $(\Gamma,\w)\in \mathcal{W}$ satisfying $\eqref{cc}$ is called an \textbf{integral multigraph}.
We denote by $\;\;\mathcal{I}\;\;$ the family of integral multigraphs associated to a fixed
$S^1$-action on $(\M,\J)$ and by $\mathcal{M}$  the family of integral multigraphs labeled by their corresponding magnitudes. The elements of $\mathcal{M}$ are denoted  by $(\Gamma,\m)$.
\end{defin}
The main property of integral multigraphs is that the corresponding magnitudes have integer values. This is an easy consequence of \cite[Lemma 2.6 ]{T1},
which is stated for symplectic manifolds, but also holds for general  almost complex manifolds.
\begin{lemma}\label{modulo}
Let $P$ and $Q$ be fixed points of the $S^1$-action which lie in the same connected component of $\M^{\Z_l}$, for some integer $l>1$.
Then the weights of the isotropy action of $S^1$ on $T_P\M$ agree with the ones on $T_Q\M$ modulo $l$. 
\end{lemma}
So we have the following result.
\begin{prop}
For every $(\Gamma,\m)\in \mathcal{M}$, the magnitude $\m$ has integer values, i.e. $\m\colon E\to \Z$.
\end{prop}
\begin{proof}
By definition of integral multigraph, for every edge $e$ such that $\w(e)>1$, the endpoints $\ii(e)$ and $\te(e)$ of $e$ lie in the same connected component of $\M^{\Z_{\w(e)}}$.
Hence, by Lemma~\ref{modulo}, the weights of the isotropy representation of $S^1$ on $T_{\ii(e)}\M$ agree with the ones on $T_{\te(e)}\M$ modulo $\w(e)$.
The conclusion then follows immediately from the fact  that $\frac{\cc_1^{S^1}(P)}{x}$ is the sum of the weights of the $S^1$-representation on $T_P\M$,
for every $P\in \M^{S^1}$.
\end{proof}
Two important subsets of $\mathcal{I}$ are those given by  positive integral multigraphs and non-negative integral multigraphs.
\begin{defin}\label{P0+}
An integral multigraph $(\Gamma,\w)\in \mathcal{I}$ is called a \textbf{non-negative (resp. positive)} if the corresponding 
magnitude satisfies $\m\colon E\to \Z_{\geq 0}$ (resp. $\m\colon E\to \Z_{> 0}$).  

We denote by $\mathcal{I}_{\geq 0}$ (resp. $\mathcal{I}_{> 0}$) the family of non-negative (resp.  positive) multigraphs,
and by $\mathcal{M}_{\geq 0}$ (resp. $\mathcal{M}_{> 0}$) the family of non-negative (resp. positive) multigraphs labeled by their corresponding magnitudes.  

We say that the $S^1$-space $(\M,\J,S^1)$ satisfies {\bf property} ${\bf (\mathcal{P}_0^+)}$ if the corresponding family of non-negative multigraphs   $\mathcal{I}_{\geq 0}$ is nonempty.
\end{defin}

In Sections~\ref{algo} and \ref{mnfp} we will investigate in more detail conditions on the $S^1$-action that guarantee that it
can be represented by a non-negative or positive  multigraph (see in particular Remark~\ref{monotone..} and Section~\ref{refinement}).  
\begin{exm}
Let us consider the 
circle action on $\M=S^2 \times S^2$ which rotates one sphere with speed $a$ and the other with speed $b$, where $a,b\in \mathbb{Z}_{>0}$ are relatively prime. There are four fixed points 
$$P_0=(S,S),\,\, P_1=(S,N),\,\, P_2=(N,S) \quad \text{and} \quad P_3=(N,N),$$
where $S$ and $N$ are the north and south poles of $S^2$. There are two spheres fixed by the action of the subgroup $\Z_a$ (namely $S^2\times \{S\}$ and $S^2\times \{N\}$) and two spheres fixed by the action of $\Z_b$ (namely $\{S\}\times S^2$ and $\{N\} \times S^2$). The isotropy weights at each fixed point and the number $\lambda_i$ of negative weights at $P_i$ are  
\begin{center}
\begin{tabular}{|l|| l|l|}
\hline
$P_0:$ & $(w_{01},w_{02})=(a,b)$ & $\lambda_0=0$ \\ \hline
$P_1:$ & $(w_{11},w_{12})=(-b,a)$ & $\lambda_1=1$ \\ \hline
$P_2:$ & $(w_{21},w_{22})=(-a,b)$ & $\lambda_2=1$ \\ \hline
$P_3:$ & $(w_{31},w_{32})=(-b,-a)$ & $\lambda_3=2$ \\ \hline
\end{tabular}\;.
\end{center} 
The values of the equivariant first Chern class $\cc_1^{S^1}(P_i)$ at each fixed point are
\begin{center}
\begin{tabular}{|l|| c|}
\hline
$\cc_1^{S^1}(P_0)=$ & $(a+b)x$ \\ \hline
$\cc_1^{S^1}(P_1)=$ & $(a-b)x$ \\ \hline
$\cc_1^{S^1}(P_2)=$ & $(b-a)x$ \\ \hline
$\cc_1^{S^1}(P_3)=$ & $-(a+b)x$ \\ \hline
\end{tabular}
\end{center}
and the multisets $\W_+,\W_-,I$ and $J$ are
\begin{center}
\begin{tabular}{|c|| l|}
\hline
$\W_+$ & $\{a,a,b,b\}$ \\ \hline
$\W_-$ & $\{-a,-a,-b,-b\}$ \\ \hline
$I$ & $\{(0,1),(0,2),(1,2),(2,2)\}$ \\ \hline
$J$ & $\{ (1,1),(2,1),(3,1),(3,2)\}$ \\ \hline
\end{tabular}\;.
\end{center} 
There are only four possible pairings $f_l:I \to J$, $l=1,\ldots,4$, as in \eqref{cond action}
\begin{center}
\begin{tabular}{|l|| l| c |l|| l| c |l|| l|c |l|| l|}\cline{1-2}\cline{4-5}\cline{7-8}\cline{10-11}
$f_1(0,1)=$ &$(2,1)$ & & $f_2(0,1)=$ &$(3,2)$  & & $f_3(0,1)=$ &$(2,1)$ & & $f_4(0,1)=$ &$(3,2)$\\ \cline{1-2}\cline{4-5}\cline{7-8}\cline{10-11}
$f_1(0,2)=$ & $(3,1)$ & & $f_2(0,2)=$ &$(3,1)$ & & $f_3(0,2)=$ &$(1,1)$ & & $f_4(0,2)=$ &$(1,1)$\\ \cline{1-2}\cline{4-5}\cline{7-8}\cline{10-11}
$f_1(1,2)=$ & $(3,2)$ & & $f_2(1,2)=$ &$(2,1)$ & & $f_3(1,2)=$ &$(3,2)$ & & $f_4(1,2)=$ &$(2,1)$\\ \cline{1-2}\cline{4-5}\cline{7-8}\cline{10-11}
$f_1(2,2)=$ & $(1,1)$ & & $f_2(2,2)=$ &$(1,1)$ & & $f_3(2,2)=$ &$(3,1)$ & & $f_4(2,2)=$ &$(3,1)$\\ \cline{1-2}\cline{4-5}\cline{7-8}\cline{10-11}
\end{tabular}
\end{center}
yielding the four graphs $\Gamma_l$ in Figure~\ref{Graphs}.
\begin{figure}[h!]
\!\!\!\includegraphics[scale=0.75] {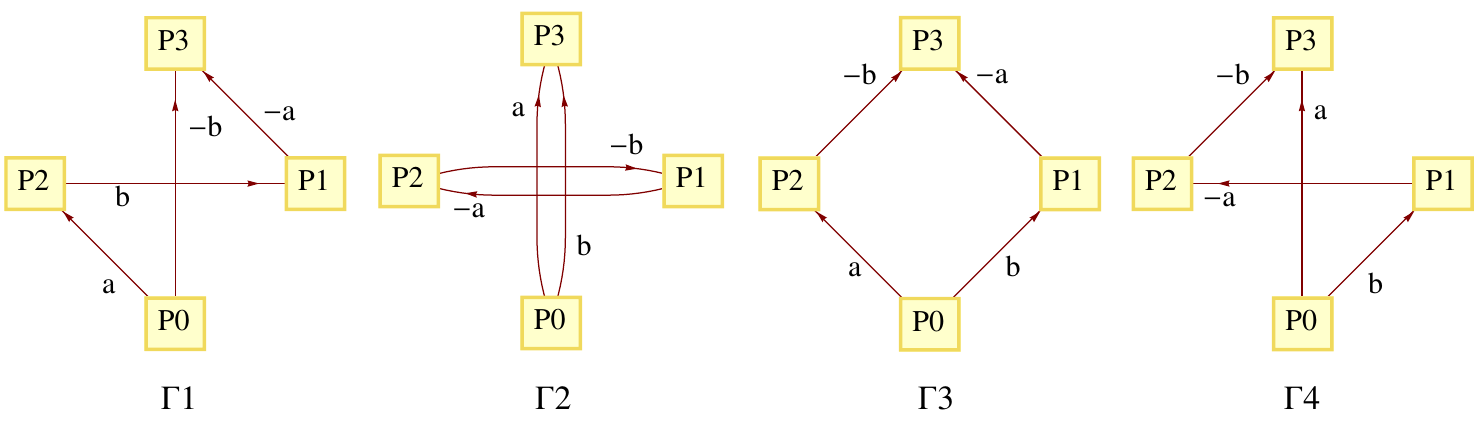}
\caption{Circle actions on $S^2\times S^2$}
\label{Graphs}
\end{figure}

Note that, if $a$ and $b$ are greater than $1$, the only integral multigraph is $\Gamma_3$ which is a positive integral multigraph. Indeed,
writing $e_1:=e_{01}$, $e_2:=e_{02}$, $e_3:=e_{12}$, $e_4:=e_{22}$, the associated magnitudes $\m\colon E \to \Q$ are
\begin{center}
\begin{tabular}{|c|| c| c |c|| c| c |c|| c| c |c|| c| }  

\cline{1-2}\cline{4-5}\cline{7-8}\cline{10-11}
$f_1$ & &  &  $f_2$ & & & $f_3$  &  & &  $f_4$  &  
\\ \cline{1-2}\cline{4-5}\cline{7-8}\cline{10-11}

$\m(e_{1})$ & $ 2 $ & & 
$\m(e_{1})$ &  $ \frac{2(a+b)}{a} $ & &
$\m(e_{1})$ & $ 2 $ & &
$\m(e_{1})$ & $ \frac{2(a+b)}{a}$
\\ \cline{1-2}\cline{4-5}\cline{7-8}\cline{10-11}

$\m(e_{2})$ &  $ \frac{2(a+b)}{b}$  &  &
$\m(e_{2})$ & $ \frac{2(a+b)}{b} $ & &
$\m(e_{2})$ & $ 2 $ & & 
$\m(e_{2})$ & $ 2 $ 
\\ \cline{1-2}\cline{4-5}\cline{7-8}\cline{10-11}

$\m(e_{3})$ &   $ 2$  & &
$\m(e_{3})$ &  $\frac{2(a-b)}{a} $ & &
$\m(e_{3})$ & $ 2 $ & &
$\m(e_{3})$ &  $ \frac{2(a-b)}{a} $ 
\\ \cline{1-2}\cline{4-5}\cline{7-8}\cline{10-11}

$\m(e_{4})$  &  $ \frac{2(b-a)}{b}$  & & 
$\m(e_{4})$ &  $ \frac{2(b-a)}{b} $ & & 
$\m(e_{4})$ &  $ 2 $ &  &
$\m(e_{4})$ &  $ 2 $  
\\ \cline{1-2}\cline{4-5}\cline{7-8}\cline{10-11}
\end{tabular}\;.
\end{center}
Note that  $\sum_{e\in E}\m(e)=8$ in all cases. On the other hand, the numbers $N_p$ of fixed points with $p$ negative weights are $N_0=N_2=1$ and $N_1=2$ and so, by Theorem~\ref{thm1}, we obtain
$$
\int_\M \cc_1 \cc_{n-1}= \sum_{p=0}^2 N_p[6p(p-1)-1]=-1-2+11=8,
$$
confirming the result in Proposition~\ref{thm2}.

Finally, for $\Gamma_i=(V,E_i)$, we obtain the matrices $A(\Gamma_i)$ given by 
\begin{align*}
A(\Gamma_1)=
\left(
\begin{array}{rrrr}
2 & 1 & 0 & -1\\
1 & 2 & 1 & 0\\
0 & 1 & 2 & -1\\
-1 & 0 & -1 & 2
\end{array}
\right)\quad\quad
&
A(\Gamma_2)=
\left(
\begin{array}{rrrr}
2 & 2 & 0 & 0\\
2 & 2 & 0 & 0\\
0 & 0 & 2 & -2\\
0 & 0 & -2 & 2\\ 
\end{array}
\right)\\
\\
A(\Gamma_3)=
\left(
\begin{array}{rrrr}
2 & 1 & 0 & -1\\
1 & 2 & -1 & 0\\
0 & -1 & 2 & 1\\
-1& 0 & 1 & 2\\ 
\end{array}
\right)\quad\quad
&
A(\Gamma_4)=
\left(
\begin{array}{rrrr}
2 & 1 & 0 & 1\\
1 & 2& -1 & 0\\
0 & -1 & 2 & -1\\
1& 0 & -1 & 2\\ 
\end{array}
\right).\\
\end{align*}
\end{exm}

\section{An algorithm to determine linear relations among the weights 
}\label{algo}
After having introduced the necessary definitions in  Section~\ref{dmg} we can rewrite Question~\ref{q3}
in the following way.
\begin{question}\label{q5}
Let $n,N$ and $\lambda_i$, for $i=0,\ldots,N$, be positive integers, where $0\leq\lambda_i\leq n$ for all $i$.
Does there exist a compact almost complex  manifold 
$(\M,\J)$ of dimension $2n$ with a $\J$-preserving $S^1$-action  with isolated
fixed points $P_0,\ldots,P_N$, such that the number of negative weights
at $P_i$ is $\lambda_i$?
If so, can we determine the corresponding family $\mathcal{I}$ of integral labeled multigraphs?
\end{question}
Instead of trying to determine $\mathcal{I}$, it is more convenient
to determine $\mathcal{M}$,  the family of integral multigraphs labeled by their magnitudes $\m$.
In fact, by Theorem~\ref{sum of m}, the sum of the magnitudes only depends on $n,N$ and $\lambda_i$ and so Proposition~\ref{magnitude} gives linear relations among the weights.

Notice that if the answer to Question~\ref{q5} is affirmative, the first necessary condition on $n,N$ and the $\lambda_i$s is a direct consequence of Lemma~\ref{pairs}.
Indeed, since $\lvert \W_+\rvert =\lvert \W_-\rvert =\displaystyle\frac{\lvert \W\rvert }{2}$, we must have
\begin{equation}\label{lambdas}
\sum_{i=0}^N(n-\lambda_i)=\sum_{i=0}^N\lambda_i=\frac{(N+1)n}{2}.
\end{equation}
From now on we will assume that $n$, $N$ and the $\lambda_i$s  satisfy \eqref{lambdas}.

The linear relations among the weights can be determined as follows.
First, let us define two families $\mathcal{N}$
and $\mathcal{N}_{\geq 0}$ of labeled multigraphs associated to the integers $n$, $N$ and $\lambda_i$, for $i=0,\ldots,N$. For that, consider the sets
$$I^\prime=\left\{(i,k)\in\{0,\ldots,N\}\times \{1,\ldots,n\}\mid k\le n-\lambda_i \right\}$$ 
and 
$$J^\prime=\left\{(i,k)\in\{0,\ldots,N\}\times \{1,\ldots,n\}\mid k\le \lambda_i \right\}.$$
Observe that by \eqref{lambdas} we have 
$\lvert I^\prime \rvert=\lvert J^\prime \rvert$. Let $f\colon I^\prime \to J^\prime$ be a bijection between the two sets. Using $f$ we can construct a multigraph 
 $\Gamma=(V,E)$, where
\begin{itemize}
\item $V$ is a set of $N+1$ vertices $P_0,\ldots,P_N$,
\item for every $(i,k)\in I^\prime$ with $f(i,k)=(j,l)$, one considers 
an oriented edge  $e_{ik}$ from $P_i$ to $P_j$ so that the edge set is  $E=\{e_{ik}\mid (i,k)\in I^\prime\}$. 
\end{itemize}
Let
$\Gamma^\prime=(V,E^\prime)$ be the graph obtained from $\Gamma=(V,E)$ by deleting its cycles,
i.e. $E^\prime=E\setminus E^{\circlearrowleft}$. Moreover,
let $\Gamma_1\cup\cdots\cup\Gamma_l$ be the decomposition of $\Gamma^\prime$ into
its connected components, and let  $A(\Gamma_i)$ be the matrices associated to $\Gamma_i$
as in \eqref{matrix}.
Consider an ordering on the edges as in \eqref{ordering of E} and let $\mathcal{N}_{\Gamma}$ be the family of maps
$\n\colon E\to \Z$ satisfying
\begin{eqnarray}
\label{cond1}
\displaystyle\sum_{e\in E} \n(e) =\sum_{p=0}^n N_p[6p(p-1)+\frac{5n-3n^2}{2}]\\
\label{cond2}
\Big(  \ns\left(A(\Gamma_i)-\diag\left(\n(E_i)\right)\right)\Big)\cap \Z_{>0}^{\lvert E_i\rvert }\neq \emptyset
\end{eqnarray}
for every $i=1,\ldots,l$.
Then we define  $\mathcal{N}$ as the set of pairs $(\Gamma,\n)$, where $\Gamma$ is associated to a bijection $f:I^\prime \to J^\prime$  as above and $\n\in \mathcal{N}_{\Gamma}$. Moreover,
we denote by $\mathcal{N}_{\geq 0}$ the subset of $\mathcal{N}$ 
given by pairs $(\Gamma,\n)\in \mathcal{N}$ such that $\n\colon E\to \Z_{\geq 0}$. 

\begin{rmk}
Condition \eqref{cond2} implies that $\det \left(A(\Gamma_i)-\diag\left(\n(E_i)\right)\right)=0$ for every $i$, yielding  polynomial
equations in the $\n(e)$s. However, \eqref{cond2} is much stronger.
For example, let $B_h$ be the $h$-th row of $A(\Gamma_i)-\diag\left(\n(E_i)\right)$, for some $i$. Then
it is easy to see that $B_h$ must contain entries of opposite signs for all $h=1,\ldots,\lvert E_i \rvert$. 
More precisely $B_h$ looks like
$$(b_{h,1},\ldots,b_{h,h-1},2-\n(e_h),\ldots,b_{h,\lvert E_i\rvert}),$$ 
where $b_{h,k}$ is an integer;
so,  if $b_{h,k}\geqslant 0$ for all $k\neq h$, we need $\n(e_h)\geqslant 2$.
Taking linear  combinations of rows in such a way that all the constant coefficients have the same sign, one can get
lower bounds for linear combinations of the $\n(e)$s. However, these estimates do not seem to be optimal. 
It would be interesting to know whether one can get better restrictions on the set of the $\n(e)$s satisfying \eqref{cond2}. 
\end{rmk}

Using the above notation we have the following result.
\begin{theorem}\label{NC}
Let $n$, $N$ and $\lambda_i$, for $i=0,\ldots,N$, be non negative integers satisfying $0\leq\lambda_i\leq n$ for $0\leq i\leq N$ as well as \eqref{lambdas}. 

If there exists a compact almost complex manifold 
$(\M, \J)$ of dimension $2n$ and a $\J$-preserving $S^1$-action on $\M$  with isolated
fixed points $P_0,\ldots,P_N$, such that the number of negative weights
at $P_i$ is $\lambda_i$, then
$\mathcal{M}\subset\mathcal{N}$ and $\mathcal{M}_{\geq 0}\subset \mathcal{N}_{\geq 0}$.

Moreover, for every $(\Gamma,\w)\in \mathcal{I}$ (resp. $\mathcal{I}_{\geq 0}$) there exists $(\Gamma,\n)\in \mathcal{N}$ (resp. $\mathcal{N}_{\geq 0}$) such that
\begin{equation*}
\w(E)\in  \Big(\ns(A\left(\Gamma)-\diag\left(\n(E)\right)\right)\Big)\cap \Z_{>0}^{\lvert E\rvert}.
\end{equation*}
\end{theorem}
\begin{proof}
Let $\{w_{i1},\ldots,w_{in}\}$ be the set of weights of the $S^1$-representation on $T_{P_i}\M$, for all $i=0,\ldots,N$,
and let $\W$, $\W_+$ and $\W_-$ be the multisets defined in Section~\ref{dmg}. We can order the weights in
such a way that $I^\prime$ and $J^\prime$ are the index sets of $\W_+$ and $\W_-$, i.e. $I^\prime=\{(i,k)\mid w_{ik}\in \W_+\}$
and $J^\prime=\{(i,k)\mid w_{ik}\in \W_-\}$. Then the conclusion follows from the definitions along with Proposition~\ref{magnitude} and Theorem~\ref{sum m}.
\end{proof}

Hence, to look for linear relations among the isotropy weights, we have to take every possible multigraph obtained from a bijection $f:I^\prime \to J^\prime$ and look for possible functions $\n:E\to \Z$ satisfying \eqref{cond1} and \eqref{cond2}, which would correspond to the magnitudes of the multigraphs. To make this task more efficient we point out a few facts that, in some cases, imply that the number of such functions $\n$ is finite. Moreover, when we cannot assure that this number is finite we can impose additional conditions on the action, making this number finite at the cost of restricting the class of circle actions considered.
Before doing so, we consider the following technical lemma.
\begin{lemma}\label{l dividing}
Let $(\M,\J)$ be a compact almost complex manifold with a $\J$-preserving $S^1$-action with isolated fixed points.
Let $(\Gamma,\w)\in \mathcal{I}$ be an integral multigraph associated to this action. Then, for every $w\in H_{S^1}^2(\M;\Z)$, we have
$$
\frac{w(\ii(e))-w(\te(e))}{\w(e)x} \in \Z \quad\mbox{for every edge $e\in E$}.
$$
\end{lemma}
\begin{proof}
As remarked in Section~\ref{eclb}, to every element $w\in H_{S^1}^2(\M;\Z)$, we can associate an equivariant complex line bundle
$\LL^{S^1}$ whose first equivariant Chern class is $w$ and such that $\LL^{S^1}(P)=t^{w(P)}$ for every fixed point $P$.
Let $e$ be an edge of $(\Gamma,\w)\in \mathcal{I}$ and suppose that $\w(e)>1$. By definition of integral multigraph,
$\ii(e)$ and $\te(e)$ are in the same connected component $N$ of $\M^{\Z_{\w(e)}}$, the submanifold of $\M$ fixed by $\Z_{\w(e)}$.
Thus the $S^1$-modules $\LL^{S^1}\left(\ii(e)\right)$ and $\LL^{S^1}\left(\te(e)\right)$ are equivalent as $\Z_{\w(e)}$-modules, and the conclusion follows immediately.
\end{proof}

As a result of this Lemma, for every equivariant class $w\in H^2_{S^1}(\M;\Z)$ that ``divides" the first equivariant Chern class $\ce$, i.e.  a class $w$ satisfying 
\begin{equation}\label{eq:C}
\ce=Cw
\end{equation}
for some constant $C\in \Z\setminus \{0\}$, we have
$$
\frac{w(\ii(e))-w(\te(e))}{\w(e)}=\frac{\m(e)}{C}\in \Z 
$$
for every edge $e\in E$. Hence, by \eqref{sum m}, we conclude that this constant $C$ must divide 
\begin{equation}\label{RHS}
\sum_{p=0}^n N_p[6p(p-1)+\frac{5n-3n^2}{2}].
\end{equation}
In many cases 
the largest positive integer $C$ satisfying this property can be explicitly  written in terms of the weights at the fixed points of index $0$ and $2$.
Indeed, 
assume now, and for the remaining of the section, that $(\M,\omega)$ is a compact symplectic manifold
with a compatible almost complex structure $\J$, and that the $S^1$-actions considered are Hamiltonian and preserve $\J$.

Let $P_0$ be  the fixed point of index 0, and $P_1^1,\ldots, P_1^k$ the fixed points of index $2$.
Assume that there exist degree-$2$ generators of $H_{S^1}^2(\M;\Z)$, $\tau_1^1,\ldots,\tau_1^k$, satisfying
\begin{equation}\label{canonical 2}
\tau_1^i(P)=0\quad\mbox{ for all}\quad P\in \M^{S^1}\setminus\{P_1^i\}\;\;\;\mbox{ such that }\lambda(P)\leq 1,
\end{equation}
which we call \emph{canonical classes}. These classes do not always exist, but, if they do exist, they are unique (see \cite[Lemma 2.7]{GT}).
Then it is easy to see that 
\begin{equation*}
\ce=\sum_{i=1}^m\alpha_i\tau_1^i+\beta,
\end{equation*}
for some $1\leq m\leq k$, with 
\begin{equation*}
\alpha_i=\displaystyle\frac{\ce(P_1^i)-\ce(P_0)}{\Lambda_{P_1^i}^-}\in \Z\setminus \{0\}\quad\text{and}\quad\beta=\ce(P_0).
\end{equation*}
Indeed, from Lemma 3.8 in \cite{T1} we have that $r(\cc_1^{S^1})=\cc_1\neq 0$, thus implying that $m\geq 1$; the explicit
expression of the $\alpha_i$s and $\beta$ is an easy consequence of the definitions.
So the largest positive constant $C$ dividing \eqref{RHS} for which $\ce/C$ is still an integer class
is 
\begin{equation*}
C=\gcd\{\lvert \alpha_1\rvert ,\ldots,\lvert \alpha_m\rvert \}.
\end{equation*}

Let us now make the additional assumption that there is a multigraph  with an edge   $e_i$ from  $P_0$ to   $P_1^i$, for every $i=1,\ldots,m$. Then $\m(e_i)=\alpha_i$ for every $i=1,\ldots,m$
and so 
$$C=\gcd\{\lvert \m(e_1)\rvert ,\ldots,\lvert \m(e_m)\rvert \}.$$ 
Note that we can always take such a multigraph when 
there exists a sphere connecting $P_0$ to $P_1^i$ fixed by some subgroup of $S^1$, for every $i=1,\ldots,m$.

\begin{exm}\label{flag}
Consider the standard $T^3$-action on $\C^3$ given by
$$
(\xi_1,\xi_2,\xi_3)\cdot (z_1,z_2,z_3)=(\xi_1 z_1, \xi_2 z_2, \xi_3 z_3)\;.
$$
This action descends to an action
 on the complete flag manifold $\mathcal{F}l(\C^3)$, where the diagonal circle $S^1=\{(\xi,\xi,\xi)\in T^3\}$
 acts trivially. Take  the action of the quotient group $T^3/S^1\simeq T^2$. This action has $6$ fixed points, given 
by flags in the coordinate lines of  $\C^3$. These are indexed by permutations $\sigma \in  S_3$ on $3$ letters, where  the fixed point corresponding to $\sigma$ is given by
$$
\langle 0 \rangle \subset \langle f_{\sigma (1)} \rangle  \subset \langle  f_{\sigma (1)},  f_{\sigma (2)} \rangle \subset  \langle f_{\sigma (1)}, f_{\sigma (2)}, f_{\sigma (3)} \rangle  = \C^3,
$$
the brackets indicate the span of the vectors, and $\{f_1,f_2,f_3\}$ is the standard basis of $\C^3$. By taking a generic 
circle $S^1\subset T^3$ we obtain a circle action with one fixed point of index $0$, two fixed points of index $2$, $P_1^1$ and $P_1^2$, two fixed points of index $4$ and one fixed point of index $6$. The weights 
at the minimum are then $\{m,n,m+n\}$ and  the weights at the index-$2$ fixed  points are respectively $\{-m,n,m+n\}$ and $\{-n,m,m+n\}$, where
$m$ and $n$ are coprime integers which depend on the circle $S^1\subset T^2$ chosen.
Let $\tau_1^1$ and $\tau_1^2$ be generators of $H^2_{S^1}(\M;\Z)$. Then they can be chosen in such a way that they satisfy property \eqref{canonical 2},
and it is easy to check that $\ce=2(\tau_1^1+\tau_1^2)+2(m+n)$, yielding $C=2$.
\end{exm}

Taking the above remarks into consideration let us  see how to proceed in general.

\vspace{.3cm}
${\bf 1.}\,\,$  \emph{If none of the $\m(e)$s is negative} we have the following algorithm.

\begin{alg}\label{alg:1}
If none of  the $\m(e)$s is negative, then, for each possible multigraph $\Gamma$ with connected components $\Gamma_1,\ldots,\Gamma_l$, we look for partitions of
$$
\sum_{p=0}^n N_p[6p(p-1)+\frac{5n-3n^2}{2}]
$$
into $\lvert E^\prime\rvert=\lvert E\rvert - \lvert E^{\circlearrowleft} \rvert$ positive integers  $\m(e)$. Then we add the zeros $\m(e)=0$ whenever the edge $e$ is a cycle and choose, among  the resulting  sequences of $\lvert E \rvert$ nonegative numbers, those  for which
$$
\Big(  \ns\left(A(\Gamma_i)-\diag\left(\m(E_i)\right)\right)\Big)\cap \Z_{>0}^{\lvert E_i\rvert }\neq \emptyset
$$
for every $i=1,\dots,l$.
\end{alg}

When none of  the $\m(e)$s is negative we have an upper bound for any constant $C$ as in \eqref{eq:C}.
In fact, since $\frac{\m(e)}{C}$ is an integer for each edge $e\in E$, we have
$$\sum_{e\in E^\prime}\frac{\m(e)}{C}\geq \lvert E^\prime\rvert =\frac{n(N+1)}{2} -\lvert E^{\circlearrowleft}\rvert,$$
where $N+1$ is the total number of fixed points. Observe that
$|E^{\circlearrowleft}|<|E|$; in fact, since the action is Hamiltonian, none of the edges
starting at the minimum $P_0$ can be a cycle. So $|E'|>0$.
Moreover it's easy to see that $|E^{\circlearrowleft}|$ is bounded by $\sum_{i=0}^N\min\{\lambda_i,n-\lambda_i\}$.
By using \eqref{sum m} together with the fact that $C$ is a positive integer, we obtain
\begin{equation}\label{ub C}
1\leq C\leq \frac{2\sum_{p=0}^n N_p[6p(p-1)+\frac{5n-3n^2}{2}]}{n(N+1)-2\lvert E^{\circlearrowleft}\rvert}.
\end{equation}

\begin{exm} 
If $\dim(\M)=6$, $N_0=N_3=1$ and $N_1=N_2=2$, we have 
$$ \lvert E^{\circlearrowleft} \rvert \leq 4, \quad \sum_{e\in E}\m(e)=24\quad \text{and} \quad \lvert E \rvert=9, 
$$
and so $C=1$ or $2$ whenever $\lvert E^{\circlearrowleft} \rvert=0$, $C=1,2$ or $3$ when $\lvert E^{\circlearrowleft}\rvert=1$ or $2$, and  $C=1,2,3$ or $4$ if  $\lvert E^{\circlearrowleft} \rvert=3$ or $4$. Note that Example~\ref{flag} falls into this case.
\end{exm}

If we  can choose a basis of $H_{S^1}^2(\M;\Z)$ where the degree-$2$ generators 
are canonical classes, and there exists an edge between the point of index zero and each point of index $2$,
we can use the bound in \eqref{ub C} to improve Algorithm~\ref{alg:1}.

\begin{alg}\label{alg:2}
Assume that we can choose a basis of $H_{S^1}^2(\M;\Z)$ where the degree-$2$ generators 
are canonical classes. 
Suppose that there exists a multigraph $\Gamma$ such that none of the $\m(e)$s is negative and which has an edge $e_i$ from $P_0$ to $P_1^i$, for each $i$ such that 
$\ce(P_0)/x>\ce(P_1^i)/x$; let $e_1,\ldots,e_m$ be these edges.
Let $\Gamma_1,\ldots,\Gamma_l$ be the connected components of $\Gamma$. Then for each divisor $C$ of 
$$
\sum_{p=0}^n N_p[6p(p-1)+\frac{5n-3n^2}{2}],
$$
satisfying
$$
1\leq C\leq \frac{2\sum_{p=0}^n N_p[6p(p-1)+\frac{5n-3n^2}{2}]}{n(N+1)-2 \lvert E^{\circlearrowleft} \rvert},
$$
we look for partitions of
$$
\frac{1}{C}\sum_{p=0}^n N_p[6p(p-1)+\frac{5n-3n^2}{2}]
$$
into $\lvert E^\prime\rvert =\lvert E \rvert - \lvert E^{\circlearrowleft} \rvert$ positive integers, $\lu(e):=\frac{\m(e)}{C}$ with $\gcd\{\lu(e_1),\ldots,\lu(e_m)\}=1$. Then we consider the corresponding integers $\m(e)=C\cdot \lu(e)$, add the zeros $\m(e)=0$ whenever the edge $e$ is a cycle,
and choose, among the resulting sequences of $\lvert E \rvert$ numbers, the ones  for which
$$
\Big(  \ns\left(A(\Gamma_i)-\diag\left(\m(E_i)\right)\right)\Big)\cap \Z_{>0}^{\lvert E_i\rvert }\neq \emptyset,
$$
for every $i=1,\dots,l$.
\end{alg}

\vspace{0.5cm}

\begin{rmk}\label{monotone..}{\bf (Non-negative multigraphs)}
One cannot expect in general that non-negative multigraphs exist. Indeed, if the number of fixed points is big, meaning that the $N_p$s are greater than $1$,
the right hand side of \eqref{sum m}, i.e. $\int_\M \cc_1\cc_{n-1}$, can be \emph{negative}. Consequently,
requiring all the $\m(e)$s to be non-negative is a very restrictive  assumption. For example, if $\dim \M=4$, we have $N_0=N_2=1$ and $\int_\M \cc_1^2=10-N_1$.

However, if given an $S^1$-action we can choose a multigraph such that for every edge $e\in E$ there exists an isotropy $2$-sphere $S_e^2$ having
$\ii(e)$ and $\te(e)$ as south and north pole (e.g. for GKM manifolds, see Section~\ref{refinement} {\bf v)}), 
then this multigraph is positive if $\cc_1$ is positive on each of these spheres, which happens for example 
\begin{itemize}
\item 
If $(\M,\omega)$ is \emph{monotone}, i.e. the symplectic form satisfies  $\cc_1=C[\omega]$ for some positive constant $C$.
\item 
More generally, when $(\M,\omega)$ is \emph{symplectic Fano}, i.e. if,
given an almost complex structure $\J$ compatible with $\omega$, we have $\cc_1(A)>0$ for every $A\in H_2(\M)$ which can be represented by a $\J$-holomorphic curve.
\item If there exist classes $\sigma^1,\ldots,\sigma^k$ in $H^2(\M;\Z)$ such that
$\cc_1=\sum_{i=1}^k\beta_i\sigma^i$, where $\beta_i\in \Z_{>0}$ for every $i$, and for every
isotropy sphere $S_e^2$, there exists $j$ such that $\int_{S_e^2} \sigma^j>0$ and $\int_{S_e^2} \sigma^i\geq 0$ for $j\neq i$.
\end{itemize}
\end{rmk}

Let us now see how to obtain a finite algorithm from \eqref{sum m} even when we cannot assume the $\m(e)$s to be non-negative.
\\$\;$\\

${\bf 2.}\,\,$
\emph{If the $\m(e)$s are also negative} we proceed as follows. Consider the map $\Psi$ that to each fixed point $P$ associates the sum of the weights at $P$,
\begin{eqnarray}
 \Psi\colon &  \M^{S^1} &\longrightarrow  \Z \nonumber\\
          & P               &\longmapsto               \sum_{i=1}^n w_{iP}=\frac{\ce(P)}{x},    \label{Psi}  
 \end{eqnarray}
 and let 
 \begin{equation*}
 D=\max\{\lvert \Psi(P)\rvert \}_{P\in \M^{S^1}}. 
 \end{equation*}
 Note that $D> 0$ since the fact that the  action is Hamiltonian yields $\Psi(P_0)>0$.
 Then 
 $$\lvert \m(e)\rvert \leq \left\lvert \frac{\ce(\ii(e))-\ce(\te(e))}{x} \right\rvert \leq 2 D$$
 and  so we can still construct an algorithm to obtain necessary conditions on the isotropy weights, now with a prescribed bound on $\Psi$.
 \begin{alg}\label{alg:3}  Fixing  a positive integer $D$, for each multigraph $\Gamma$ with connected components $\Gamma_1,\ldots,\Gamma_l$, we look for partitions of 
\begin{equation*}
\sum_{p=0}^n N_p[6p(p-1)+\frac{5n-3n^2}{2}] 
\end{equation*} 
into $\lvert E \rvert$ integers, $\m(e)$, with $\lvert \m(e)\rvert \leq 2D$ and choose those for which
$$
\Big(  \ns\left(A(\Gamma_i)-\diag\left(\m(E_i)\right)\right)\Big)\cap \Z_{>0}^{\lvert E_i\rvert }\neq \emptyset
$$
for every $i=1,\dots,l$.

Note that some of the above integers may be zero depending on the existence of cycles in $\Gamma$.
\end{alg}

Moreover, suppose that there exists an edge $e_1^i$ between the point of index $0$ and every index-$2$ fixed point $P_1^i$
such that $\Psi(P_0)\neq \Psi(P_1^i)$; call these edges $e_1,\ldots,e_m$.
Suppose that there exists canonical classes $\tau_1^i$ satisfying
\eqref{canonical 2}; then we can again improve this algorithm.

\begin{alg}\label{alg:4}
Fixing  a positive integer $D$, for each multigraph $\Gamma$ with connected components $\Gamma_1,\ldots,\Gamma_l$ and for each integer $C\in [1,D]\cap \Z$ 
dividing 
$$
\sum_{p=0}^n N_p[6p(p-1)+\frac{5n-3n^2}{2}],
$$
we look for partitions of 
\begin{equation*}
\frac{1}{C} \sum_{p=0}^n N_p[6p(p-1)+\frac{5n-3n^2}{2}] 
\end{equation*} 
into $\lvert E \rvert$ integers, $\lu(e):=\frac{\m(e)}{C}$, with $\lvert \lu(e)\rvert \leq \frac{2D}{C}$ where $\gcd\{\lu(e_1),\ldots,\lu(e_m)\}=1$.  Then we choose those for which
$$
\Big(  \ns\left(A(\Gamma_i)-\diag\left(\m(E_i)\right)\right)\Big)\cap \Z_{>0}^{\lvert E_i\rvert }\neq \emptyset
$$
for every $i=1,\dots,l$.

Note that some of the above integers may be zero depending on the existence of cycles in $\Gamma$.
This algorithm yields necessary conditions for all $S^1$-Hamiltonian actions with 
$$\max\left\{ \left \lvert \sum_{i=1}^nw_{iP}\right\rvert\right\}_{P\in \M^{S^1}}\leq D.$$
\end{alg} 
\begin{rmk}\label{C bound}
The upper bound for $C$ in the above algorithm can be improved
when $\Psi$ is injective. Indeed in this case we can use 
Theorem~\ref{hattori} (2), which is due to Hattori \cite{Ha}, to prove that $C$
must satisfy $1\leq C \leq N+1$, where $N+1$ is the number of fixed points.
\end{rmk}

\section{Minimal number of fixed points}\label{mnfp}

Suppose now that the number of fixed points is minimal. Then there are several important consequences which are explored next.
\subsection{An explicit basis for the equivariant cohomology}\label{hamiltonian minimal}
Suppose that $(\M,\omega)$ is a compact symplectic manifold of dimension $2n$ with a Hamiltonian $S^1$-action 
and minimal number of fixed points $P_0,\ldots,P_n$. Thus 
$$H^i(\M;\Z)=H^i(\C P^n;\Z)$$ 
for all $i$, and $N_p=1$ for every $p=0,\ldots,n$ (see Section~\ref{ec symplectic}).
Let us  order the fixed points $P_0,\ldots,P_n$ is such a way that $\lambda(P_i)=i$ for
all $i=0,\dots,n$. 

Tolman shows in \cite{T1} that, in this case, it is possible to recover the equivariant
cohomology ring $H_{S^1}^*(\M;\Z)$ (and hence $H^*(\M;\Z)$) from the isotropy representation of $S^1$ at the fixed
point set, and she gives an explicit basis for $H_{S^1}^*(\M;\Z)$ and $H^*(\M;\Z)$.
Here we review the main results of this construction, together with other useful facts.

Let $\cc^{S^1}\in H_{S^1}^*(\M;\Z)$ be the total equivariant Chern class of the tangent bundle and $\cc$ the total ordinary Chern class, that is, 
$$\cc^{S^1}=\sum_{j=0}^n \cc_j^{S^1}\quad \text{and} \quad \cc=\sum_{j=0}^n \cc_j.$$
We recall that, for every fixed point $P_i$, the restriction of the $j$-th equivariant Chern class to $P_i$ is given by
$$\cc_j^{S^1}(P_i)=\sigma_j(w_{i1},\ldots,w_{in})x^j,$$ 
where $\sigma_j$ denotes the $j$-th elementary symmetric polynomial,
$w_{i1},\ldots,w_{in}$ are the isotropy weights at $P_i$, and $x$ is the degree-$2$ generator of $H_{S^1}^*(\{P_i\};\Z)=\Z[x]$.
In particular, $\cc_1^{S^1}(P_i)=(\sum_{i=1}^nw_{in})x$ (see Section~\ref{ec}). 

The next result combines Proposition $3.4$ and Lemmas $3.8$ and $3.23$ of \cite{T1}.
\begin{prop}[Tolman]\label{psi e c}
Let the circle act on a compact symplectic manifold $(\M,\omega)$ of dimension $2n$ with moment map
$\psi\colon\M \to \R$, and let $P_0,\ldots,P_n$ be its fixed points, where $\lambda(P_i)=i$ for every $i$.
Then $\cc_1\neq 0$. Moreover,\\$\;$
\begin{equation}
\label{m.m increasing}
\psi(P_i)<\psi(P_j)  \quad\mbox{if and only if}\quad i<j\;,
\end{equation}
and
\begin{equation}
\label{c1 decreasing}
\displaystyle\frac{\cc_1^{S^1}(P_i)-\cc_1^{S^1}(P_j)}{x}>0  \quad\mbox{if and only if}\quad i<j\;,
\end{equation}
where $x$ is the generator of $H^*(\C P^{\infty};\Z)$.
\end{prop}
In particular, the equivariant symplectic form and the first Chern class restricted to the fixed point set are injective.

The following theorem combines Corollaries $3.14$ and $3.19$ of \cite{T1}.
\begin{thm}[Tolman]\label{bases}
Let the circle act on a compact symplectic manifold $(\M,\omega)$ of dimension $2n$ with moment map
$\psi\colon\M \to \R$, and let $P_0,\ldots,P_n$ be its fixed points, where $\lambda(P_i)=i$ for every $i$.
\begin{enumerate}
\item 
As a $H^*(\C P^{\infty};\Z)=\Z[x]$ module, $H_{S^1}^*(\M;\Z)$ is freely generated by $\tau_0,\tau_1,\ldots,\tau_n$, where
\begin{equation}\label{basis equivariant cohomology}
\tau_i=\frac{1}{C_i}\prod_{j=0}^{i-1}\left(\ce-\ce(P_j)\right)\quad\mbox{and}\quad C_i=\frac{\prod_{j=0}^{i-1}\left(\ce(P_i)-\ce(P_j)\right)}{\Lambda_i^-}\;.
\end{equation}
(In particular, $\tau_i\in H_{S^1}^{2i}(\M;\Z)$ for all $i$.)
\item
As a group, $H^*(\M;\Z)$ is freely generated by $\taut_0,\taut_1,\ldots,\taut_n$, where
\begin{equation}\label{basis cohomology}
\taut_i=r(\tau_i)=\frac{1}{C_i}\cc_1^i\;.
\end{equation} 
(In particular, $\taut_i\in H^{2i}(\M;\Z)$ for all $i$.)
\end{enumerate}
\end{thm} 
Here we consider  the empty product as being equal one so that $\tau_0=\taut_0=1$. 
Note that $\tau_i(P_i)=\Lambda_i^-$ and $\tau_i(P_j)=0$ for every $j\leq i-1$.

Since $\iota^*(\cc^{S^1})$ is determined by the weights at the fixed points,
using \eqref{basis equivariant cohomology} we can  explicitly compute $\iota^*(\tau_i)$. Moreover, by Lemma~\ref{coefficients},
we can compute the equivariant Chern classes in terms of the $\tau_i$s and the
ordinary Chern classes in terms of the $\taut_i$s.
\begin{rmk}\label{Ci's}
Since $\Lambda_i^-=(-x)^{i}\displaystyle\prod_{w_{ij}<0}\lvert w_{ij}\rvert$, it follows from \eqref{c1 decreasing} that
$C_i\in \Q_{>0}$ and then Theorem~\ref{bases} implies that it is a positive \emph{integer} for all $i$.
In the next sections the constant 
\begin{equation}\label{C1}
C_1=\frac{\ce(P_1)-\ce(P_0)}{\Lambda_1^-}\in \Z_{>0}
\end{equation}
will be particularly important. 
Also, observe that by \eqref{basis cohomology}, $C_i$ is the unique positive integer such that
$\displaystyle\frac{\cc_1^i}{C_i}$ is an integral generator of $H^{2i}(\M;\Z)$ for all $i=0,\ldots,n$.

Let $(\Gamma,\w)$ be an integral multigraph associated to the $S^1$-action on $(\M,\omega)$ (see Definition~\ref{img}).
For every edge $e\in \Gamma$, let
\begin{equation*}
\lu(e):=\frac{\tau_1(\ii(e))-\tau_1(\te(e))}{\w(e)x}=\frac{\m(e)}{C_1}\;.
\end{equation*}
As an easy consequence of Lemma~\ref{l dividing},  \eqref{sum m min} of Theorem~\ref{sum of m}
and the definitions above, we have the following result.
\begin{prop}\label{C1p}
Let $(\M,\omega)$ be a compact symplectic manifold of dimension $2n$ with a Hamiltonian $S^1$-action and $n+1$
fixed points.
Let $C_1$ be the positive integer such that $\displaystyle\frac{\cc_1}{C_1}$ is a generator
of $H^2(\M;\Z)$. Then
\begin{equation}\label{C1 divides}
C_1\;\;\;\mbox{divides}\;\;\;\frac{1}{2}n(n+1)^2.
\end{equation}
More precisely, 
let $(\Gamma,\w)$ be an integral multigraph associated to the $S^1$-action on $(\M,\omega)$. Then $\lu(e)$ is an integer
for every $e\in E$ and
\begin{equation*}
\sum_{e\in E}\lu(e)=\frac{1}{2}\frac{n(n+1)^2}{C_1}\;.
\end{equation*}
\end{prop}
Another important property regarding the constant $C_1$ will be proved in Section~\ref{hattori results}.
As we will see in Proposition~\ref{symp hattori} \eqref{bound C1} we also have
$$
1\leq C_1\leq n+1\;.
$$
(see also Remark~\ref{C bound}).
\end{rmk}
\subsubsection{Reversing the circle action.}\label{reversing flow}
Reversing the circle action, we obtain another basis for $H_{S^1}^*(\M;\Z)$ as a $\Z[x]$-module.
More precisely, the elements of this basis are $\tau^\prime_0,\ldots,\tau^\prime_n$, where $\tau^\prime_i\in H_{S^1}^{2i}(\M;\Z)$ is given by
 $$
 \tau^\prime_i=\frac{1}{C^\prime_i}\prod_{j=n-i+1}^{n}(\ce-\ce(P_j))\quad\mbox{with}\quad C^\prime_i=\frac{\prod_{j=n-i+1}^{n}(\ce(P_{n-i})-\ce(P_j))}{\Lambda_{n-i}^+}\;.
$$
Notice that $\tau^\prime_i(P_{n-i})=\Lambda_{n-i}^+$ and $\tau^\prime_i(P_j)=0$ for every $j\geq n-i+1$.

Moreover, as a group, $H^*(\M;\Z)$ is freely generated by $\taut^\prime_0,\ldots,\taut^\prime_n$, where $i$, $\taut^\prime_i\in H^{2i}(\M;\Z)$ is given by
$$
\taut^\prime_i=r(\tau^\prime_i)=\frac{\cc_1^i}{C^\prime_i}
$$
for all $i$.
It is easy to see that $C^\prime_i$ is a positive integer for all $i$ and then, since $\frac{\cc_1^i}{C^\prime_i}$ is an integral generator of $H_{S^1}^{2i}(\M;\Z)$, we have, by Remark~\ref{Ci's}, that
\begin{equation}\label{symmetries}
C_i=C^\prime_i
\end{equation}
for all $i$ and hence 
\begin{equation}\label{taus}
\taut_i=\taut^\prime_i
\end{equation}
for all $i$. Using \eqref{taus} and 
the ABBV Localization formula (cf. Corollary~\ref{abbv discrete}) we have
\begin{equation}\label{duality}
\int_\M \taut_i\taut_{n-i}=\int_\M \taut_i\taut^\prime_{n-i}=\int_\M \tau_i\tau^\prime_{n-i}=\sum_{j=0}^n\frac{\tau_i(P_j)\tau_{n-i}(P_j)}{\Lambda_j}=\frac{\Lambda_i^-\Lambda_{i}^+}{\Lambda_i}=1,
\end{equation}
for all $i=0,\ldots,n$.

\begin{rmk}
Observe that equations \eqref{symmetries} impose many polynomial relations among the
weights at the different fixed points. 
\end{rmk}

\subsection{Positive multigraphs} Let us now study conditions that would ensure the existence of a positive multigraph in the case where the number of fixed points is minimal.

Requiring  all the $\m(e)$s to be positive,  is equivalent to requiring
\begin{equation*}
\frac{\ce(\ii(e))-\ce(\te(e))}{x}> 0\quad\mbox{ for every $e\in E$.}
\end{equation*}
Since the number of fixed points is minimal, 
$(\M,\omega)$ is monotone. Indeed, since
$H^2(\M;\Z)=\Z$, we can choose $\omega$ such that
$\cc_1=C_1[\omega]$, with $[\omega]$ the generator of $H^2(\M;\Z)$ and $C_1$ a positive integer.
Moreover, we can choose the moment map so that  $\ce=C_1[\omega-\psi\otimes x]$, and then 
all the $\m(e)$s are positive if and only if 
$$
\psi(\ii(e))<\psi(\te(e)) \quad \text{for every $e\in E$,}
$$
or, equivalently, if and only if
\begin{equation}
\label{eq:index}
\lambda(\ii(e))<\lambda(\te(e)) \quad \text{for every $e\in E$}
\end{equation}
(cf. Proposition~\ref{psi e c}).

We conclude that, in this situation, we have the following result.
\begin{prop}\label{spm}
Let $(\M,\omega)$ be a compact symplectic manifold with a Hamiltonian $S^1$-action and a minimal number of fixed points.
Let $\Gamma=(V,E)$ be a multigraph associated to the $S^1$-action. 
Then $\Gamma$ is positive (resp. non-negative) if and only if its (directed) edges connect points to points with a greater  (resp. greater or equal index).
 \end{prop}
  If we can guarantee the existence of positive  multigraphs for every action within a certain class of circle actions, then we can use Algorithm 1. of Section~\ref{algo}.
It is then natural to ask whether these multigraphs  exist. There are many situations where we can guarantee the existence of such multigraphs associated to a given circle action. Here we present a few.

\vspace{0.5cm}

{\bf i)} Whenever $\M$ is four dimensional with $3$ isolated fixed points (see Section~\ref{dim4}).

\vspace{0.3cm}
{\bf ii)} Whenever $\M$ is six dimensional with $4$ fixed points: Ahara \cite{Ah} and
Tolman \cite{T1} prove the existence of a positive multigraph for every such Hamiltonian circle action.

\vspace{0.3cm}
{\bf iii)} When $\M$ is $8$-dimensional with $5$ fixed points and the $S^1$-action extends to a $T^2$-action: 
\begin{prop}\label{t2}
Suppose that $(\M,\omega)$ is $8$-dimensional and the $S^1$-action has $5$ fixed points.
If this action
 extends to an effective Hamiltonian $T^2$-action, then there exists a  positive multigraph associated to the circle action. 
 \end{prop}
 \begin{proof}
 To show this we just have to prove that, for each of these actions, there  exists a  multigraph satisfying \eqref{eq:index}. 

For that, let us consider the \emph{x-ray} of $(\M,\omega,\phi)$, where $\phi$ is the $T^2$-moment map,  given by the closed orbit type stratification $\mathcal{X}$ of $\M$, together with  the convex polytopes $\phi(X)$ for each $X\in \mathcal{X}$  \cite{T2}. Let $P_i$ be the $S^1$-fixed point with $\lambda(P_i)=i$. Each line in the x-ray is the image of a connected component of the set of points with a given $1$-dimensional stabilizer. If there is no line in the x-ray containing the images of  $P_0$ and $P_1$  then there are at most $3$ lines in the x-ray through $\phi(P_1)$. 

If there were $3$ of these lines, then there would  exist a $4$-dimensional manifold  $X\in \mathcal{X}$, fixed by a circle inside $T^2$, containing $P_1$ and one additional fixed point. If there were $2$ lines then we would either have two $4$-dimensional manifolds $X_1,X_2\in \mathcal{X}$ (each fixed by a different  circle inside $T^2$) containing $P_1$, or a $6$-dimensional manifold  $X\in \mathcal{X}$ fixed by a circle inside $T^2$, containing $P_1$. In the first case, one of the manifolds $X_1$, $X_2$ could only have two fixed points and, in the latter, the  manifold $X$ would have  at most three fixed points. Finally, if there existed just one line through $\phi(P_1)$, it would be  the image of an $8$-dimensional manifold with at most $4$ fixed points.

Since the minimal number of fixed points on a $2m$-dimensional  $S^1$-Hamiltonian manifold  is $m+1$ all the above cases are impossible. Therefore, we conclude that there must exist a line in the x-ray containing the images of $P_0$ and $P_1$. Hence, there exists a  manifold $X\in \mathcal{X}$, a component of the set of points with a certain $1$-dimensional stabilizer, which contains $P_0$ and $P_1$. Note that, since the $T^2$-action is effective, we have  $\dim X \leq 6$. Moreover,  $P_0$ and $P_1$ are  fixed points of the restriction of the original $S^1$-action on $\M$ to $X$, respectively of index $0$ and $2$. By the classification of Hamiltonian $S^1$-actions on $4$-dimensional manifolds with isolated fixed points \cite{K} and the classification of   Hamiltonian $S^1$-actions on $6$-dimensional manifolds with  a minimal number of fixed points \cite{T1}, we conclude that there exists a multigraph for this restricted action on $X$ with an edge connecting $P_0$ and $P_1$.  

Similarly, we can conclude that there  exists a  manifold $X^\prime \in \mathcal{X}$ which is a component of the set of points with a certain $1$-dimensional stabilizer, containing $P_4$ and $P_3$, and that there is a multigraph for the restricted $S^1$-action on $X^\prime$ with an edge connecting $P_3$ and $P_4$.

Since any multigraph for the original $S^1$-action on $\M$ restricts to multigraphs for the restrictions of the action to $X$ and $X^\prime$, we conclude that there exists a multigraph for the $S^1$-action on $\M$ which has an edge connecting $P_0$ to $P_1$ and another one connecting $P_3$ and $P_4$, thus satisfying 
\begin{equation}\label{eq:index2}
\lambda(\ii(e))\leq \lambda(\te(e)) \quad \text{for every $e\in E$.}
\end{equation}

This argument can be easily adapted to prove the existence of a multigraph satisfying \eqref{eq:index2} with a strict inequality. For that, let us  consider all multigraphs with an edge connecting $P_0$ to $P_1$ and an edge connecting $P_3$ to $P_4$. Then we just have to show that among these there is a multigraph with no cycles at $P_2$. 

Let us consider again the x-ray of  $(\M,\omega,\phi)$.  If there are $4$ lines through  $\phi(P_2)$ then there exist four $2$-dimensional $S^1$-manifolds $X_i$ with different  $1$-dimensional   stabilizers, having $P_2$ in their fixed-point set. Since the only compact $2$-dimensional manifolds admitting a Hamiltonian circle action are spheres, where the circle acts by rotation, they all must have an additional fixed point different from $P_2$.  Hence, the multigraphs for the restricted $S^1$-actions on $X_i$ consist of only one edge with $P_2$ at one of the endpoints and another fixed point at the other. We conclude that there exists a multigraph for the original $S^1$-action on $\M$ that has $4$ edges with $P_2$ as an endpoint and  a different fixed point at the other end (thus with no cycles at $P_2$). 
 
 If there are $3$ lines through $P_2$ on the x-ray then we  either have two $4$-dimensional manifolds $X_1,X_2\in \mathcal{X}$ (each fixed by a different  circle inside $T^2$)  containing $P_2$, or a $6$-dimensional manifold  $X\in \mathcal{X}$ and a $2$-dimensional manifold $X^\prime$ fixed by two different circles inside $T^2$, with both $X$ and $X^\prime$ containing $P_2$. In the first case, both manifolds $X_1$, $X_2$  have exactly $3$ fixed points and, in the latter, the  manifold $X$ has at most $4$ fixed points while $X^\prime$ has $2$. 
 
 By the classification of Hamiltonian $S^1$-actions on $4$-dimensional manifolds with isolated fixed points \cite{K} and the classification of   Hamiltonian $S^1$-actions on $6$-dimensional manifolds with  a minimal number of fixed points \cite{T1}, we conclude that, in both situations, there exist multigraphs for the restricted $S^1$-actions on these submanifolds that have no cycles at $P_2$. We conclude that in all cases there exists a multigraph for the original $S^1$-action on $\M$ that has no cycles at $P_2$ and so this multigraph is necessarily positive. 
\end{proof}
\vspace{0.3cm}
{\bf iv)} Whenever $\M$ is $8$-dimensional with $5$ fixed points and none of the isotropy weights of the action is $1$:
\begin{prop}\label{not pm1}
Suppose that $(\M,\omega)$ is $8$-dimensional and the $S^1$-action has $5$ fixed points.
If none of the isotropy weights of the action is $1$ then there exists a positive multigraph associated to the $S^1$-action.
\end{prop}
\begin{proof}
Observe that, by Lemma~\ref{pairs}, none of the weights is $-1$.
A similar argument to the one in the proof of Proposition~\ref{t2} can be carried out, now using the closed orbit type stratification $\mathcal{X}$ of the $S^1$-manifold $\M$, where  we consider the elements $X$ of $ \mathcal{X}$ that are  connected components of the set of points with a given finite cyclic group $\Z_k$ as stabilizer ($k\neq 1$). 

Let us assume that the negative weight at $P_1$ is $-k$ with $k>1$. We will show that there is at least one weight at $P_0$ that is equal to $k$. If this is not the case then the  number of weights at $P_0$ which are multiples of $k$ (but different from $k$) has to be at least two but at most $3$ (since the action is effective). Let $X$ be the (connected) component of the points fixed by $\Z_k$ containing $P_0$ and $P_1$. If there were two weights at $P_0$ which were multiples of $k$ then $X$   would be a $4$-dimensional manifold admitting an effective $S^1\cong S^1/\Z_k$-action. This manifold would then have more than two non-trivial chains of gradient spheres\footnote{A chain of gradient spheres is a sequence of gradient spheres $S_1,\ldots,S_l$ such that the south pole of $S_0$ is a minimum for the moment map, the north pole of $S_{i-1}$ is the south pole of $S_i$ for each $1<i\leq l$, and the north pole of $S_l$ is a maximum for the moment map. A chain is non-trivial if it contains more than one sphere, or if it contains one sphere whose stabilizer is non-trivial.}, which is impossible by Karshon's classification results on Hamiltonian $S^1$-actions on $4$-manifolds \cite[Proposition $5.13$]{K}. (Indeed, the negative weight at $P_1$ for this effective action is $-1$ and none of the weights at $P_0$ is $1$.) If, on the other hand,  there were $3$ weights at $P_0$ which were multiples of $k$, then $X$ would be $6$-dimensional. Then the effective Hamiltonian $S^1\cong S^1/\Z_k$-action on $X$ would have $4$ fixed points (since when $\dim \M=6$ and the number of fixed points is not minimal one must have, by \eqref{eq:reverse}, at least $6$ fixed points). By Tolman's classification of  Hamiltonian actions on $6$-manifolds with a minimal number of fixed points \cite{T1}, none of these manifolds has a negative weight $-1$ at the point of index-$2$ and no weight equal to $1$ at the minimum.  We conclude that, if the negative weight at $P_1$ is $-k$ with $k>1$, then there is at least one weight at $P_0$  that is equal to $k$. 

Similarly, we can conclude that, if the positive weight at $P_3$ is $k>1$, then there is at least one weight at $P_4$ that is equal to $-k$ and so, given an $S^1$-action satisfying the assumptions above, we can always choose a multigraph that is non-negative. To show that we can always choose one that is positive we still have to show that we can choose a multigraph with no cycles at $P_2$.

For that, let us assume that we have chosen a multigraph with an edge $e_{01}$  connecting $P_0$ to $P_1$ and an edge connecting $P_3$ to $P_4$ and that, at $P_2$, the action has  a weight $k$ and a weight $-k$ with $k>1$. Let  $X$ be the (connected) component of the points fixed by $\Z_k$ containing $P_2$ and consider the effective Hamiltonian  $S^1\cong S^1/\Z_k$-action on $X$. This action has $P_2$ as an index-$2$ fixed point and so the index-$0$ point has to be $P_0$ or $P_1$. Let $W_0=\{w_{01},\ldots,w_{04}\}$ be the multiset of weights at $P_0$ with $w_{01}$ the weight associated to the edge $e_{01}$, and let   $W_1=\{w_{11},w_{12},\ldots,w_{14}\}$ be the multiset of weights at $P_1$ with  $w_{11}=-w_{01}$. If there is no weight in $W_0\setminus \{w_{01}\}\cup W_1\setminus \{w_{10}\}$ that is equal to $k$, then, again by Karshon's classification of Hamiltonian $S^1$-actions on $4$-manifolds \cite{K} and Tolman's classification of Hamiltonian $S^1$-actions on $6$-manifolds with a minimal number of fixed points \cite{T1}, we get a contradiction. The result then follows.
\end{proof}
{\bf v)} Whenever the $S^1$-action extends to a ``GKM action" on $\M$ with a minimal number of fixed points: Then there exists a natural multigraph associated to $\M$ that, when
$\M$ has a minimal number of fixed points, is the \emph{complete graph} on the set of vertices given by the fixed points (hence positive).

Indeed, let $(\M,\omega,\phi)$ be a compact symplectic manifold with a Hamiltonian $T$-action and isolated fixed points, where $\dim(T)>1$.
For every $K\subseteq T$ let $\M^K$ be the points fixed by $K$, $\mathfrak{t}$ the Lie algebra of $T$, and $\mathfrak{t}^*$ its dual.
Then $\M$ is a \textbf{GKM (Goresky-Kottwitz-MacPherson) manifold} \cite{GKM} if the weights $\alpha_{1},\ldots,\alpha_{n}\in \mathfrak{t}^*$ of the isotropy representation of $T$ on $T_P\M$
are pairwise linearly independent, for every fixed point $P$. This is equivalent to saying that, for every codimension one subgroup $K\subset T$,
the connected components of $\M^K$ are either points or 2-spheres. 
This class of spaces include two families of well-known spaces, coadjoint orbits of a simple Lie group $G$ endowed with the
action of a maximal torus $T\subset G$, and toric symplectic manifolds. 

Fix $P\in \M^T$ and let 
$\alpha_{1},\ldots,\alpha_{n}\in \mathfrak{t}^*$ be the weights of the isotropy representation of $T$ on $T_P\M$.
Let
$S_{i}$ be the 2-sphere fixed by $\exp(\ker \alpha_i)\subset T$ for every $i=1,\ldots,n$. Then the $T$-action on $S_{i}$
has two fixed points $P$ and $Q$, and the weight of the $T$ representation on $T_Q(S_{i})$ is $-\alpha_i$. 
Assume that the $S^1$-action we are starting from is the action of a generic circle inside $T$, and let $\psi\colon \M\to \R$
be its moment map.
Then, we can associate to each of these two spheres an edge going from $P$ to $Q$, with $\psi(P)<\psi(Q)$, which is equivalent to saying that the weight of the $S^1$-action
at $P$ is positive.
Since $\alpha_1,\ldots,\alpha_n$ are pairwise linearly independent, it follows that $S_i\cap S_j=\{P\}$ for every $1\leq i<j\leq n$, and, if $\M$
has exactly $n+1$ fixed points, the graph
constructed above must be a complete, hence  positive, multigraph.
(For a classification of GKM manifolds with a minimal number of fixed points see \cite{Mo}.)

\begin{rmk}\label{l=1}
Notice that for every positive multigraph with a minimal number of fixed points, there must exist an edge $e_{01}$ from $P_0$ to $P_1$, and an edge $e_{n\,n+1}$
from $P_{n}$ to $P_{n+1}$. Moreover, from the definitions in Sections~\ref{dmg} and \ref{hamiltonian minimal},  it follows that $\m(e_{01})=\m(e_{n\,n+1})=C_1$, and so $\lu(e_{01})=\lu(e_{n\,n+1})=1$.
\end{rmk}
\begin{rmk}\label{upper bound c1^n}
Let $(\M,\omega)$ be a $2n$-dimensional compact symplectic manifold with a Hamiltonian $S^1$-action and fixed points $P_0,\ldots,P_n$.
Suppose that there exists an integral multigraph $\Gamma$ associated to the action which is the complete graph on $n+1$ vertices
(hence positive).
Then it is easy to see that, by the ABBV Localization formula, for every $i=0,\ldots,n$, we have
\begin{equation}\label{c1^n}
\int_\M\cc_1^n=\int_\M \prod_{j\neq i}\left(\cc_1^{S^1}-\cc_1^{S^1}(P_j)\right)=\frac{\prod_{j\neq i}\left(\cc_1^{S^1}(P_i)-\cc_1^{S^1}(P_j)\right)}{\Lambda_i}=\prod_{h=1}^n\m(e_h),
\end{equation}
where $e_1,\ldots,e_n$ is the set of edges ending or starting at $P_i$, and each $\m(e_h)$ is a positive integer.
By  \eqref{sum m min} in Theorem~\ref{sum of m}, we have  
\begin{align*}
\sum_{h=1}^n\m(e_h) & \leq \sum_{e\in E}\m(e) - \left(\frac{n(n+1)}{2}-n \right)\\ 
& \leq \frac{n(n+1)^2}{2} - \frac{n(n-1)}{2} = \frac{n(n^2+n+2)}{2},
\end{align*} 
since the graph has $n(n+1)/2$ edges and $m(e)\geq 1$.
Then we have
\begin{equation*}
\int_\M \cc_1^n=\prod_{h=1}^n\m(e_h)\leq \left( \frac{n^2+n+2}{2}\right)^n.
\end{equation*}
The same argument can be carried out when the action has an integral multigraph $\Gamma$ for which there exists a vertex $P_i$ 
connected to the other $n$ vertices of $\Gamma$ through $n$ (undirected) edges of $\Gamma$. 
It would be interesting to generalize this estimate to the case in which the multigraph does not necessarily satisfy the above property (see also Remark~\ref{225}).
\end{rmk}

\subsection{A refinement of the algorithm}\label{refinement}
In all the above situations or, in general, whenever we are able to guarantee that, within a given class of  $S^1$-actions with minimal number of fixed points, there always exists a  non-negative multigraph for each action, we can use Algorithm~\ref{alg:1} with slight improvements.

\setcounter{alg}{0}
\begin{alg}[B] \label{alg:5} For each possible multigraph $\Gamma$ with connected components $\Gamma_1,\ldots,\Gamma_l$, we look for partitions of
$$
\frac{1}{2}n(n+1)^2
$$
into $\lvert E^\prime\rvert=\lvert E\rvert - \lvert E^{\circlearrowleft} \rvert$ positive integers  $\m(e)$. Then we add zeros $\m(e)=0$ whenever the edge $e$ is a cycle and choose, among  the resulting  sequences of $\lvert E \rvert$ nonegative numbers, those  for which
$$
\Big(  \ns\left(A(\Gamma_i)-\diag\left(\m(E_i)\right)\right)\Big)\cap \Z_{>0}^{\lvert E_i\rvert }\neq \emptyset
$$
for every $i=1,\dots,l$.
\end{alg}

If, in addition, we know that there exists an edge $e_{01}$  from $P_0$ to $P_1$ and/or an edge $e_{n\,n+1}$ between $P_n$ and $P_{n+1}$ (for instance when the multigraph is positive), then we can use  Remark~\ref{l=1} to further improve the algorithm, adapting Algorithm~\ref{alg:2} to this situation. Here note that, when we have a minimal number of fixed points, the number of cycles $\lvert  E^{\circlearrowleft}  \rvert$ of a non-negative multigraph satisfies
\begin{align*}
\lvert  E^{\circlearrowleft}  \rvert & \leq \sum_{p=0}^n \min\{p,n-p\} = \sum_{p=1}^{\lfloor n/2 \rfloor} p + \sum_{p=\lfloor n/2 \rfloor + 1}^{n-1} (n-p) \\ & = \lfloor n/2 \rfloor^2-(n-1)\lfloor n/2 \rfloor + \frac{n^2-n}{2} \\
 & = \left\{ \begin{array}{ll} \frac{n^2}{4}, & \text{if $n$ is even} \\ \\ \frac{n^2-1}{4}, & \text{if $n$ is odd} \end{array} \right. \;,
\end{align*}
and so the constant $C$ in \eqref{ub C} satisfies
$$
1\leq C \leq \frac{n(n+1)^2}{n(n+1)-2 \lvert  E^{\circlearrowleft}  \rvert } \leq  \left\{ \begin{array}{ll} 2n+1, & \text{if $n$ is even} \\ \\ 2n, & \text{if $n$ is odd} \end{array} \right.\;.
$$
However, by Proposition~\ref{psi e c}, we know that the map $\Psi$ defined in \eqref{Psi} is injective and so we can use Theorem~\ref{hattori} to obtain a better upper bound for $C$,
$$
1\leq C \leq n+1,
$$
which is also the same as the one obtained when $\lvert  E^{\circlearrowleft}  \rvert =0$.
\begin{alg}[B] \label{alg:6} Given a graph  $\Gamma$ with connected components $\Gamma_1,\ldots,\Gamma_l$ with an edge $e_{01}$  from $P_0$ to $P_1$ (and/or an edge $e_{n\,n+1}$ between $P_n$ and $P_{n+1}$), 
for each divisor $C$ of 
$$
\frac{1}{2}n(n+1)^2
$$
satisfying
$$
1\leq C\leq n+1,
$$
we look for partitions of
$$
\frac{1}{2C}n(n+1)^2
$$
into $\lvert E^\prime \rvert =\lvert E \rvert - \lvert E^{\circlearrowleft} \rvert$ positive integers, $\lu(e):=\frac{\m(e)}{C}$ with $\lu(e_{01})=1$ (and/or  $\lu(e_{n\, n+1})=1$). Then we consider the corresponding integers $\m(e)=C\cdot \lu(e)$, add  zeros $\m(e)=0$ whenever the edge $e$ is a cycle,
and choose, among the resulting sequences of $\lvert E \rvert$ numbers, the ones  for which
$$
\Big(  \ns\left(A(\Gamma_i)-\diag\left(\m(E_i)\right)\right)\Big)\cap \Z_{>0}^{\lvert E_i\rvert}\neq \emptyset,
$$
for every $i=1,\dots,l$.
\end{alg}

Suppose now that we cannot guarantee the existence of non-negative multigraphs for every circle action considered and so we have to deal with the possibility that 
the $\m(e)$s might be negative.
Then we can use Algorithms~\ref{alg:3} and \ref{alg:4} in Section~\ref{algo}, noting  that, in the case of Algorithm~\ref{alg:4}, the map $\Psi$ defined in \eqref{Psi} is now injective by Proposition~\ref{psi e c},
and so, by Remark~\ref{C bound}, we again have $1\leq C \leq n+1$.
\subsection{Known examples of $S^1$-Hamiltonian manifolds with minimal number of fixed points}\label{known examples}
Let us now describe the known examples of Hamiltonian circle actions with a minimal number of fixed points.
However, since the first three examples arise as coadjoint orbits of simple Lie groups, we first review a few general facts.

Let $G$ be a compact simple Lie group with $\mathfrak{g}=Lie(G)$, and let $T\subset G$ be a maximal torus
with $\mathfrak{t}=Lie(T)$ and $\mathfrak{t}^*=Lie(T)^*$. 
Let $\Delta\subset \mathfrak{t}^*$ be the set of roots of $G$, and $\Delta_0\subset \Delta$ a choice of simple roots.
We use the Killing form to regard 
$\mathfrak{t}^*$ as a subspace of $\mathfrak{g}^*$.
The coadjoint orbit through a point $p_0\in \mathfrak{t}^*$, $O_{p_0}$, is a compact
manifold with a natural symplectic structure given by the Kostant-Kirillov symplectic form $\omega$.
It also inherits a natural Hamiltonian $G$-action, whose moment map is given by the inclusion
$O_{p_0}\hookrightarrow \mathfrak{g}^*$. Thus, restricting the action to a generic circle $S^1\subset T$, the compact symplectic manifold
$(O_{p_0},\omega)$ has a Hamiltonian $S^1$-action
with moment map $\psi\colon O_{p_0}\hookrightarrow\mathfrak{g}^*\to Lie(S^1)^*\simeq \R^*$, where the second
map is the projection map.
\begin{exm}\label{cpn}{\bf (The complex projective space)}\\
Let $G=SU(n+1)$, and $T^n\subset G$ the torus of diagonal matrices. Let $\{x_0,\ldots,x_n\}$ be the standard basis
of $(\R^{n+1})^*$. Then 
a standard choice of simple roots is 
$$\Delta_0=\{x_i-x_{i+1}=\alpha_i\mid 0\leq i\leq n-1\}.$$
Let $p_0$ be a generic point in $\cap_{i=1}^{n-1}\mathcal{H}_{\alpha_i}$, where $\mathcal{H}_{\alpha_i}\subset \mathfrak{t}^*$ is the hyperplane 
orthogonal to $\alpha_i$. Then $O_{p_0}$ is isomorphic to $\C P^n$.
Let $S^1\subset T$ be a generic circle generated by $\xi=(\xi_0,\ldots,\xi_n)\in \mathfrak{t}$
such that $\alpha_i(\xi)=\xi_i-\xi_{i+1}\in \Z_{>0}$
 for every $i=0,\ldots,n-1$.
 A standard computation shows that the $S^1$-fixed points are 
$$\{P_i=[0,\ldots,\overbrace{1}^{i-th},\ldots,0]\mid i=0,\ldots,n\},$$ 
 the set of weights at $P_i$ is given by
$$\{\xi_i-\xi_j\mid 0\leq j\neq i\leq n\},$$
 and $\lambda(P_i)=i$ for every $i=0,\ldots,n$.
 
A natural choice of integral multigraph is the complete graph on $n+1$ vertices, which is
also the GKM graph associated to $\C P^n$ equipped with the $T^n$-action (see Section~\ref{refinement}, {\bf v)}).
 
It is well known that, in this case, the unique positive integer $C_1$ for which $\cc_1/C_1$ is a generator of
 $H^2(\C P^n;\Z)$ is $n+1$, which agrees with \eqref{C1}.
\end{exm}
\begin{exm}\label{gr}{\bf (The Grassmannian of oriented two planes in $\R^{2n+1}$)}\\
Let $G=SO(2n+1)$ and $T^n$ a maximal torus and let us identify $\mathfrak{t}^*$ with $(\R^n)^*$ with
standard basis $\{x_0,\ldots,x_{n-1}\}$. 
Then a choice of simple roots is given by 
$$\Delta_0=\{x_i-x_{i+1}=\alpha_i,\;i=0,\ldots,n-2\}\cup\{x_{n-1}=\alpha_{n-1}\}.$$
Let $p_0$ be a generic point in  $\cap_{i=1}^{n-1}\mathcal{H}_{\alpha_i}$, where $\mathcal{H}_{\alpha_i}\subset \mathfrak{t}^*$ is the hyperplane 
orthogonal to $\alpha_i$. Then $O_{p_0}$ is isomorphic to $Gr_2^+(\R^{2n+1})$, the Grassmannian of oriented two planes in
$\R^{2n+1}$, which is a symplectic manifold of dimension $2(2n-1)$. Let $S^1\subset T$ be a generic circle generated by $\xi=(\xi_0,\ldots,\xi_{n-1})\in \mathfrak{t}$
such that $\alpha_i(\xi)\in \Z_{>0}$
 for every $i=0,\ldots,n-1$.
Then a standard computation shows that  this action has $2n$ fixed points, which can be identified with the elements $y_i$s of $\mathfrak{t}^*$ given by
$$
y_0=-x_0,\;\ldots,y_{n-1}=-x_{n-1},\;y_n=x_{n-1},\;\ldots,\;y_{2n-1}=x_0\,.
$$
Moreover,
the weights at $y_i$ are given by 
$$\{(y_j-y_i)(\xi)\mid j\neq i,\;j\neq 2n-1-i\}\cup\{-y_i(\xi)\},$$ 
and $\lambda(y_i)=i$
for every $i$.
Observe that, since $Gr_2^+(\R^{2n+1})$ is a Hamiltonian $S^1$-manifold with a minimal number of fixed points,
by the argument in Section~\ref{ec symplectic}, we have 
$$H^j(Gr_2^+(\R^{2n+1});\Z)=H^j(\C P^{2n-1};\Z)$$ 
for every $j=0,\ldots,2(2n-1)$. However the cohomology ring is different.

As in the previous example, a natural choice of integral multigraph is the complete graph on $2n$ vertices, which is
 also the GKM graph associated to $Gr_2^+(\R^{2n+1})$ equipped with the $T^n$-action (see Section~\ref{refinement}, {\bf v)}).
 
It is well known that, in this case, the unique positive integer $C_1$ such that $\cc_1/C_1$ is a generator of
 $H^2(Gr_2^+(\R^{2n+1});\Z)$ is $n$, which agrees with the definition of $C_1$ in \eqref{C1}.
 \end{exm}
  \begin{exm}\label{g2}
 Let $G_2$ be the exceptional simple Lie group and $T^2$ a maximal torus. Let $p_0$ be a generic point in $\mathcal{H}_{\alpha_1}$, where $\alpha_1$
 denotes the short simple root, and $\mathcal{H}_{\alpha_1}$ the hyperplane orthogonal to it.
 Then $O_{p_0}$ is a $10$-dimensional coadjoint orbit with an $S^1\subset T^2$-Hamiltonian action and $6$ fixed points.
 The $T^2$-action is also GKM, and hence a natural choice of multigraph is the complete multigraph on six vertices (for more details see \cite{Mo}).
  \end{exm}
\begin{exm}\label{fano}{\bf(The Fano manifolds $V_5$ and $V_{22}$)}\\
Let $(\M^{6},\omega)$ be a $6$-dimensional compact symplectic manifold with a Hamiltonian $S^1$-action and $4$ fixed points. Let $C_1$ be the unique positive integer $C_1$ such that $\cc_1/C_1$ is a generator of
$H^2(\M;\Z)$ (see Remark~\ref{Ci's}). As we will see in \eqref{bound C1} in Proposition~\ref{symp hattori},
$C_1$ can be either $1,2,3$ or $4$. As examples of manifolds with $C_1$ equal to $3$ and $4$
we have respectively $Gr_2^+(\R^5)$ and $\C P^3$ (see Examples~\ref{gr} and \ref{cpn}).
For the remaining two cases we have two Fano $3$-folds, which are known as $V_5$ and $V_{22}$. Indeed, McDuff \cite{M}
proved that they can be explicitly given a symplectic structure, 
and they possess a Hamiltonian  $S^1$-action. These Fano manifolds have the following properties:
\begin{itemize}
\item $V_5$: Here $C_1=2$, the cohomology ring is $$H^*(V_5;\Z)=\Z[x_1,x_2]/(x_1^2-5x_2,x_2^2),$$ and the isotropy weights
at the fixed points are
\begin{equation}\label{v5}
\{\{1,2,3\},\;\{-1,1,4\},\;\{-1,-4,1\},\;\{-1,-2,-3\}\}\,.
\end{equation}
\item $V_{22}$: Here $C_1=1$, the cohomology ring is $$H^*(V_{22};\Z)=\Z[x_1,x_2]/(x_1^2-22x_2,x_2^2),$$
and the isotropy weights
at the fixed points are
\begin{equation}\label{v22}
\{\{1,2,3\},\;\{-1,1,5\},\;\{-1,-5,1\},\;\{-1,-2,-3\}\}\,.
\end{equation}
\end{itemize}
There is exactly one  positive multigraph associated to these two actions, which is shown in Case III of Figure~\ref{graphs6}.
\end{exm}

\section{Implications of Hattori's results in the symplectic category}\label{hattori results}
In this section we recall some of the results obtained by Hattori in \cite{Ha},
and show how they can be used to obtain relations among the Chern numbers of a compact
symplectic manifold $(\M,\omega)$ endowed with a Hamiltonian circle action and minimal number
of fixed points. In particular, in Theorem~\ref{m not 2}, we derive their consequences in the case in which $\M$ is $8$-dimensional. 

Let $(\M,\J)$ be a compact almost complex manifold of dimension $2n$, equipped with an $S^1$-action which preserves $\J$, and has isolated fixed points
$P_0,\ldots,P_N$. 
Then the set of weights of the $S^1$-representation on $T_{P_i}\M$ is well-defined
for all $i$. Let $\{w_{i1},\ldots,w_{in}\}$ be the multi-set of weights at $P_i$.

Let $\LL$ be an admissible complex line bundle over $\M$ (see Section~\ref{eclb}).
Then for any lift $\LL^{S^1}$ there exist integers $a_0,\ldots,a_N$ such that
\begin{equation}\label{eq:ai}
\LL^{S^1}(P_i)=t^{a_i}\quad\mbox{for every}\quad i=0,\ldots,N\;,
\end{equation}
which are determined up to a constant. 
Following the terminology in \cite{Ha}, we say that 
an admissible line bundle $\LL$ is called \emph{fine} if the 1-dimensional representations 
$\LL^{S^1}(P_0),\ldots,\LL^{S^1}(P_N)$ are pairwise distinct, i.e. if $a_i\neq a_j$ for every
$i\neq j$. Moreover, 
 it is called \emph{quasi-ample} if it is fine and $\int_{\M}\cc_1(\LL)^n\neq 0$.

Let $\LL$ be a fine complex line bundle and $\LS$ an equivariant extension.
For every $i=0,\ldots,N$
 define
\begin{equation}\label{varphi}
\varphi_i(t):= \frac{\displaystyle\prod_{j\neq i}(1-t^{a_i-a_j})}{\displaystyle\prod_{k=1}^{n}(1-t^{w_{ik}})}\;.
\end{equation}
Notice that 
\begin{equation}\label{varphi index}
\varphi_i(t^{-1})=\ind_{S^1}\left(\prod_{j\neq i} \left(1-(\LL^{S^1})^{-1}t^{a_j}\right)\right)\in \Z[t,t^{-1}]\,,
\end{equation}
and so $\varphi_i(t)\in \Z[t,t^{-1}]$.

In the following we recall Theorems $4.2$ and $5.7$ in \cite{Ha}, specializing them to the
case of almost complex manifolds as well as other useful facts proved in \cite{Ha}.
\begin{theorem}[Hattori]\label{hattori}
Let $(\M,\J)$ be a compact almost complex manifold of dimension $2n$, equipped with an $S^1$-action which preserves $\J$ and has isolated fixed points
$P_0,\ldots,P_N$. Let $\LL$ be a fine complex line bundle, and for every $i=0,\ldots,N$, let $\varphi_i(t)\in \Z[t,t^{-1}]$ be as in
\eqref{varphi}. Then there exists a unique sequence $r_0(t),\ldots,r_N(t)$ of elements of $\Z[t,t^{-1}]$ such that
$$
\varphi_i(t)=r_0(t)+r_1(t)t^{a_i}+\cdots + r_N(t)t^{Na_i}\quad\mbox{for all}\quad i\;.
$$
Moreover, the $r_s(t)$s satisfy the following properties:
\begin{enumerate}
\item\label{r_0} 
$$
r_0(t)=\td(\M)=\sum_{i=0}^N\binom{\lambda_i}{n}\;,
$$
where $\lambda_i$ is the number of negative weights at $P_i$.
\item\label{k_0} If there exists $k_0\in \Z_{\geqslant 0}$ and $d\in \Z$ such that
\begin{equation}\label{D1}
\sum_{k=1}^nw_{ik}=k_0a_i+d\quad\mbox{for all}\quad i=0,\ldots,N\;,
\end{equation}
then $k_0\leqslant N+1$.
\item In \eqref{k_0}, if $k_0>0$, then, setting $l_0=N+1-k_0$, we have
$$
r_s(t)=0\quad\mbox{for all}\quad s>l_0\;,
$$
and
$$
r_{l_0-s}(t)=(-1)^{N-n}r_s(t^{-1})t^{-(d+\sum a_j)}\quad\mbox{for}\quad s\leqslant l_0\;.
$$
\item In \eqref{k_0}, if $k_0=0$, then $r_0=0$ and
$$
r_{N+1-s}(t)=(-1)^{N-n}r_s(t^{-1})t^{-(d+\sum a_j)}\quad\mbox{for}\quad 1\leqslant s\leqslant N\;.
$$
\end{enumerate}
\end{theorem}
More explicitly, the functions $r_s(t)$s are given by
\begin{equation}\label{def expl r}
r_s(t)=(-1)^s\sum_{i=0}^N\frac{\displaystyle\sum_{j_1<\cdots <j_s,j_{\nu}\neq i}t^{-(a_{j_1}+\cdots +a_{j_s})}}{\prod_{k=1}^n(1-t^{w_{ik}})}.
\end{equation}
For every $s=0,\ldots,N$,
let $k_s(\LS,t)$ be the $S^1$-equivariant bundle associated to $\LL^{S^1}$ defined inductively on $s$ by
\begin{equation}\label{k0}
k_0(\LL^{S^1},t)=1
\end{equation}
and
\begin{align}
k_s(\LL^{S^1},t)=\displaystyle\sum_{j_1<\cdots< j_s} & (t^{a_{j_1}}-\LL^{S^1})\cdots (t^{a_{j_s}}-\LL^{S^1}) \nonumber \\
\label{k ind}&-\displaystyle\sum_{\nu=1}^s\binom{N-s+\nu}{\nu}(-\LL^{S^1})^{\nu}k_{s-\nu}(\LL^{S^1},t)\;.
\end{align} 
Then, for every $s=0,\ldots,N$, the functions $r_s(t)$s satisfy 
\begin{equation}\label{explicit r}
r_s(t^{-1})=(-1)^s\ind_{S^1}\left(k_s(\LS,t)\right)\;.
\end{equation}

\begin{theorem}[Hattori]\label{hattori2}
Let $(\M,\J)$ be a compact almost complex manifold of dimension $2n$, equipped with an $S^1$-action which preserves $\J$ and has isolated fixed points. Assume the Euler characteristic of $\M$ is $n+1$. If $\LL$ is a quasi-ample complex line bundle satisfying \eqref{D1} with $k_0=n+1$ then the multisets of  weights of the $S^1$-action at each fixed point $P_i$ are given by
$$
\{w_{ik}\} = \{ a_i-a_j\}_{j\neq i},
$$
where $a_0,\ldots,a_n$ are defined up to a constant as in \eqref{eq:ai}. In particular, the multisets of weights coincide with those of $\C P^n$ with the standard circle action described in Example~\ref{cpn} (with  $\LL$  the hyperplane bundle).
\end{theorem}

In the following, we derive the consequences of Theorem~\ref{hattori} when
$(\M,\omega)$ is a compact symplectic manifold of dimension $2n$
with a Hamiltonian $S^1$-action
and a minimal number of fixed points $P_0,\ldots,P_n$.
As usual, we endow $(\M,\omega)$ with an almost complex structure $\J$
compatible with $\omega$, which is invariant under the $S^1$-action.

\begin{defin}\label{ipp}
Let $(\M,\omega)$ be a compact symplectic manifold with a Hamiltonian $S^1$-action with a
minimal number of fixed points. Let $\psi\colon \M\to \R$ be the moment map.
We say that $\omega-\psi\otimes x$ is 
\begin{itemize}
\item[(i)] \textbf{integral} if $[\omega-\psi\otimes x]\in H_{S^1}^2(\M;\Z)$ (hence $[\omega]\in H^2(\M;\Z)$);
\item[(ii)] \textbf{primitive} if it is integral and $[\omega]$ is a generator of $H^2(\M;\Z)$;
\item[(iii)] \textbf{positive} if $\cc_1$ is a positive multiple of $[\omega]$.
\end{itemize}
\end{defin}

\begin{lemma}\label{w=tau}
Let $(\M,\omega)$ be a compact symplectic manifold with a Hamiltonian
$S^1$-action and moment map $\psi\colon \M\to \R$.
If the number of fixed point is minimal, it is not restrictive to assume that
$\omega-\psi\otimes x$ is primitive and positive.
More precisely, let $\tau_1$ and $\taut_1$ be the classes defined in \eqref{basis equivariant cohomology} and \eqref{basis cohomology}.
Then we can assume that $[\omega-\psi\otimes x]=\tau_1$,
and so $[\omega]=\taut_1$ and $\cc_1=C_1[\omega]$. 
\end{lemma}
\begin{proof}
Since $H^2(\M;\Z)=\Z$ and $[\omega]\neq 0$, we can rescale the symplectic form $\omega$ in such a way
that $[\omega]$ is an integral generator of $H^2(\M;\Z)$. So $[\omega]=\pm \widetilde{\tau}_1$.
Since the kernel of the restriction map \eqref{restriction} is the ideal generated by $x$, we have
 that $[\omega-\psi\otimes x]=\pm \tau_1+P(x)$, where $P(x)$ is a constant polynomial in $x$.
 Hence, modulo translating the moment map $\psi$, we can assume that
 $[\omega-\psi\otimes x]=\pm\tau_1\in H^2_{S^1}(\M;\Z)$.
 By definition of $\tau_1$ we have that $\tau_1(P_0)=0$ and $\tau_1(P_1)=k\,x$, where $k\in \Z_{<0}$,
 and $\psi(P_0)<\psi(P_1)$. We can then conclude that $[\omega-\psi\otimes x]=\tau_1$. 
\end{proof}

The next proposition is an easy consequence of Theorem~\ref{hattori}.
\begin{prop}\label{symp hattori}
Let $(\M,\omega)$ be a compact symplectic manifold of dimension $2n$ with a Hamiltonian $S^1$-action
and moment map $\psi\colon\M\to \R$. Assume that
there is a minimal number of
fixed points
$P_0,\ldots,P_n$ and that $\omega-\psi\otimes x$ is primitive and positive.
Let $C_1$ be the positive integer defined in \eqref{C1}. Then,
\begin{itemize}
\item[(i)] there exists a quasi-ample complex line bundle $\LL$ such that $\cc_1^{S^1}(\LS)=[\omega-\psi\otimes x]$ and $\LS(P_i)=t^{-\psi(P_i)}$;
\item[(ii)] there exists $d\in \Z$ such that
\begin{equation}\label{D}
\sum_{k=1}^nw_{ik}=-C_1\psi(P_i)+d\quad\mbox{for all}\quad i=0,\ldots,n\;,
\end{equation}
and 
\begin{equation}\label{bound C1}
1\leqslant C_1 \leqslant n+1\;.
\end{equation}
\end{itemize}
\end{prop}
\begin{proof}
(i) The existence of a complex line bundle $\LL$ such that the $S^1$-action lifts to $\LL$
and $\cc_1^{S^1}(\LS)=[\omega-\psi\otimes x]\in H^2_{S^1}(\M;\Z)$ 
 is a direct consequence of Theorem~\ref{lift clb}.
It follows that $\LS(P_i)=t^{-\psi(P_i)}$, and hence $\LS$ is fine by \eqref{m.m increasing}. 
Finally, $\LL$ is quasi-ample because $\int_\M\cc_1(\LL)^n=\int_\M [\omega]^n\neq 0$.
\\
(ii) By Lemma~\ref{w=tau} we can assume that $[\omega-\psi\otimes x]=\tau_1$, thus
$$
\left( \sum_{k}w_{ik} \right)x = \ce(P_i) = C_1\tau_1(P_i)+\ce(P_0)=\left( -C_1\psi(P_i)+d \right) x\;.
$$
Then, by $(2)$ in Theorem~\ref{hattori},  we have that $1\leqslant C_1\leqslant n+1$.
\end{proof}
\begin{rmk}
The bundle  $\LL$ is usually called the \emph{pre-quantization line bundle for} $(\M,\omega)$, and its equivariant extension
$\LL^{S^1}$ the \emph{$S^1$-equivariant pre-quantization line bundle for }$(\M,\omega,\psi)$ (cf. \cite{GKS}).
\end{rmk}

Combining Theorem~\ref{hattori} with Proposition~\ref{symp hattori} we obtain the following result.
\begin{corollary}\label{equations sym}
Let $(\M,\omega)$ be as in Proposition~\ref{symp hattori} and let $\LL^{S^1}$ be the
$S^1$-equivariant pre-quantization line bundle.
For every $s=0,\ldots,n$, let $r_s(t)\in\Z[t,t^{-1}]$ be the element associated
to $\LS$, as defined in \eqref{def expl r}, and let $l_0=n+1-C_1$. 
Then
\begin{align}
\bullet&  \label{sum varphi} \int_\M[\omega]^n  =\;\;r_0(1)+r_1(1)+\cdots +r_{l_0}(1);\\
\bullet & \label{r0=1}  \,\, r_0(1)  =\;\;\td(\M)=1;\\
\bullet  &\,\, r_{s}(1)  = \;\;r_{l_0-s}(1) \quad\mbox{for all}\quad 0\leqslant s\leqslant l_0\quad\mbox{and}\quad
\label{sym r} r_s(t)=\;\;0\quad\mbox{for all}\quad s>l_0\;.
\end{align}
\end{corollary}
\begin{proof}
By Proposition~\ref{symp hattori}, we can apply Theorem~\ref{hattori} to $\LL$. Thus
\eqref{r0=1} and \eqref{sym r} 
are direct consequences of Theorem~\ref{hattori}, which also implies that
$$
\varphi_i(1)=r_0(1)+r_1(1)+\cdots +r_{l_0}(1)\,,
$$
where $\varphi_i(t)\in \Z[t,t^{-1}]$ is defined in \eqref{varphi}.
Hence, we just  have to prove  that $\varphi_i(1)=\int_\M[\omega]^n$.
Consider the commutative diagram \eqref{K commutes}.
By \eqref{varphi index} we have that
\begin{equation*}
\label{varphi(1)}\varphi_i(1)= r\left(\ind_{S^1}\left( \prod_{j\neq i}\left(1-(\LL^{S^1})^{-1}t^{a_j}\right)\right)\right)=\ind\left(1-\LL^{-1}\right)^n.
\end{equation*}
Using the Atiyah-Singer formula, and the fact that
$\ch(1-\LL^{-1})=\sum_{k=1}^{\infty}(-1)^{k+1}\frac{[\omega]^k}{k!}$, we obtain
\begin{equation*}
\varphi_i(1)
=\ind\left(1-\LL^{-1}\right)^n
=\int_{\M}\ch(1-\LL^{-1})^n\ttot(\M)=\int_\M[\omega]^n\,,
\end{equation*}
which completes the proof (cf. \cite[Lemma 3.6]{Ha}).
\end{proof}
By an argument similar to the proof of Corollary~\ref{equations sym}, 
equations
\eqref{r0=1} and \eqref{sym r} can also be turned into equations
involving the
Chern numbers of the manifold.
Namely, by \eqref{explicit r} we have that, for every $i=0,\ldots,n$,
\begin{align}
\label{r(1)}r_i(1)= & (-1)^ir\left(\ind_{S^1}\left(k_i(\LS,t)\right)\right)=(-1)^i\ind\left(k_i(\LL,1)\right)\;.
\end{align}
However, the explicit computation of $\ind(k_i(\LL,1))$ is harder to do in general.
In Theorem~\ref{m not 2} we will compute these values in the case in which $\M$
is $8$-dimensional. For that we first have to prove the following results.
\begin{lemma}\label{integral square}
Let $(\M,\omega)$ be a compact symplectic manifold of dimension $2n$
with a Hamiltonian $S^1$-action with $n+1$ fixed points. 
Let $\tau$ be an element of $H^2(\M;\Z)$.

If $n$ is \emph{even} then
$$
\quad \int_\M\tau^n=Q^2\,,
$$
where $Q$ is the unique integer such that $\tau^{n/2}=Q\taut_{n/2}$ (see Theorem~\ref{bases}).
\end{lemma}
\begin{proof}
Le $\{\taut_i\}_{i=0}^n$ be the basis of $H^*(\M;\Z)$ defined in \eqref{basis cohomology} and
 $\{\taut^\prime_i\}_{i=0}^n$ the basis of $H^*(\M;\Z)$ obtained by reversing the flow (cf. Section~\ref{reversing flow}).
Then, by \eqref{taus}, $\taut_{n/2}=\taut^\prime_{n/2}$, and hence \eqref{duality} implies that
$\int_\M \taut_{n/2}^2=\int_\M\taut_{n/2}\taut^\prime_{n/2}=1$, and the result follows.
\end{proof}
Corollary~\ref{equations sym} together with Lemma~\ref{integral square} imply the following proposition.
\begin{prop}\label{integrals}
Let $(\M,\omega)$ be as in Proposition~\ref{symp hattori}. Then,
\begin{itemize}
\item[(i)] if $C_1=n+1$ we have $\int_\M [\omega]^n=1$;\\$\;$
\item[(ii)] if $C_1=n$ we have $\int_\M [\omega]^n=2$ and $n$ is \emph{odd}.
\end{itemize}
\end{prop}
\begin{proof}
(i) If $C_1=n+1$, then by \eqref{sym r} in Corollary~\ref{equations sym}, we have $r_s(t)=0$ for all $s>0$,
and then, by \eqref{sum varphi} and \eqref{r0=1}, it follows that $\int_\M [\omega]^n=r_0(1)=1$.

(ii) If $C_1=n$, then, by Corollary \eqref{equations sym}, we have  
$r_s(1)=0$ for all $s>1$ and $r_0(1)=r_1(1)$.
Combining  \eqref{sum varphi}  and \eqref{r0=1}  we have that
$\int_\M[\omega]^n=r_0(1)+r_1(1)=2$, and, since $[\omega]\in H^2(\M;\Z)$, Lemma \eqref{integral square} implies that
$n$ is odd.
\end{proof}
Recall that, when $(\M,\omega)$ is a compact symplectic manifold of dimension $2n$ with a Hamiltonian
$S^1$-action with a minimal number of fixed points, we have
$H^i(\M;\Z)=H^i(\C P^n;\Z)$ for every $i$ (see Section~\ref{ec symplectic}).
The next proposition 
shows what is the minimal information required to compute
the ring structure of $H^*(\M;\Z)$
and the total Chern class when $\M$ is of $8$-dimensional. 

For an $8$-dimensional manifold, the total Todd class is 
\begin{equation}\label{totalTodd}
\ttot(\M)=\sum_{i=0}^4T_0^i=1+\frac{\cc_1}{2}+\frac{\cc_1^2+\cc_2}{12}+\frac{\cc_1\cc_2}{24}+\frac{-\cc_1^4+4\cc_1^2\cc_2+3\cc_2^2+\cc_1\cc_3-\cc_4}{720}\; \in H^*(\M;\Q),
\end{equation}
where $T_0^i$ is the term of degree $2i$ in $H^*(\M;\Q)$ for all $i$, and the Todd genus is given by
\begin{equation}\label{Toddgenus}
\td(\M)=\int_\M \frac{-\cc_1^4+4\cc_1^2\cc_2+3\cc_2^2+\cc_1\cc_3-\cc_4}{720}
\end{equation}
(see Section~\ref{c1cn-1}).
\begin{prop}\label{cr & Cc}
Let $(\M,\omega)$ be a compact symplectic manifold of dimension $8$, with a Hamiltonian $S^1$-action
and $5$ fixed points. 
\begin{itemize}
\item[(i)]
Let $C_1$ and $C_2$ be the constants defined in \eqref{basis equivariant cohomology},
and $l:=C_2/C_1^2$. Then $l\in \Z_{>0}$ and
$$
H^*(\M;\Z)=\Z[x_1,x_2,x_3]/(x_1^2-l\,x_2,x_1x_3-x_2^2,x_1^5,x_2^3,x_3^2,x_2x_3)\;,
$$
where $x_1,x_2,x_3$ have degrees respectively $2$, $4$ and $6$.

\item[(ii)] 
There exists $m\in \Q$ such that $\cc_2=m\,x_1^2$, and the total Chern class is given by
$$
\cc(T\M)=1+C_1x_1+(l\,m)x_2+(50/C_1)x_3+5x_1x_3\;.
$$
Moreover $C_1$, $m$ and $l$ satisfy
\begin{equation}\label{T04}
l^2(-C_1^4+4C_1^2m+3m^2)-675=0\;.
\end{equation}
\end{itemize}
\end{prop}
In particular, $H^*(\M;\Z)\simeq H^*(\C P^4;\Z)$ as rings if and only if $l=1$, and the total
Chern class $\cc(T\M)$ agrees with the one of $\C P^4$, if and only if $C_1=5$ and $l\,m=10$. 
\begin{proof}
(i) Let $x_1:=\taut_1$, $x_2:=\taut_2$ and $x_3:=\taut_3$ (see Theorem~\ref{bases}).
Then it is easy to see that $x_1^2=(C_2/C_1^2)x_2\,$. Moreover, since $C_2$ is positive,
it follows that $l\in \Z_{>0}$.
Since $H^8(\M;\Z)=\Z$, it follows from \eqref{duality}  that $x_1\,x_3=x_2^2=\taut_4$,
and $x_1^5=x_2^3=x_3^2=x_2x_3=0$ by dimensional reasons.
Moreover, using \eqref{duality}, it is  easy to see that these relations imply
$x_1^3=l^2\,x_3$ and $x_1\,x_2=l\,x_3$.

(ii) By (i) $x_1^2\neq 0$, so there exists $m\in \Q$ such that $\cc_2=m\,x_1^2=(m\,l)x_2$ (which implies in particular that $m\,l$ is an integer).
Let $\alpha\in \Z$ be such that $\cc_3=\alpha\, x_3$. Then by  \eqref{eq:c1cn-1minimal}  in Corollary~\ref{corollary1} and \eqref{duality}, we have 
$$\int_\M\cc_1\cc_3=C_1\alpha\int_\M x_1x_3=C_1\alpha=50.$$
Finally, let $\beta\in \Z$ be such that $\cc_4=\beta \,x_1x_3$. Then, by \eqref{cn} and \eqref{duality}, we have  
$$\int_\M\cc_4=\beta=5.$$
Now by \eqref{r0=1} in  Corollary~\ref{equations sym}   and \eqref{Toddgenus}, we have that
\begin{equation}\label{eqTodd}
\int_\M\frac{-\cc_1^4+4\cc_1^2\cc_2+3\cc_2^2+\cc_1\cc_3-\cc_4}{720}=\td(\M)=1
\end{equation}
and so, by \eqref{duality}, we have
$\int_\M x_1^4=l^2$, and \eqref{T04} immediately follows.
\end{proof}
We are now ready to prove the main theorem of this section.
\begin{theorem}\label{m not 2}
Let $(\M,\omega)$ be a compact symplectic manifold of dimension $8$, with a Hamiltonian $S^1$-action
with moment map $\psi\colon \M\to \R$,
and $5$ fixed points. Suppose that $[\omega-\psi\otimes x]$ is primitive and positive,
so that $\cc_1=C_1[\omega]$.

Then $C_1$ is either $1$ or $5$.
Moreover, 
the cohomology ring $H^*(\M;\Z)$ and the total Chern class $\cc(T\M)$
agree with the ones of $\C P^4$ if and only if $C_1=5$.
\end{theorem}
\begin{proof}
Let $[\omega-\psi\otimes x]=\tau_1\in H^2_{S^1}(\M;\Z)$ be the equivariant symplectic
form and $\LS$ the $S^1$-equivariant line bundle such that $\cc_1^{S^1}(\LS)=[\omega-\psi\otimes x]$
(see Lemma~\ref{w=tau} and Proposition~\ref{symp hattori}).
For every $s=0,\ldots,n$, let $k_s(\LL^{S^1},t)$ be the bundles associated to $\LL^{S^1}$ as defined in \eqref{k0} and \eqref{k ind}.
Thus we have that
\begin{align*}
k_0(\LL,1)= &\;\;1\quad\mbox{and} \\
 k_s(\LL,1)=& \binom{n+1}{s}(1-\LL)^s-\sum_{\nu=1}^s\binom{n-s+\nu}{\nu}(-\LL)^\nu k_{s-\nu}(\LL,1),
\end{align*}
yielding
\begin{align}
k_0(\LL,1)=&\;\;1\;,\nonumber\\
k_1(\LL,1)=&\;\; n+1-\LL\;,\nonumber\\
k_2(\LL,1)=& \;\;\frac{n(n+1)}{2}-(n+1)\LL+\LL^2\;,\nonumber\\
k_3(\LL,1)=& \;\;\frac{n(n^2-1)}{6}-\frac{n(n+1)}{2}\LL+(n+1)\LL^2-\LL^3\;,\nonumber\\
\label{ks}k_4(\LL,1)=&\;\; \frac{n(n^3-2n^2-n+2)}{24}+\frac{n(1-n^2)}{6}\LL+\frac{n(n+1)}{2}\LL^2-(n+1)\LL^3+\LL^4.
\end{align}

Since $\M$ is $8$-dimensional, using the fact that $\ch(\LL)=\sum_{k=0}^4\frac{[\omega]^k}{k!}$ and \eqref{ks}, it is easy to verify that
\begin{align}
\ch(k_0(\LL,1))=&\;\; 1\;,\nonumber\\
\ch(k_1(\LL,1))=&\;\; 4-[\omega]-\frac{1}{2}[\omega]^2-\frac{1}{6}[\omega]^3-\frac{1}{24}[\omega]^4\;,\nonumber\\
\ch(k_2(\LL,1))=& \;\;6-3[\omega]-\frac{1}{2}[\omega]^2+\frac{1}{2}[\omega]^3+\frac{11}{24}[\omega]^4\;,\nonumber\\
\ch(k_3(\LL,1))=& \;\; 4-3[\omega]+\frac{1}{2}[\omega]^2+\frac{1}{2}[\omega]^3-\frac{11}{24}[\omega]^4\;,\nonumber\\
\label{chern k}\ch(k_4(\LL,1))= &\;\; 1-[\omega]+\frac{1}{2}[\omega]^2-\frac{1}{6}[\omega]^3+\frac{1}{24}[\omega]^4\;.
\end{align}
Using the Atiyah-Singer formula and \eqref{r(1)}, we have
$$
r_i(1)=(-1)^i\int_\M\ch\left(k_i(\LL,1)\right) \ttot(\M)\;
$$
which, together with \eqref{chern k} and the expression of the total Todd class $\ttot(\M)$ in \eqref{totalTodd}, implies that
\begin{align*}
r_0(1)=&\;\;\int_\M T_0^4=\int_\M\frac{-\cc_1^4+4\cc_1^2\cc_2+3\cc_2^2+\cc_1\cc_3-\cc_4}{720},\\
r_1(1)=&\;\;\int_\M \left(-4\,T_0^4+\frac{[\omega]\cc_1\cc_2+[\omega]^2(\cc_1^2+\cc_2)+2[\omega]^3\cc_1+[\omega]^4}{24}\right),\\
r_2(1)=&\;\;\int_\M \left(6\,T_0^4+\frac{-3[\omega]\cc_1\cc_2-[\omega]^2(\cc_1^2+\cc_2)+6[\omega]^3\cc_1+11[\omega]^4}{24}\right),\\
r_3(1)=&\;\;\int_\M\left( -4\,T_0^4+\frac{3[\omega]\cc_1\cc_2-[\omega]^2(\cc_1^2+\cc_2)-6[\omega]^3\cc_1+11[\omega]^4}{24}\right),\\
r_4(1)=&\;\;\int_\M \left(T_0^4+\frac{-[\omega]\cc_1\cc_2+[\omega]^2(\cc_1^2+\cc_2)-2[\omega]^3\cc_1+[\omega]^4}{24}\right).
\end{align*}
Since $\omega-\psi\otimes x$ is primitive and positive, 
by Lemma~\ref{w=tau} and Proposition~\ref{cr & Cc} we have
$[\omega]=\taut_1=x_1$. Let $l\in \Z_{>0}$ and $m\in \Q$ be defined as in Proposition~\ref{cr & Cc}. Then $\int_\M[\omega]^4=l^2$.
Using \eqref{eqTodd}, $\int_\M \cc_4=5$ and
$\int_\M\cc_1\cc_3=50$ (see \eqref{cn} and \eqref{eq:c1cn-1minimal} in Corollary~\ref{corollary1}), we have that
\begin{align}
 r_1(1)=&\;-4+\frac{l^2}{24}(C_1m+C_1^2+m+2C_1+1),\nonumber\\
r_2(1)=&\;6+\frac{l^2}{24}(-3C_1m-C_1^2-m+6C_1+11),\nonumber\\
r_3(1)=&\;-4+\frac{l^2}{24}(3C_1m-C_1^2-m-6C_1+11),\nonumber\\
\label{exp_r} r_4(1)=&\;1+\frac{l^2}{24}(-C_1m+C_1^2+m-2C_1+1).
\end{align}
By  \eqref{bound C1} in Proposition~\ref{symp hattori}, we know that $1\leq C_1\leq 5$. Moreover, by \eqref{C1 divides} in Proposition~\ref{C1p}, 
$C_1$ divides $50$, so $C_1$ cannot be  $3$ nor $4$.

If $C_1=2$, Corollary~\ref{equations sym} gives
$$r_1(1)=r_2(1),\quad r_3(1)=1,\quad  r_4(1)=0$$
and $l^2=2+2r_1(1)$.
It is easy to check that, using the expressions in \eqref{exp_r}, all these conditions give the same equation, namely
\begin{equation}\label{IC2}
24+l^2-l^2m=0\;.
\end{equation}
Combining \eqref{T04} and \eqref{IC2}, we get that $m=\displaystyle\frac{97\pm\sqrt{97}}{48}$, which is impossible,
since $m$ must be rational. Hence, $C_1$ is either 1 or 5.

If $C_1=5$, Proposition~\ref{integrals} (i) implies that $\int_\M[\omega]^4=1$.
On the other hand, as we observed before,  
$$\int_\M[\omega]^4=\int_\M x^4=l^2$$
with $l\in \Z_{>0}$, and so $l=1$.
By \eqref{sym r} in Corollary~\ref{equations sym}, we have $r_s(1)=0$ for all $s>0$, and, using one of these equations 
together with the expression of $r_s(1)$ given in \eqref{exp_r},
we get $m=10$. 
By Proposition~\ref{cr & Cc} we can conclude that, if $C_1=5$, the cohomology ring and Chern
classes are standard (i.e. they agree with those of $\C P^4$).

If $C_1=1$, by Proposition~\ref{cr & Cc} (ii),  it follows immediately that the total Chern class is not standard.
In order to prove that the cohomology ring is not standard,  by Proposition~\ref{cr & Cc} (i), we need to prove that $l\neq 1$.
It is sufficient to observe that
for $l=1$ \eqref{T04} does not have any rational solutions.
\end{proof}
\begin{rmk}\label{Petrie 4}
Notice that Theorem~\ref{m not 2} also proves the Petrie conjecture when $(\M,\omega)$
is an $8$-dimensional compact symplectic manifold with a Hamiltonian $S^1$-action and 5 fixed points.
However this result is not new, see \cite{Ja}.  
\end{rmk}
\begin{rmk}\label{225}
It is natural to ask whether the equations in Corollary~\ref{equations sym} give more information when $C_1=1$. Unfortunately
 they are all identities, and the only meaningful equation is \eqref{T04}, which, in this case, is
\begin{equation}\label{C1=1}
3m^2l^2+4ml^2-l^2=675\;.
\end{equation}
However, we can find a lower bound for $\int_\M [\omega]^4=l^2$. 
In fact, it is easy to see that the first values of $l$ for which \eqref{C1=1} has rational
solutions are $l=15, 25, 40, 60$... thus implying that $\int_\M[\omega]^4\geq 225$.
In the K\"ahler case, by a Chern inequality following from the Calabi Conjecture, the only
possible values of $\int_\M[\omega]^4=\int_\M \cc_1^4$ are 225 and 625 (see \cite{W,Y}).
It would be particularly interesting to know whether one can get an upper bound using
symplectic techniques (see also Remark~\ref{upper bound c1^n}).
\end{rmk}
\begin{rmk}
When $\M$ is $6$-dimensional we have,
$1\leq C_1\leq 4$
and,
by Proposition~\ref{integrals},
\begin{itemize}
\item[$\bullet$] if $C_1=4$, then$\int_\M\cc_1^3=64$;
\item[$\bullet$] if $C_1=3$,  then $\int_\M\cc_1^3=54$.
\end{itemize}
Moreover, by \eqref{cn} and \eqref{eq:c1cn-1minimal} in Corollary~\ref{corollary1},  we have $\int_\M\cc_1\cc_2=24$ and $\int_\M\cc_3=4$. 

Hence, if $C_1=4$, the Chern classes are standard, i.e. they agree with the ones of $\C P^3$ and
when $C_1=3$ they agree with the ones of $Gr_2^+(\R^5)$, the Grassmannian of oriented
$2$-planes in $\R^5$.

However, when $C_1$ is either $1$ or $2$, the equations given by Corollary~\ref{equations sym}
are all identities.
\end{rmk}

\section{Minimal number of fixed points: classification results}\label{classification results}
Let $(\M,\omega)$ be a compact symplectic manifold of dimension $2n$ equipped with a Hamiltonian $S^1$-action with a minimal number of fixed points $P_0, P_1, \ldots, P_n$, with $\lambda(P_i)=i$. We will now apply our algorithms towards a classification of these actions. Assuming that the $S^1$-action satisfies $(\mathcal{P}_0^+)$, we use Algorithms~\ref{alg:5}(B) and \ref{alg:6}(B) according to the existence of an edge from $P_0$ to $P_1$ and/or an edge from $P_{n-1}$ to $P_n$, and we obtain a list that necessarily contains all possible isotropy weights. Then we take into account  several simple properties satisfied by the isotropy weights in order to reduce the number of possibilities:
\begin{itemize}
\item  at each fixed point, the isotropy weights must be coprime integers;
\item at the fixed point $P_i$ the first $i$ weights are negative and the others are positive;
\item by \eqref{C1} and \eqref{symmetries} we must have
\begin{equation}\label{eq:equal}
\frac{\ce(P_1)-\ce(P_0)}{\Lambda_1^-}=\frac{\ce(P_{n-1})-\ce(P_n)}{\Lambda_{n-1}^+},
\end{equation}
where $\Lambda_1^-=w_{11}x$, with $w_{11}$ the unique negative weight at $P_1$, and $\Lambda_{n-1}^+=w_{n-1\,n} x$, with $w_{n-1\,n}$ the unique positive weight at $P_{n-1}$; moreover, the number in \eqref{eq:equal} must be a positive divisor of $\frac{1}{2}n(n+1)^2$ smaller or equal to $n+1$; note that \eqref{eq:equal} can be rewritten in terms of the isotropy weights as
\begin{equation*}
\frac{\sum_{j=1}^n w_{1j} - \sum_{j=1}^n w_{0j} }{w_{11}}=\frac{\sum_{j=1}^n w_{n-1\, j} - \sum_{j=1}^n w_{n\,j} }{w_{n-1\,n}}.
\end{equation*}
\item when $\dim \M=8$ we know by Theorem~\ref{thm dim8} that the number in \eqref{eq:equal}  must be equal to $1$ or $5$; 
\item the isotropy weights must satisfy the equations in \eqref{localized chern}.
\end{itemize}
Moreover, whenever the graph has multiple edges, there are other simple properties that the isotropy weights must satisfy. These are summarized in the following technical lemmas.
\begin{lemma}\label{mecycle}
Let $S\subset E$ be a set of multiple edges between two fixed points $P$ and $Q$ (i.e. $\ii(e)=P$ and  $\te(e)=Q$ for every $e\in S$) with $\vert S \rvert=n-2$,
and assume that $n>3$. For $F\in \{P,Q\}$ let $E_F=E_{F,\ii}\cup E_{F,\te}$ be the set of edges  such that either  $\ii(e)=F$ or $\te(e)=F$. If there is an $F\in  \{P,Q\}$ such that $E_F\setminus S\subset E^{\circlearrowleft}$ then the isotropy weights corresponding to the multiple edges in $S$ must be coprime.
\end{lemma}
\begin{proof}
Define $E_F^{\circlearrowleft}$ to be $E_F\cap E^{\circlearrowleft}$ and
let us  assume, without loss of generality, that $E_P\setminus S =E_P^{\circlearrowleft}$. If 
$$\gcd_{e\in S} \{ \lvert w(e)\rvert \}=k>1$$ 
and $r$ is the absolute value of the weight corresponding to the cycle $e\in E_P^{\circlearrowleft}$, then, since the action is effective, we must have $\gcd\{k,r\}=1$. Consequently, the isotropy submanifold fixed by $\Z_k$ is a $2(n-2)$-submanifold with an effective $S^1\cong S^1/\Z_k$-Hamiltonian action with only two fixed points which is  impossible.
\end{proof}
\begin{lemma}\label{lemma:8.2}
Let $S\subset E$ be a set of multiple edges between two fixed points $P$ and $Q$ and let $\ell=\lvert S \rvert$ and $n>2$. 
\begin{enumerate}
\item If $\ell \geq n-1$ then   $\gcd_{e\in S}\{ \lvert w(e)\rvert\}=1$.
\item For every integer $2\leq r\leq \ell$ and every subset $\widetilde{S}\subset S$ with $\lvert \widetilde{S} \rvert=r$, there exist edges $e_1\in E_P\setminus \widetilde{S}$ and $e_2\in E_Q \setminus \widetilde{S}$ whose isotropy weights are multiples of 
$g:=\gcd_{e\in \widetilde{S}} \{ \lvert w(e)\rvert \}$.
\end{enumerate}
\end{lemma}
\begin{proof}
If $\ell=n$ then, since the action is effective, we must have $\gcd_{e\in S}\{ \lvert w(e)\rvert\}=1$. If $\ell=n-1$ and  $\gcd_{e\in S}\{ \lvert w(e)\rvert\}=k>1$ then, denoting by $r$ the absolute value of the weight corresponding to the edge in $E_P\setminus S$, we have $\gcd\{k,r\}=1$ since the action is effective. Then the isotropy submanifold fixed by $\Z_k$ is a $2(n-1)$-submanifold
 with an effective $S^1\cong S^1/\Z_k$ Hamiltonian action with only two fixed points which is  impossible and so we must have $k=1$.

To prove $(2)$ we see that if there exists a subset $\widetilde{S} \subset S$ for which there is no edge in $E_P\setminus \widetilde{S}$ with weight a multiple of $g:=\gcd_{e\in \widetilde{S}} \{ \lvert w(e)\rvert \}$, then there is an isotropy submanifold of dimension $2\lvert \widetilde{S} \rvert$ fixed by  $\Z_g$ with only two fixed points which is impossible. Similarly, we conclude the same for $Q$.
\end{proof}

We first run Part I of the  \texttt{Mathematica} file \texttt{MinimalIW.nb} to generate the list of possible non-negative multigraphs. Then we run Part II of this file to produce the list of the determinants of the matrices $A(\Gamma_i)-\diag(\m(E_i))$ in \eqref{eq:Matricesi}. Note that, when the graph has more than one connected component, this file yields the sum of the squares of determinants of these matrices since this sum is zero if and only if all the determinants are zero. Then we run the  \texttt{C++} files  \texttt{NewPartitions$\#$.cpp}, where $\#$ denotes the number of the multigraph in the list above, to obtain partitions of  $\frac{1}{2}n(n+1)^2$ into $n(n+1)/2$ non-negative numbers $\m(e)$ (according to Algorithms~\ref{alg:5}(B) and \ref{alg:6}(B), depending on the existence of an edge from $P_0$ to $P_1$ and/or an edge from $P_{n-1}$ to $P_n$), for which all the determinants of  the matrices $A(\Gamma_i)-\diag(\m(E_i))$ are zero. Then we run Part III of   \texttt{MinimalIW.nb}  to sort these partitions according to the rank of the matrices, dividing them into two sets: those that originate matrices of rank $\lvert E \rvert -1$ and those that originate matrices of lower rank. In the first case, Part III of the file  \texttt{MinimalIW.nb}  also selects  those partitions that  originate matrices with nullspaces intersecting $\Z_{>0}^{\lvert E \rvert}$.  Part IV of   \texttt{MinimalIW.nb} considers the first set of partitions,  producing  a list of the corresponding  isotropy weights 
$$
\w(E) \in \Big(  \ns\left(A(\Gamma)-\diag\left(\m(E)\right)\right)\Big)\cap \Z_{>0}^{\lvert E \rvert},
$$ 
and checking if they can satisfy  the polynomial equations  in \eqref{localized chern}. 
 Part V of \texttt{MinimalIW.nb} considers the second set of partitions. It begins by selecting those that yield matrices with nullspaces intersecting  $\Z_{>0}^{\lvert E \rvert}$ (Part V a.), and producing the list of the corresponding isotropy weights (Part V b.). Then it selects those that possibly verify the properties listed in the beginning of this section as well as Lemmas~\ref{mecycle} and \ref{lemma:8.2}. At each step, the resulting lists of isotropy weights are saved in different files so that it is easy to verify which isotropy weights are discarded at each test performed. All the relevant  files can be downloaded from \texttt{http://www.math.ist.utl.pt/$\sim$lgodin/MinimalActions.html}.

In the following we list the results obtained when the dimension of $\M$ is $4$, $6$ or $8$. Note that $\texttt{b}[\,i\,]$ denotes a positive integer for every $i$.

\subsection{Dimension $4$}\label{dim4}
\begin{figure}[h]\label{graphsdim4}
\begin{center}
\includegraphics[scale=0.6] {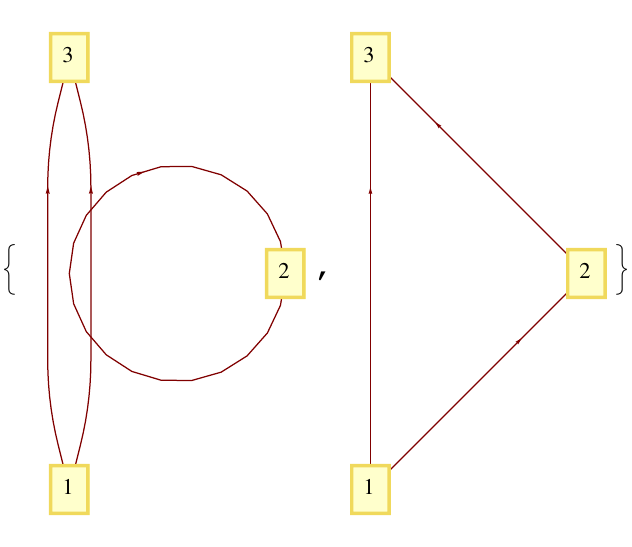}
\caption{Dimension $4$}
\end{center}
\end{figure}

When $\M$ is $4$-dimensional, it is easy to see that the only two possible multigraphs
that can arise are those in Figure~\ref{graphsdim4}.
Thus, by Proposition~\ref{spm}, they are  non-negative. 
For the first multigraph we see that the weights at the point of index 2 are $1$ and $-1$. Indeed, since the action is effective,
if this were not the case we would necessarily have the second multigraph.
Moreover, by Corollary~\ref{abbv discrete} applied to $\mu=1$, it is easy to see that the weights at the minimum must be
$1$ and $2$, and that the ones at the maximum are $-1$ and $-2$. Thus this set of weights can also be obtained with the second multigraph.

\begin{figure}[h]
\begin{center}
\includegraphics[scale=.9] {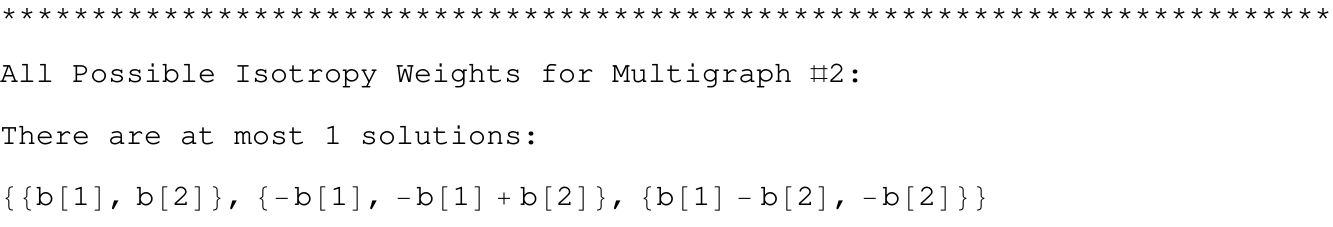}
\caption{Possible weights in dimension $4$}
\label{wdim4}
\end{center}
\end{figure}

Therefore we can run Algorithm~\ref{alg:6}(B), obtaining the  set of weights in Figure~\ref{wdim4}.

This is  the same multiset of isotropy weights as in the standard $S^1$-action on $\C P^2$ described in Example~\ref{cpn}, with \texttt{b}$[\,i\,]=\xi_0-\xi_i$, $i=1,2$.
Hence, by Theorem~\ref{bases}, the (equivariant) cohomology ring and Chern classes of the manifold are the same as the one of $\C P^2$
(with this $S^1$-action).
 Note that this result agrees with the one obtained by Karshon's classification for $4$-dimensional $S^1$-Hamiltonian manifolds \cite[Section 6.3]{K}.
  
\subsection{Dimension $6$}
\begin{figure}[h!]
\begin{center}
\includegraphics[scale=.65] {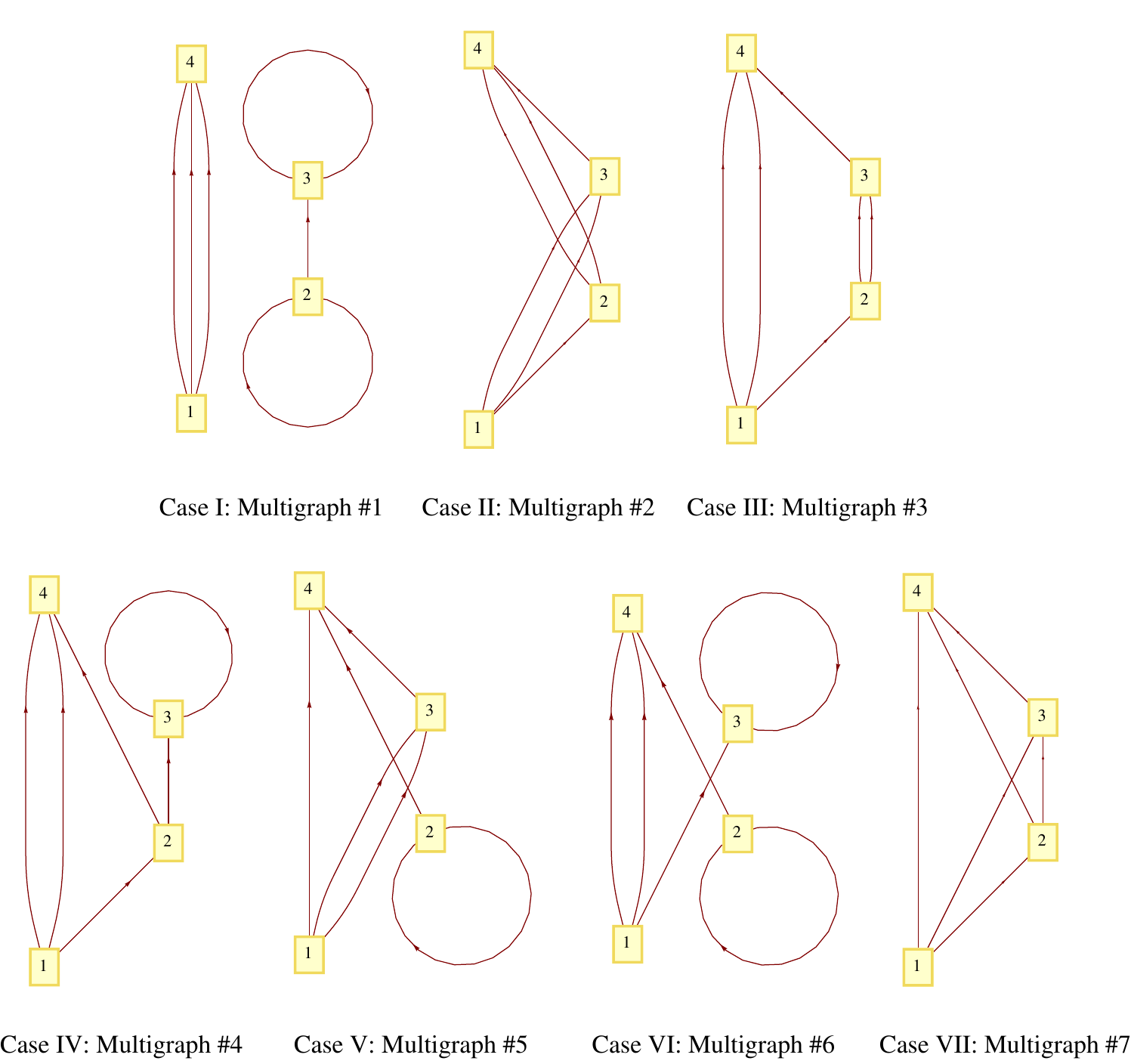}
\caption{Dimension $6$}
\label{graphs6}
\end{center}
\end{figure}

Let $(\M,\omega)$ be a six dimensional compact symplectic manifold with an $S^1$-Hamiltonian action and fixed points
$P_0,\ldots,P_3$, with $\lambda(P_i)=i$.
Running  Algorithms~\ref{alg:5}(B) and \ref{alg:6}(B) for the $7$ non-negative multigraphs in Figure~\ref{graphs6}, we obtain that all of them  may, in principle, admit possible solutions. However, as we will see next, all the solutions can be obtained by considering only the multigraphs with no cycles $\# 3$ and  $\# 7$ from Cases $3.$ and $7.$ below. Note that Tolman in \cite{T1} also rules out the existence of multigraphs with cycles  so we could have run our algorithm only for positive multigraphs. We opted to consider all non-negative multigraphs in order to show that our methods also rule out the existence of weights specific to multigraphs with cycles.

\begin{figure}[h!]
\begin{center}
\includegraphics[scale=0.8] {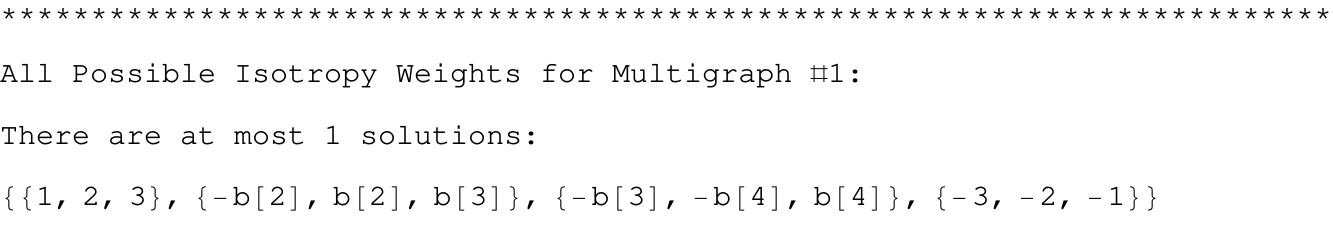}
\caption{Possible weights for Case 1}
\label{wc1}
\end{center}
\end{figure}

\textbf{Case 1:} 
For the multigraph in Case $1$ of Figure~\ref{graphs6} the weights given by our algorithm are the ones in Figure~\ref{wc1}. Using the ABBV Localization formula, with $\mu=1$ and $\mu=\cc_1^{S^1}$  and this  set of weights, we obtain
\begin{align*}
0=\int_\M 1& = \sum_{i=0}^3 \frac{1}{\prod_{j=1}^3 w_{ij}}= \frac{1}{6} - \frac{1}{\texttt{b}[2] ^2\texttt{b}[3] }+\frac{1}{\texttt{b}[3]\texttt{b}[4]^2 } -\frac{1}{6} =\frac{ \texttt{b}[2]^2-\texttt{b}[4]^2}{\texttt{b}[2] ^2\texttt{b}[3] \texttt{b}[4] ^2}
\end{align*}
and
\begin{align*}
0=\int_\M \cc_1^{S^1} & =  \sum_{i=0}^3\left( \frac{\sum_{j=1}^3 w_{ij} }{\prod_{j=1}^3 w_{ij}}\right)=1 - \frac{1}{\texttt{b}[2] ^2}-\frac{1}{\texttt{b}[4]^2 }  +1 =2-\frac{\texttt{b}[2] ^2+\texttt{b}[4]^2 }{\texttt{b}[2] ^2 \texttt{b}[4] ^2},
\end{align*}
and so $\texttt{b}[2]=\texttt{b}[4]=1$. Then, using the fact that
$$
\frac{\ce(P_1)-\ce(P_0)}{\Lambda_1^-}=6-\texttt{b}[3]
$$
must be a divisor of $24$ no larger that $4$, we conclude that $\texttt{b}[3]=2,3,4$ or $5$. If $\texttt{b}[3]=2$, the resulting multiset of weights is a particular case of the one in Case $7$ II. Note that it is precisely the set of weights of the $S^1$-action on $\C P^3$ described in Example~\ref{cpn} with $\xi_0-\xi_1=1$, $\xi_0-\xi_2=2$ and $\xi_0-\xi_3=3$. If $\texttt{b}[3]=3$, the resulting multiset of weights is a particular case of the one obtained in Case $7$ I. It is precisely the set of weights of the $S^1$-action on $Gr_2^+(\R^5)$ described in Example~\ref{gr}, by taking $\xi_0=2$ and $\xi_1=1$. If $\texttt{b}[3]=4$, the multiset of weights obtained is precisely \eqref{v5} of the $S^1$-action described in Example~\ref{fano} for the Fano manifold $V_5$. Finally, if $\texttt{b}[3]=5$, the multiset of weights obtained is precisely  \eqref{v22} of the $S^1$-action described in Example~\ref{fano} for the Fano manifold $V_{22}$.

\begin{figure}[h!]
\begin{center}
\includegraphics[scale=0.8] {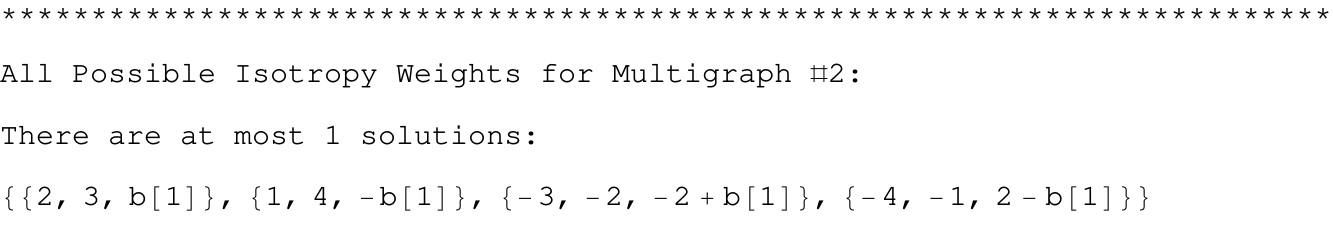}
\caption{Possible weights for Case 2}
\label{wc2}
\end{center}
\end{figure}

\textbf{Case 2:} 
For the multigraph in Case $2$ of Figure~\ref{graphs6} the weights given by our algorithm are the ones in Figure~\ref{wc2}. Using the ABBV Localization formula, with $\mu=1$ and this  set of weights, we obtain
$$
0=\int_\M 1= \frac{1}{6\texttt{b}[1] } - \frac{1}{4\texttt{b}[1] }+\frac{1}{6(\texttt{b}[1]-2)}  - \frac{1}{4(\texttt{b}[1] -2)} =-\frac{\texttt{b}[1]-1}{6\texttt{b}[1](\texttt{b}[1]-2) },
$$
and so $\texttt{b}[1] =1$. Since, on the other hand, the weight $\texttt{b}[1]-2$ at $P_2$ must be positive, we conclude that this case is impossible.

\begin{figure}[h!]
\begin{center}
\includegraphics[scale=0.8] {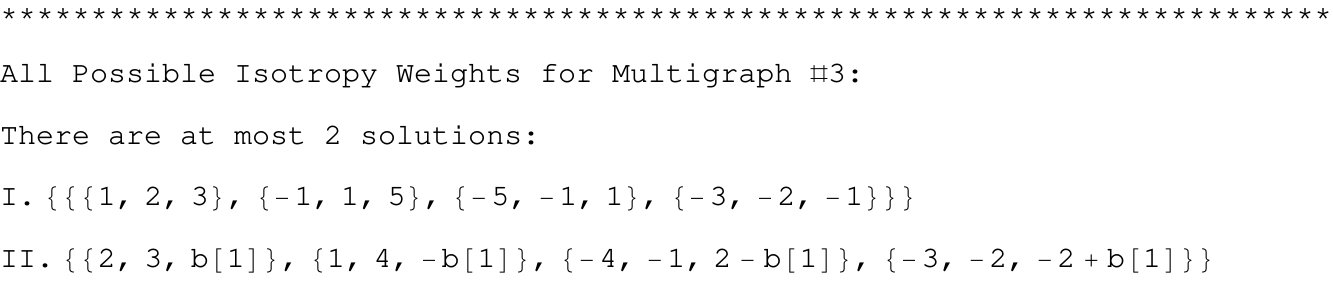}
\caption{Possible weights for Case 3}
\label{wc3}
\end{center}
\end{figure}

\textbf{Case 3:} For the multigraph in Case $3$ of Figure~\ref{graphs6} the weights given by our algorithm are the ones in Figure~\ref{wc3}.

{\bf I.} The set of weights in Figure~\ref{wc3} I.  is precisely the set of weights \eqref{v22} of the $S^1$-action described in Example~\ref{fano} for the Fano manifold $V_{22}$.

{\bf II.}  Using the ABBV Localization formula, with $\mu=1$ and the set of weights in Figure~\ref{wc3} II. we have
$$
0=\int_\M 1= \frac{1}{6\texttt{b}[1] } - \frac{1}{4\texttt{b}[1] }+\frac{1}{4(2-\texttt{b}[1])}  - \frac{1}{6(2-\texttt{b}[1] )} =-\frac{1-\texttt{b}[1]}{6\texttt{b}[1](2-\texttt{b}[1]) },
$$
and so $\texttt{b}[1] =1$. Then this is precisely the set of weights \eqref{v5} of the $S^1$-action described in Example~\ref{fano} for the Fano manifold $V_5$.

\begin{figure}[h!]
\begin{center}
\includegraphics[scale=0.8] {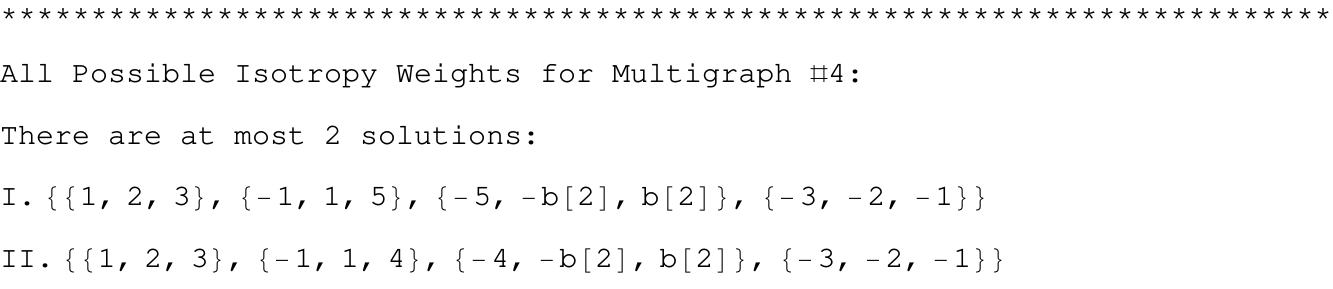}
\caption{Possible weights for Case 4}
\label{wc4}
\end{center}
\end{figure}

\vspace{.3cm}
\textbf{Case 4:} For the multigraph in Case $4$ of Figure~\ref{graphs6} the weights given by our algorithm are the ones in Figure~\ref{wc4}. Using the ABBV Localization formula, with $\mu=1$ and this set of weights, we obtain that $\texttt{b}[2] =1$ in both cases and so we get the same sets of weights as in Cases $3$ I. and $3$ II. respectively.

\begin{figure}[h!]
\begin{center}
\includegraphics[scale=0.8] {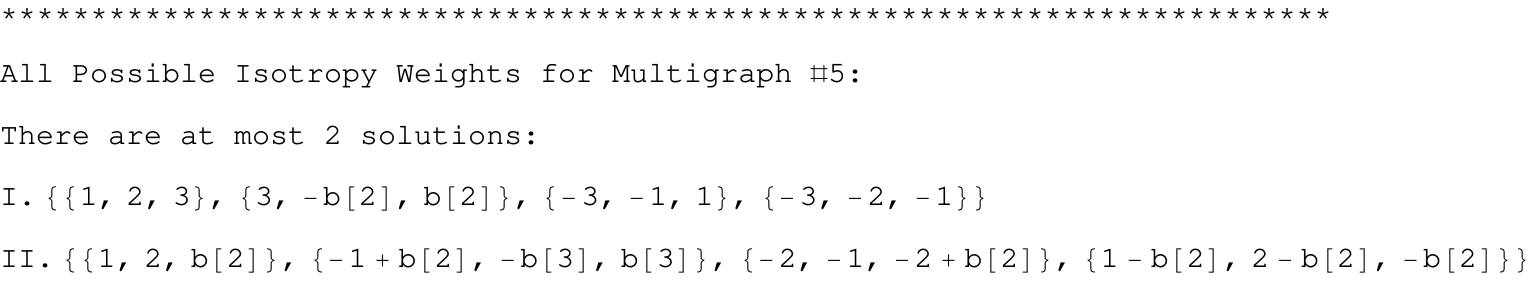}
\caption{Possible weights for Case $5$}
\label{wc5}
\end{center}
\end{figure}

\vspace{.3cm}
\textbf{Case $5$:} For the multigraph in Case $5$ of Figure~\ref{graphs6} the weights given by our algorithm are the ones in Figure~\ref{wc5}. 

{\bf I.} Using the ABBV Localization formula, with $\mu=1$ and the set of weights in Figure~\ref{wc5} I., we obtain that $\texttt{b}[2] =1$.
Then the multiset of weights is a particular example of Case $7$ I. 

{\bf II.}  Using the ABBV Localization formula, with $\mu=1$ and the set of weights in Figure~\ref{wc5} II., we obtain that $\texttt{b}[3] =1$. Then this set of weights is a particular example of Case $7$ II. 

\begin{figure}[h!]
\begin{center}
\includegraphics[scale=0.8] {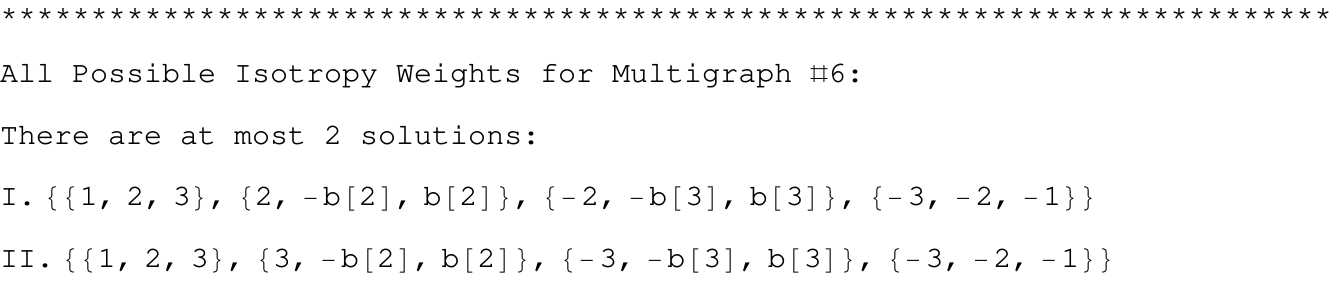}
\caption{Possible weights for Case 6}
\label{wc6}
\end{center}
\end{figure}

\vspace{.3cm}
\textbf{Case 6:} For the multigraph in Case $6$ of Figure~\ref{graphs6} the weights given by our algorithm are the ones in Figure~\ref{wc6}. Using the ABBV Localization formula, with $\mu=1$ and with $\mu=\cc_1^{S^1}$ we obtain that $\texttt{b}[2] =\texttt{b}[3] =1$ in both cases. 

Then,

{\bf I.}  the set of weights in Figure~\ref{wc6} I. is a particular example of Case $7$ II. 

{\bf II.} the set of weights in Figure~\ref{wc6} II. is a particular example of Case $7$ I.

\begin{figure}[h!]
\begin{center}
\includegraphics[scale=0.8] {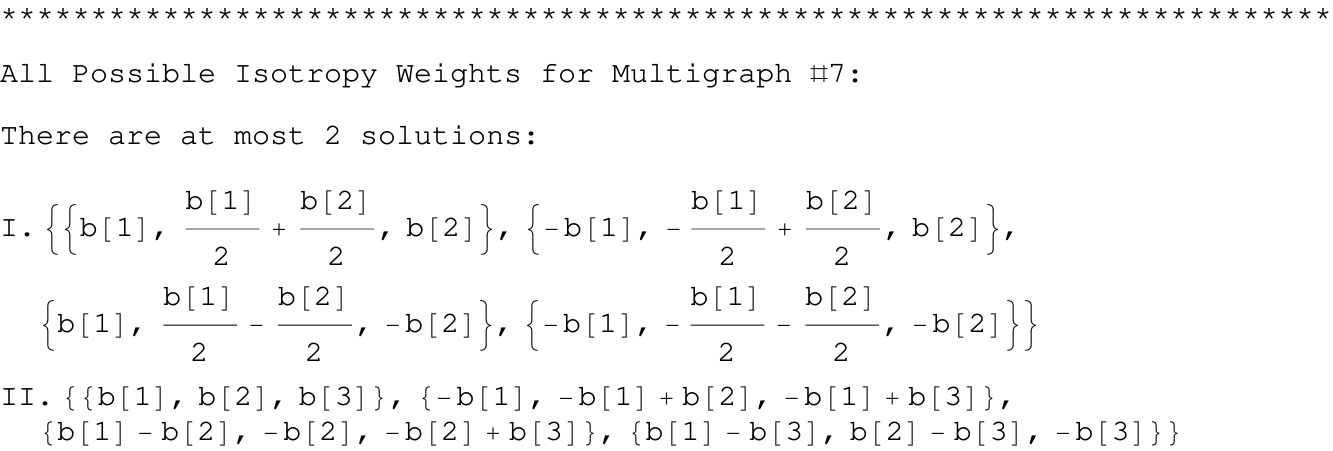}
\caption{Possible weights for Case 7}
\label{wc7}
\end{center}
\end{figure}

\vspace{.3cm}
\textbf{Case 7:} For the multigraph in Case $7$ of Figure~\ref{graphs6} the weights given by our algorithm are the ones in Figure~\ref{wc7}.

{\bf I.} The set of weights in Figure~\ref{wc7} I. is precisely the set of weights of the $S^1$-action on $Gr_2^+(\R^5)$ described in Example~\ref{gr}, by taking $\xi_0=\frac{\text{\texttt{b}$[\,1\,]$+\text{\texttt{b}$[\,2\,]$}}}{2}$ and $\xi_1=\frac{\text{\texttt{b}$[\,2\,]$-\text{\texttt{b}$[\,1\,]$}}}{2}$.

{\bf II.} The set of weights in Figure~\ref{wc7} II. is precisely the set of weights of the $S^1$-action on $\C P^3$ described in Example~\ref{cpn} with
\texttt{b}$[\,i\,]=\xi_0-\xi_i$, $i=1,2,3$.

\subsection{Dimension $8$.}\label{class dim 8}

\begin{figure}[h!]
\begin{center}
\includegraphics[scale=.8]{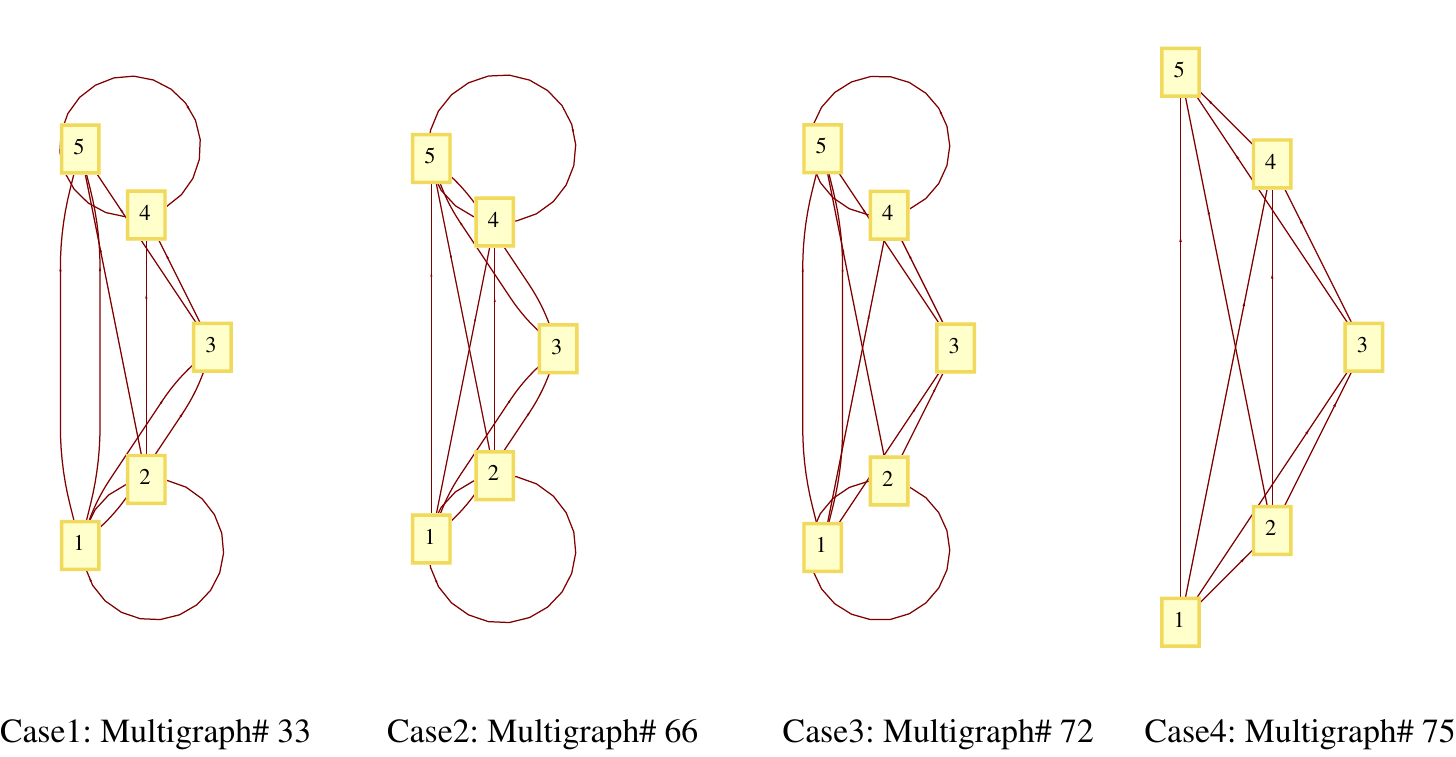}
\caption{Dimension 8}
\label{graphs8}
\end{center}
\end{figure}

Running Algorithms~\ref{alg:5}(B) and \ref{alg:6}(B) for the $75$ non-negative multigraphs on an $8$-dimensional manifold with $5$ fixed points, we obtain that only  the multigraphs in Figure~\ref{graphs8} may, in principle, admit possible solutions. However, as we will see next they can all be easily ruled out except for multigraph$\#$ $75$ considered in Case $4$. below.

\begin{figure}[h!]
\begin{center}
\includegraphics[scale=.8]{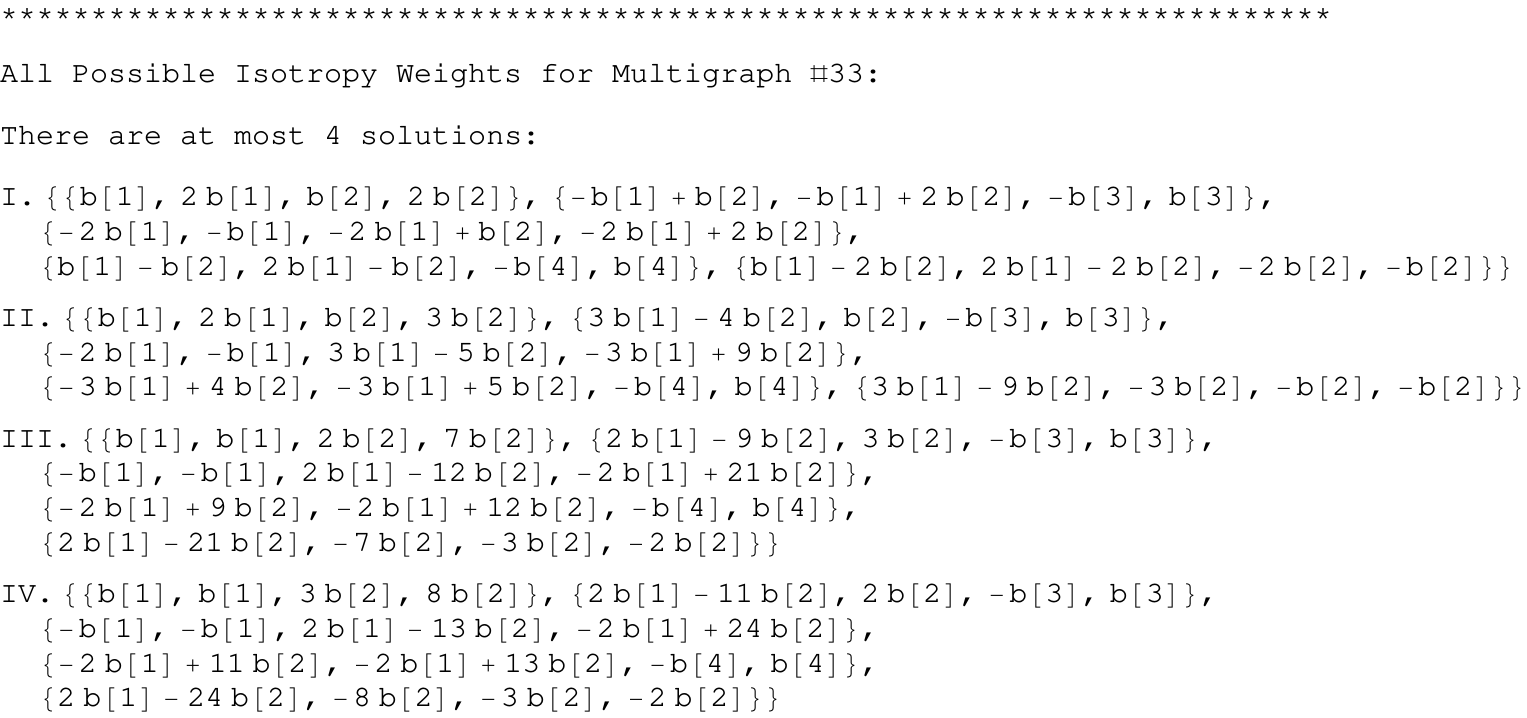}
\caption{Possible isotropy weights for Case 1.}
\label{wc81}
\end{center}
\end{figure}

\textbf{Case $1$:} For the multigraph in Case $1$ of Figure~\ref{graphs8} the weights given by our algorithm are the ones in Figure~\ref{wc81}. Note that this multigraph has two pairs of multiple  edges: two edges from $P_0$ to $P_4$ and two edges from $P_0$ to $P_2$.

{\bf I.} For the set of weights in Figure~\ref{wc81} I. we see that the weights corresponding to the two edges from $P_0$ to $P_2$ are $\texttt{b}[1]$ and $2\texttt{b}[1]$, and the ones corresponding to the two edges from $P_0$ to $P_4$ are  $\texttt{b}[2]$ and $2\texttt{b}[2]$. Since the action is effective, we have $\gcd\left\{\texttt{b}[1],\texttt{b}[2]\right\}=1$. 
Then, by Lemma~\ref{lemma:8.2}  $(2)$, we either have
$$
\texttt{b}[1]=1,\quad \text{or} \quad  \gcd\left\{\texttt{b}[1],2\texttt{b}[2]\right\}=\texttt{b}[1],
$$
and so
$$
\texttt{b}[1]=1, \quad \text{or} \quad  \texttt{b}[1]= 2.
$$
Similarly, using the two edges from $P_0$ to $P_4$, we conclude that 
$$
\texttt{b}[2]=1, \quad \text{or} \quad  \texttt{b}[2]= 2.
$$
Since the weight $\texttt{b}[2]-2\texttt{b}[1]$ at $P_2$ must be positive this is impossible.

{\bf II.} For the set of weights in Figure~\ref{wc81} II. we see that  the weights corresponding to the two edges from $P_0$ to $P_2$ are $\texttt{b}[1]$ and $2\texttt{b}[1]$, and the ones corresponding to the two edges from $P_0$ to $P_4$ are  $\texttt{b}[2]$ and $3\texttt{b}[2]$. Since the action is effective, we have $\gcd\left\{\texttt{b}[1],\texttt{b}[2]\right\}=1$. Moreover, by Lemma~\ref{lemma:8.2}  $(2)$, we either have
$$
\texttt{b}[1]=1,\quad \text{or} \quad  \gcd\left\{\texttt{b}[1],3\texttt{b}[2]\right\}=\texttt{b}[1],
$$
and so
$$
\texttt{b}[1]=1, \quad \text{or} \quad  \texttt{b}[1]= 3.
$$
Similarly, using the two edges from $P_0$ to $P_4$ we conclude that we either have
$$
\texttt{b}[2]=1, \quad \text{or} \quad  \texttt{b}[2]= 2.
$$
Since the weights $3\texttt{b}[1]-5\texttt{b}[2]$ and $9\texttt{b}[2]-3\texttt{b}[1]$ at $P_2$ must be positive this is impossible.

{\bf III.}   For the set of weights in Figure~\ref{wc81} III. we see that the weights corresponding to the two edges from $P_0$ to $P_2$ are both $\texttt{b}[1]$, and the ones corresponding to the two edges from $P_0$ to $P_4$ are  $2\texttt{b}[2]$ and $7\texttt{b}[2]$. Since the action is effective, we have $\gcd\left\{\texttt{b}[1],\texttt{b}[2]\right\}=1$. Moreover, by Lemma~\ref{lemma:8.2}  $(2)$, we either have
$$
\texttt{b}[1]=1,\quad \text{or} \quad  \gcd\left\{\texttt{b}[1],2\texttt{b}[2]\right\}=\texttt{b}[1], \quad \text{or} \quad \gcd\left\{\texttt{b}[1],7\texttt{b}[2]\right\}=\texttt{b}[1],
$$
and so
$$
\texttt{b}[1]=1, \quad \text{or} \quad  \texttt{b}[1]= 2,  \quad \text{or} \quad  \texttt{b}[1]= 7 .
$$
Similarly, using the two edges from $P_0$ to $P_4$ we conclude that we must have
$$
\texttt{b}[2]=1,\quad \text{or} \quad  \gcd\left\{2\texttt{b}[2],\texttt{b}[1]\right\}=2\texttt{b}[2],\quad \text{or} \quad  \gcd\left\{7\texttt{b}[2],\texttt{b}[1]\right\}=7\texttt{b}[2].
$$
and so the only possibility is $\texttt{b}[2]=1$.
Since the weight $2\texttt{b}[1]-12\texttt{b}[2]$ at $P_2$ must be positive the only possibility is  $\texttt{b}[1]=7$ and  $\texttt{b}[2]=1$. In this case the (connected) set of points fixed by $\Z_7$ would contain $P_0$ and $P_4$. However, if   $\texttt{b}[1]=7$ and $\texttt{b}[2]=1$, the subgroup $\Z_7$ acts trivially on a $3$-dimensional complex subspace of $T_{P_0}\M$ and on a $2$-dimensional complex subspace at $T_{P_4}\M$ (since the  weights at $P_4$ would be $\{-2,-3,-7,-7\}$), which is impossible.

{\bf IV.}  For the set of weights in Figure~\ref{wc81} IV. we see that   the weights corresponding to the two edges from $P_0$ to $P_2$ are both $\texttt{b}[1]$, and the ones corresponding to the two edges from $P_0$ to $P_4$ are  $3\texttt{b}[2]$ and $8\texttt{b}[2]$. Since the action is effective, we have $\gcd\left\{\texttt{b}[1],\texttt{b}[2]\right\}=1$. Moreover, by Lemma~\ref{lemma:8.2}  $(2)$, we either have
$$
\texttt{b}[1]=1,\quad \text{or} \quad  \gcd\left\{\texttt{b}[1],3\texttt{b}[2]\right\}=\texttt{b}[1], \quad \text{or} \quad \gcd\left\{\texttt{b}[1],8\texttt{b}[2]\right\}=\texttt{b}[1],
$$
and so
$$
\texttt{b}[1]\in \{1, \,2,\, 3,\,4,\,8\}.
$$
Similarly, using the two edges from $P_0$ to $P_4$ we conclude that we must have
$$
\texttt{b}[2]=1,\quad \text{or} \quad  \gcd\left\{3\texttt{b}[2],\texttt{b}[1]\right\}=3\texttt{b}[2],\quad \text{or} \quad  \gcd\left\{8\texttt{b}[2],\texttt{b}[1]\right\}=8\texttt{b}[2],
$$
and so the only possibility is $\texttt{b}[2]=1$.
Since the weight $2\texttt{b}[1]-13\texttt{b}[2]$ at $P_2$ must be positive the only possibility is  $\texttt{b}[1]=8$ and  $\texttt{b}[2]=1$. 
In this case the (connected) set of points fixed by $\Z_8$ would contain $P_0$ and $P_4$. However, if  $\texttt{b}[1]=8$  and  $\texttt{b}[2]=1$, the subgroup $\Z_8$ acts trivially on a $3$-dimensional complex subspace of $T_{P_0}\M$ and on a $2$-dimensional complex subspace at $T_{P_4}\M$ (since the  weights at $P_4$ would be $\{-2,-3,-8,-8\}$), which is impossible.

We conclude that Case $1$. is impossible.

\begin{figure}[h!]
\begin{center}
\includegraphics[scale=.8]{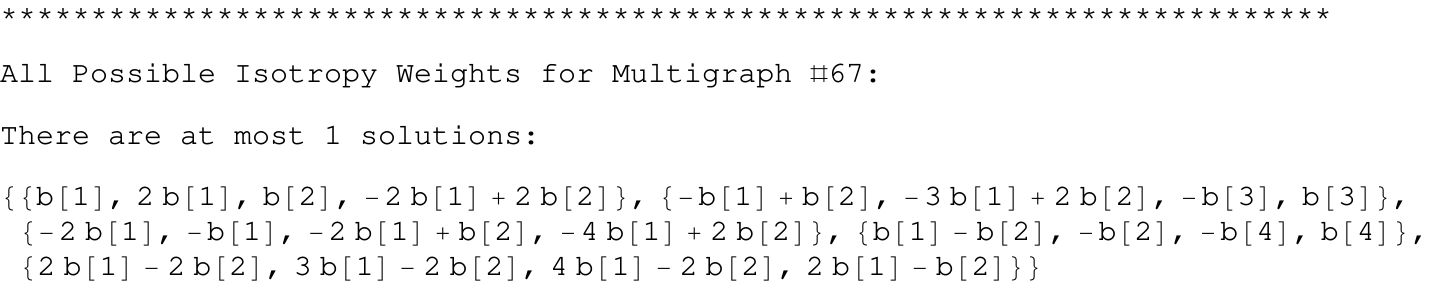}
\caption{Possible isotropy weights for Case $2$.}
\label{wc82}
\end{center}
\end{figure}

\textbf{Case $2$.}  For the multigraph in Case $2$ of Figure~\ref{graphs8} the weights given by our algorithm are the ones in Figure~\ref{wc82}. First we see that there exist two pairs of multiple  edges on this multigraph: two edges from $P_0$ to $P_2$ and two edges from $P_2$ to $P_4$. The weights corresponding to the first two  edges  are $\texttt{b}[1]$ and $2\texttt{b}[1]$, while the ones  corresponding to the last are  $\texttt{b}[2]-2\texttt{b}[1]$ and $2(\texttt{b}[2]-2\texttt{b}[1])$. Since the action is effective, we see from the weights at $P_0$ that we  need $\gcd\left\{\texttt{b}[1],\texttt{b}[2]\right\}=1$. Moreover, by Lemma~\ref{lemma:8.2}  $(2)$, we either have
$$
\texttt{b}[1]=1,\quad \text{or} \quad  \gcd\left\{\texttt{b}[1],2(\texttt{b}[2]-2\texttt{b}[1])\right\}=\texttt{b}[1],
$$
and so
$$
\texttt{b}[1]=1, \quad \text{or} \quad  \texttt{b}[1]= 2.
$$
Similarly, using the other two edges from $P_2$ to $P_4$, we have
$$
l=1, \quad \text{or} \quad  \gcd\left\{l, \texttt{b}[1]\right\}=l, \quad \text{or} \quad  \gcd\left\{l,2\texttt{b}[1]\right\}=l,
$$
with $l=\texttt{b}[2]-2\texttt{b}[1]$ and so, since
$$
\gcd\left\{l, \texttt{b}[1]\right\}=\gcd\left\{\texttt{b}[2]-2\texttt{b}[1],\texttt{b}[1]\right\}=\gcd\left\{\texttt{b}[2],\texttt{b}[1]\right\}=1
$$
and 
$$
\gcd\left\{l, 2\texttt{b}[1]\right\}=\gcd\left\{\texttt{b}[2]-2\texttt{b}[1],2\texttt{b}[1]\right\}=\gcd\left\{\texttt{b}[2],2\texttt{b}[1]\right\},
$$
we conclude that $l=1$ or $l=2$.
Therefore, since $\gcd\left\{\texttt{b}[1],\texttt{b}[2]\right\}=1$,  the only possibilities for $\texttt{b}[1]$ and $\texttt{b}[2]$ are
$$
(\texttt{b}[1],\texttt{b}[2])=(1,3), \quad  (\texttt{b}[1],\texttt{b}[2])=(1,4), \quad \text{and} \quad (\texttt{b}[1],\texttt{b}[2])=(2,5). 
$$

If $(\texttt{b}[1],\texttt{b}[2])=(2,5)$ we have a compact connected $6$-dimensional isotropy submanifold fixed by $\Z_2$ with an effective $S ^1\cong S^1/ \Z_2$-Hamiltonian action, with only $4$ fixed points, 
($P_0, P_1,P_2$ and $P_4$) and weights 
$$
\{1,2,3\},\{-\texttt{b}[3]/2,\texttt{b}[3]/2,2\}, \{-1,-2,1\},\{-1,-2,-3\}.
$$
Then, by the classification of Hamiltonian circle actions with a minimal number of fixed points on a $6$-dimensional manifold, we have that $\texttt{b}[3]=2$. Using  the ABBV Localization formula on $\M$ with $\mu=1$,
we obtain $\texttt{b}[4]=1$ and then the resulting set of weights can be obtained from the positive multigraph of Case $4$ in Figure~\ref{graphs8} (multigraph $\#$ $75$).

If $(\texttt{b}[1],\texttt{b}[2])=(1,4),$ we get, by the same methods, the set of weights for the reversed circle action of the case  $(\texttt{b}[1],\texttt{b}[2])=(2,5)$ described above and so it can again be obtained from the  multigraph of Case $4$.

If  $(\texttt{b}[1],\texttt{b}[2])=(1,3)$, we obtain the following multiset of weights:
$$
\{ \{1,2,3,4 \}, \{ 2,3,-\texttt{b}[3],\texttt{b}[3] \}, \{ -2,-1,1,2 \}, \{ -2,-3,-\texttt{b}[4],\texttt{b}[4] \}, \{-4,-3,-2,-1\}\}.
$$
Using the ABBV Localization formula, with $\mu=1$, we have
$$
0=\int_\M 1= \sum_{i=0}^4 \frac{1}{\prod_{j=1}^4 w_{ij}}=\frac{1}{24} - \frac{1}{6\texttt{b}[3] ^2}+\frac{1}{4}  - \frac{1}{6\texttt{b}[4] ^2} + \frac{1}{24},
$$
and so $\texttt{b}[3] =\texttt{b}[4]=1$. Consequently, this set of weights falls again into Case $4$ and can be obtained with multigraph  $\#$ $75$.

\begin{figure}[h!]
\begin{center}
\includegraphics[scale=.8]{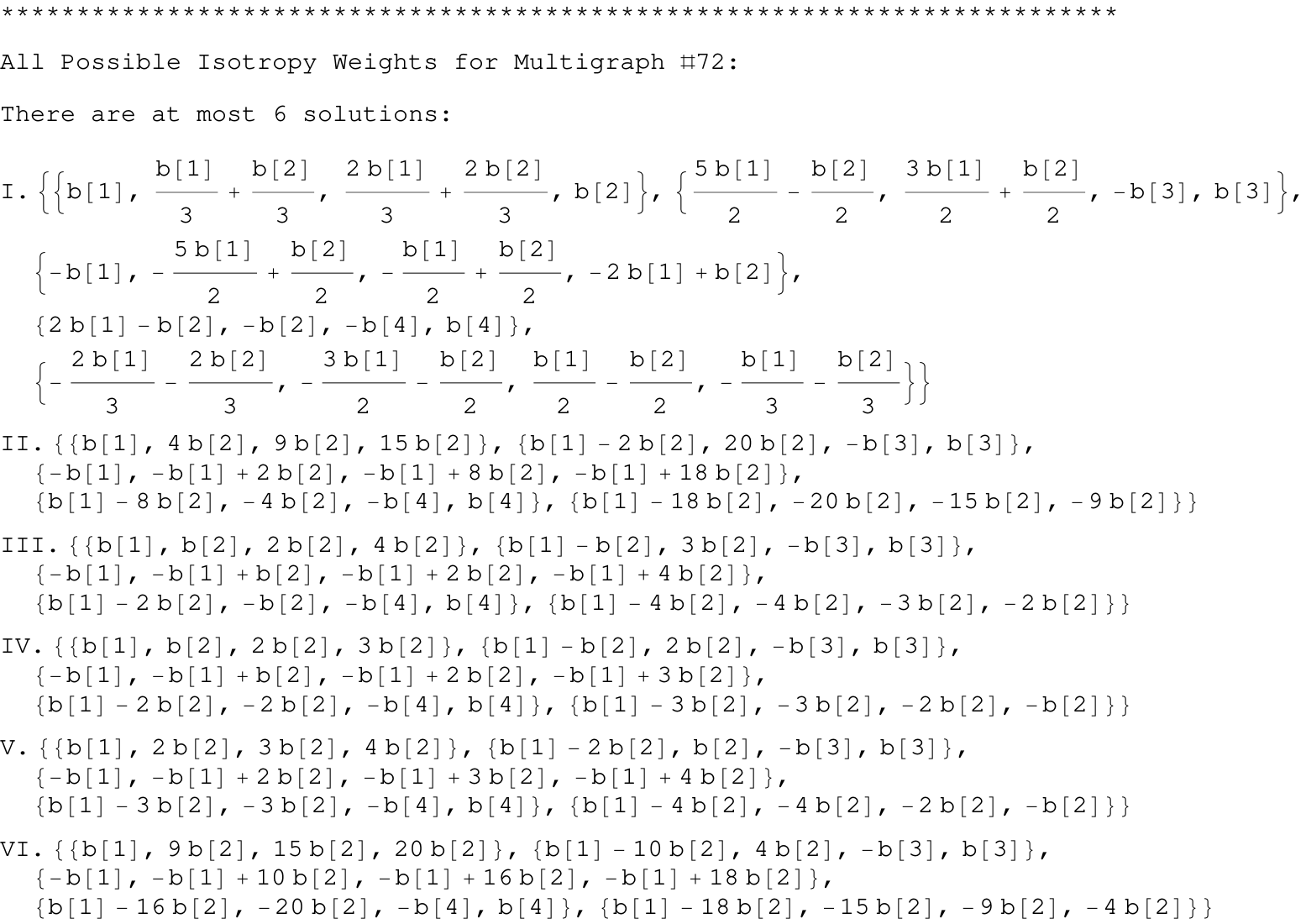}
\caption{Possible isotropy weights for Case $3$.}
\label{wc83}
\end{center}
\end{figure}

\textbf{Case $3$:} For the multigraph in Case $3$ of Figure~\ref{graphs8} the weights given by our algorithm are the ones in Figure~\ref{wc83}. Note that this multigraph has a pair of multiple  edges from $P_0$ to $P_4$. 

{\bf I.}  For the set of weights in Figure~\ref{wc83} I. we see that the weights corresponding to the two edges from $P_0$ to $P_4$ are $\frac{\texttt{b}[1] +\texttt{b}[2] }{3}$ and  $2(\frac{\texttt{b}[1] +\texttt{b}[2] }{3})$ and that the other two weights at $P_0$ are $\texttt{b}[1]$ and $\texttt{b}[2]$. Then Lemma~\ref{lemma:8.2} $(2)$ implies that we must have
$$
\gcd\{\ell,\texttt{b}[1]\} = \ell \quad \text{or} \quad \gcd\{\ell,\texttt{b}[2]\} = \ell, 
$$
with $\ell=\frac{\texttt{b}[1] +\texttt{b}[2] }{3}$. Hence we must have
$$
\texttt{b}[1] =2\texttt{b}[2]  \quad \text{and} \quad  \ell= \texttt{b}[2],
$$
or
$$
\texttt{b}[2] =2\texttt{b}[1] \quad \text{and} \quad \ell= \texttt{b}[1].
$$
Since the action is effective we conclude that, in the first case, we have $\texttt{b}[2]=1$ and $\texttt{b}[1]=2$, while, in the last, we have $\texttt{b}[1]=1$ and $\texttt{b}[2]=2$. However, both cases are impossible since the weights at $P_1$ would not be integers. 

{\bf II.}  For the set of weights in Figure~\ref{wc83} II. we see that the weights corresponding to the two edges from $P_0$ to $P_4$ are  $9\texttt{b}[2]$ and $15\texttt{b}[2]$ and that the other two weights at $P_0$ are $\texttt{b}[1]$ and $4\texttt{b}[2]$.
Since the action is effective, we have $\gcd\left\{\texttt{b}[1],\texttt{b}[2]\right\}=1$. Moreover, by Lemma~\ref{lemma:8.2}  $(2)$, we have
$$
\gcd\left\{3\texttt{b}[2],\texttt{b}[1]\right\}=3\texttt{b}[2],
$$
and so
$$
\texttt{b}[1]=3 \quad \text{and} \quad  \texttt{b}[2]=1.
$$
In this case the (connected) set of points fixed by $\Z_5$ would contain $P_0$ and $P_4$. However, the subgroup $\Z_5$ acts trivially on a $1$-dimensional complex subspace of $T_{P_0}\M$ and on a $3$-dimensional complex subspace at $T_{P_4}\M$ (since the  weights at $P_4$ and $P_0$ would respectively be $\{-15,-20,-15,-9\}$ and $\{3,4,9,15\}$), which is impossible.
 
{\bf III.}  For the set of weights in Figure~\ref{wc83} III. we see that the weights corresponding to the two edges from $P_0$ to $P_4$ are  $2\texttt{b}[2]$ and $4\texttt{b}[2]$ and that the other two weights at $P_0$ are $\texttt{b}[1]$ and $\texttt{b}[2]$.
Since the action is effective, we have $\gcd\left\{\texttt{b}[1],\texttt{b}[2]\right\}=1$. Moreover, by Lemma~\ref{lemma:8.2}  $(2)$, we have
$$
\gcd\left\{2\texttt{b}[2],\texttt{b}[1]\right\}=2\texttt{b}[2],
$$
and so
$$
\texttt{b}[1]=2 \quad \text{and} \quad  \texttt{b}[2]=1.
$$
However, in this case the weight $\texttt{b}[1]-2\texttt{b}[2]$ at $P_3$ would be $0$ which is impossible.

{\bf IV.}  For the set of weights in Figure~\ref{wc83} IV. we see that the weights corresponding to the two edges from $P_0$ to $P_4$ are  $\texttt{b}[2]$ and $3\texttt{b}[2]$ and that the other two weights at $P_0$ are $\texttt{b}[1]$ and $2\texttt{b}[2]$.
Since the action is effective, we have $\gcd\left\{\texttt{b}[1],\texttt{b}[2]\right\}=1$. If $\texttt{b}[2]=1$, then the weights $\texttt{b}[1]-\texttt{b}[2]$ and $2\texttt{b}[2]-\texttt{b}[1]$ respectively at $P_1$ and  $P_2$ cannot be simultaneously positive and so we conclude that  $\texttt{b}[2]> 1$. Then, the (connected) set of points fixed by $\Z_{\texttt{b}[2]}$ is a $6$-dimensional manifold  with an effective $S ^1\cong S^1/ \Z_{\texttt{b}[2]}$-Hamiltonian action with  $4$ fixed points (all except $P_2$) and weights 
$$
\{1,2,3\},\{-\texttt{b}[3]/\texttt{b}[2],\texttt{b}[3]/\texttt{b}[2],2\}, \{-\texttt{b}[4]/\texttt{b}[2],-2,\texttt{b}[4]/\texttt{b}[2]\},\{-1,-2,-3\}.
$$
Then, by the classification of Hamiltonian circle actions with a minimal number of fixed points on a $6$-dimensional manifold, we have that $\texttt{b}[3]=\texttt{b}[4]=\texttt{b}[2]$ and then the resulting set of weights can be obtained from the positive multigraph of Case $4$   (multigraph $\#$ $75$).

{\bf V.}  For the set of weights in Figure~\ref{wc83} V. we see that the weights corresponding to the two edges from $P_0$ to $P_4$ are  $2\texttt{b}[2]$ and $4\texttt{b}[2]$ and that the other two weights at $P_0$ are $\texttt{b}[1]$ and $3\texttt{b}[2]$. Then, just like in III., we have  
$$\texttt{b}[1]=2 \quad \text{and} \quad  \texttt{b}[2]=1.$$
However, in this case, the weight $\texttt{b}[1]-2\texttt{b}[2]$ at $P_2$ would be $0$ which is impossible.

{\bf VI.}  For the set of weights in Figure~\ref{wc83} VI. we see that the weights corresponding to the two edges from $P_0$ to $P_4$ are  $9\texttt{b}[2]$ and $15\texttt{b}[2]$ and that the other two weights at $P_0$ are $\texttt{b}[1]$ and $20\texttt{b}[2]$. Then, just like in II., we have 
 $$
\texttt{b}[1]=3 \quad \text{and} \quad  \texttt{b}[2]=1.
$$
However, in this case, the weight $\texttt{b}[1]-10\texttt{b}[2]$ at $P_1$ would be negative which is impossible.
 
\begin{figure}[h!]
\begin{center}
\includegraphics[scale=.8]{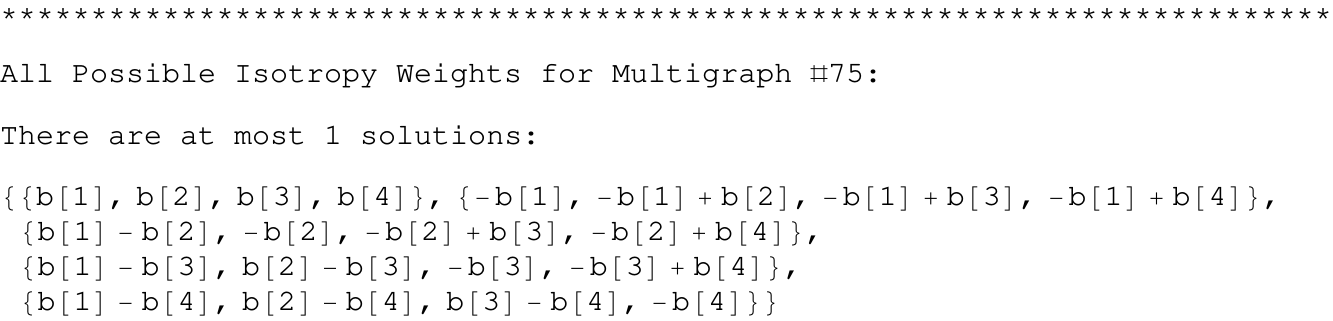}
\caption{Possible isotropy weights for Case $4$.}
\label{wc84}
\end{center}
\end{figure}

\textbf{Case $4$:}  For the multigraph in Case $4$ of Figure~\ref{graphs8} the weights given by our algorithm are the ones in Figure~\ref{wc84}. These are precisely the weights of the
$S^1$-action on $\C P^4$ described in Example~\ref{cpn}, with \texttt{b}$[\,i\,]=\xi_0-\xi_i$, $i=1,2,3,4$.

Since all the possible isotropy weights obtained can be  included in this last case and
the basis for $H_{S^1}^*(\M;\Z)$ described in Theorem~\ref{bases} and the (equivariant) Chern
classes only depend on the isotropy weights at the fixed points, we can summarize our results in the following theorem. 
\begin{thm}\label{thm dim8}
Let $(\M,\omega)$ be a compact symplectic manifold of dimension $8$, with a Hamiltonian $S^1$-action
with $5$ fixed points. 
If there exists a non-negative multigraph associated to the action, then
the isotropy weights at the fixed points agree with the one of $\C P^4$
for the standard $S^1$-action. Moreover,
 the cohomology ring and Chern classes
agree with the ones of $\C P^4$, i.e.
$$
H^*(\M;\Z)=\Z[y]/(y^5)\quad\mbox{and}\quad \cc=(1+y)^5\;,
$$
where $y$ has degree 2.
\end{thm}
Corollary~\ref{main dim 8} in the introduction is then an easy consequence
of Proposition~\ref{t2}, Proposition~\ref{not pm1} and Theorem~\ref{thm dim8}.

We can now
prove Theorem~\ref{RS1} stated in the Introduction.
\begin{proof}[Proof of Theorem~\ref{RS1}]
By Lemma~\ref{w=tau} it is not restrictive to assume that $[\omega]=y$, and by Theorem~\ref{m not 2}, we know that  $C_1$ is either 1 or 5.\\
(i)$\implies$(ii), (iii), (iv), (v): this is exactly the content of Theorem~\ref{thm dim8}. \\ (ii)$\iff$(iii) and (ii)$\iff$(iv): this follows
from Theorem~\ref{m not 2}. \\
(v)$\implies$(i): it follows easily from the definitions (see Example~\ref{cpn}).\\
(ii)$\implies$(v): this is directly implied by Theorem~\ref{hattori2}. Indeed, the Euler characteristic of $\M$ satisfies $\chi_{-1}(\M)=\int_\M \cc_4=5$ (see
\eqref{cn}). Moreover, by Lemma~\ref{w=tau}, we can choose $\omega-\psi\otimes x$ to be
primitive and positive, and, by
 Proposition~\ref{symp hattori}, the associated pre-quantization line bundle is quasi-ample and satisfies \eqref{D}. Thus Theorem~\ref{hattori2}
 implies that the isotropy weights at the fixed points agree with the ones of $\C P^4$ with the standard $S^1$-action. \\
 The remaining equivalences follow easily from the above ones.
\end{proof}

\end{document}